\documentclass{gtart_h}

\def\ifplaintex{\expandafter\ifx\csname documentclass\endcsname\relax}

\def\ifplaintex{\expandafter\ifx\csname documentclass\endcsname\relax}

\ifplaintex 
\else
\fi

\expandafter\ifx\csname epsfbox\endcsname\relax\input epsf\fi

\def\gt{{\mathsurround=0pt\it $\cal G\mskip-2mu$eometry \&\ 
$\cal T\!\!$opology}}        

\def\gtp{{\mathsurround=0pt\it $\cal G\mskip-2mu$eometry \&\ 
$\cal T\!\!$opology $\cal P\!$ublications}}  

\def\lognumber#1{\def\thelognumber{#1}}
\def\volumenumber#1{\def\thevolumenumber{#1}}
\def\papernumber#1{\def\thepapernumber{#1}}
\def\volumeyear#1{\def\thevolumeyear{#1}}

\def\pagenumbers#1#2{\def\startpage{#1}\def\finishpage{#2}}
\def\published#1{\def\publishdate{#1}}
\def\proposed#1{\def\theproposer{#1}}
\def\seconded#1{\def\theseconders{#1}}
\def\received#1{\def\receiveddate{#1}}
\def\revised#1{\def\reviseddate{#1}}
\def\accepted#1{\def\accepteddate{#1}}

\def\asciiaddress#1{\def\theasciiaddress{#1}}

\long\def\asciiabstract#1{\long\def\theasciiabstract{#1}}
\def\asciikeywords#1{\def\theasciikeywords{#1}}

\let\\\par\let\thelognumber\relax
\let\thevolumenumber\relax\let\thepapernumber\relax
\let\thevolumeyear\relax\let\thesamplenumber\relax\let\startpage\relax
\let\finishpage\relax\let\publishdate\relax\let\receiveddate\relax
\let\reviseddate\relax\let\accepteddate\relax\let\theasciititle\relax
\let\theasciiauthors\relax\let\theasciiaddress\relax
\let\theasciiabstract\relax\let\theasciikeywords\relax
\let\theasciiemail\relax\let\theshortauthors\relax\let\theshorttitle\relax

\long\def\maketitlep{   

\count0=\startpage

\gt\hfill      
\hbox to 77pt{\vbox to 0pt{\vglue -15pt\epsfbox{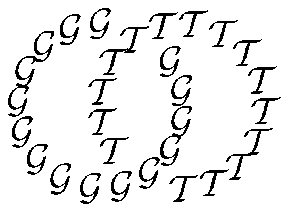}\vss}\hss}
\break
{\small\ifx\thesamplenumber\relax 
Volume \else Sample
\fi\thevolumenumber\ (\thevolumeyear)
\startpage--\finishpage\nl
Published: \publishdate}
\vglue 0.5truein plus 0.4fil minus 0.1truein

{\parskip=0pt\leftskip 0pt plus 1fil\def\\{\par\smallskip}{\ifplaintex\large
\else\Large\fi\bf\thetitle}\par\medskip}   

\vglue 0pt plus 0.1fil 

{\parskip=0pt\leftskip 0pt plus 1fil\def\\{\par}{\sc\theauthors}
\par\medskip}

\vglue 0pt plus 0.1fil 

{\small\parskip=0pt\let\newline\\
{\leftskip 0pt plus 1fil\def\\{\par}{\sl\theaddress}\par}
\expandafter\ifx\theemail\relax    
\relax\else\vglue 5pt plus 0.02fil minus 2pt\def\\{\stdspace{\rm 
and}\stdspace} 
\cl{Email:\stdspace\tt\theemail}\fi
\ifx\theurl\relax                  
\relax\else\vglue 5pt plus 0.02fil minus 2pt\def\\{\stdspace{\rm 
and}\stdspace}
\cl{URL:\stdspace\tt\theurl}\fi\par}

\vglue 7pt plus 0.3fil minus 3pt

{\bf Abstract}
\vglue 5pt plus 0.1fil minus 2pt

\theabstract

\vglue 7pt plus 0.3fil minus 3pt

{\bf AMS Classification numbers}\quad Primary:\quad \theprimaryclass

Secondary:\quad \thesecondaryclass

\vglue 5pt plus 0.3fil minus 2pt

{\bf Keywords:}\quad \thekeywords

\vglue 10pt plus 0.5fil minus 5pt

{\small  Proposed: \theproposer\hfill Received: \receiveddate\nl
Seconded: \theseconders\hfill 
\ifx\reviseddate\relax                         
Accepted: \accepteddate                        
\else
Revised: \reviseddate                          
\fi}
\eject
}       

\font\phead=cmsl9 scaled 950
\font=cmsl9 scaled 1050
\font\pnum=cmbx10 scaled 913
\font\lnum=cmbx10 
\font\pfoot=cmsl9 scaled 950
\font=cmsl9 scaled 1050
\ifplaintex
\headline{\vbox to 0pt{\vskip -4.5mm\line{\small\phead\ifnum
\count0=\startpage ISSN 1364-0380 (on line)
1465-3060 (printed) \hfill {\pnum\folio}\else\ifodd\count0\def\\{ }%
\ifx\theshorttitle\relax\thetitle\else\theshorttitle\fi\hfill{\pnum\folio}
\else\def\\{ and }{\pnum\folio}\hfill\ifx\theshortauthors\relax\theauthors
\else\theshortauthors\fi\fi\fi}\vss}}
\footline{\vbox to 0pt{\vglue 0mm\line{\small\pfoot\ifnum\count0=\startpage
\copyright\ \gtp\hfill\else
\gt, Volume \thevolumenumber\ (\thevolumeyear)\hfill\fi}\vss
}}
\else
\makeatletter
\def\@oddhead{{\small\lhead\ifnum\count0=\startpage ISSN 1364-0380 (on line)
1465-3060 (printed) \hfill {\lnum\number\count0}\else\ifodd\count0
\def\\{ }\ifx\theshorttitle\relax \thetitle \else\theshorttitle\fi\hfill
{\lnum\number\count0}\else\def\\{ and }{\lnum\number\count0}
\hfill\ifx\theshortauthors\relax 
\theauthors\else\theshortauthors\fi\fi\fi}}\def\@evenhead{\@oddhead}
\def\@oddfoot{\small\lfoot\ifnum\count0=\startpage\copyright\ \gtp\hfill\else
\gt, Volume \thevolumenumber\ (\thevolumeyear)\hfill\fi}
\def\@evenfoot{\@oddfoot}
\makeatother
\fi

\newwrite\gtoutfile
\long\gdef\makeheadfile{  
{\def\\{, }\def\s{ }
\immediate\openout\gtoutfile head.xxx
\immediate\write\gtoutfile{Proxy-for: \ifx\theasciiauthors\relax
\theauthors\else\theasciiauthors\fi\s<\ifx\theasciiemail\relax\theemail\else\theasciiemail\fi>}
\immediate\write\gtoutfile{\noexpand\\}
\immediate\write\gtoutfile{Authors: \ifx\theasciiauthors\relax
\theauthors\else\theasciiauthors\fi}
\immediate\write\gtoutfile{Title: \ifx\theasciititle\relax
\thetitle\else\theasciititle\fi}
\immediate\write\gtoutfile{Subj-class: GT or SG or MG etc}
\immediate\write\gtoutfile{MSC-class: \theprimaryclass\ifx\thesecondaryclass\relax\else, \thesecondaryclass\fi}
\immediate\write\gtoutfile{Journal-ref: Geom. Topol. \thevolumenumber
(\thevolumeyear) \startpage-\finishpage}
\immediate\write\gtoutfile{Comments: Published by Geometry and Topology at}
\immediate\write\gtoutfile{\s\s http://www.maths.warwick.ac.uk/gt/GTVol\thevolumenumber/paper\thepapernumber.abs.html}
\immediate\write\gtoutfile{\noexpand\\}
\immediate\write\gtoutfile{}
\ifx\theasciiabstract\relax
\immediate\write\gtoutfile{\theabstract}\else
\immediate\write\gtoutfile{\theasciiabstract}\fi
\immediate\write\gtoutfile{}
\immediate\write\gtoutfile{\noexpand\\}
\immediate\write\gtoutfile{}
\immediate\closeout\gtoutfile}}  

\def\maketitlepage{\maketitlep\makeheadfile}
\let\maketitle\maketitlepage

\lognumber{339}
\volumenumber{8}\papernumber{9}\volumeyear{2004}
\pagenumbers{335}{412}
\received{12 July 2003}
\revised{12 February 2004}
\published{14 February 2004}
\accepted{30 January 2004}

\proposed{Bill Dwyer}
\seconded{Thomas Goodwillie, Haynes Miller}

\usepackage{amsmath,amssymb}
\usepackage[matrix,arrow]{xy}

\numberwithin{equation}{section}
\newtheorem{theorem}[equation]{Theorem}
\newtheorem{lemma}[equation]{Lemma}

\newtheorem{corollary}[equation]{Corollary}
\theoremstyle{definition}

\newtheorem{construction}[equation]{Construction}
\newtheorem{remark}[equation]{Remark}
\newtheorem{example}[equation]{Example}
\swapnumbers
\newtheorem{numbered paragraph}[equation]{}

\DeclareMathOperator{\Ch}{Ch}
\DeclareMathOperator{\coCh}{coCh}
\DeclareMathOperator{\Hom}{Hom}

\DeclareMathOperator{\Ext}{Ext}
\DeclareMathOperator{\hExt}{{\mathbb E}xt}
\DeclareMathOperator{\map}{map}
\DeclareMathOperator{\Rg}{Ring}
\DeclareMathOperator{\der}{der}
\DeclareMathOperator{\Aut}{Aut}
\DeclareMathOperator{\End}{End}
\DeclareMathOperator{\Ho}{Ho}
\DeclareMathOperator{\holim}{holim}
\DeclareMathOperator{\colim}{colim}
\DeclareMathOperator{\f}{f}

\newcommand{\mF}{{\mathbb F}}
\newcommand{\mQ}{{\mathbb Q}}
\newcommand{\mS}{{\mathbb S}}
\newcommand{\mZ}{{\mathbb Z}}
\newcommand{\Ab}{{\mathcal A}b}
\newcommand{\B}{{\mathcal B}}
\newcommand{\C}{{\mathcal C}}
\newcommand{\D}{{\mathcal D}}
\newcommand{\F}{{\mathcal F}}
\newcommand{\I}{\mathcal I}
\newcommand{\Lc}{\mathcal L}
\newcommand{\GR}{\mathcal GR}
\newcommand{\Z}{\mathcal Z}
\newcommand{\iso}{\cong}
\newcommand{\rscript}[1]{^{\mbox{\scriptsize {#1}}}}
\newcommand{\sm}{\wedge}
\newcommand{\tensor}{\otimes}
\newcommand{\Wedge}{{\scriptstyle \vee}}

\renewcommand{\to}{\longrightarrow}
\newcommand{\varl}[2]{\stackrel{#2}{\hbox to #1{\leftarrowfill}}}
\newcommand{\varr}[2]{\stackrel{#2}{\hbox to #1{\rightarrowfill}}}

\begin{document}

\title[Formal groups and stable homotopy of commutative rings]
{Formal groups and stable homotopy\\of commutative rings}
\author{Stefan Schwede}
\address{Mathematisches Institut, Universit\"at Bonn\\53115 Bonn, Germany}
\asciiaddress{Mathematisches Institut, Universitaet Bonn\\53115 Bonn, Germany}
\email{schwede@math.uni-bonn.de}

\begin{abstract}
We explain a new relationship between formal group laws and ring spectra 
in stable homotopy theory. We study a ring spectrum denoted 
$DB$ which depends on a commutative ring $B$ and
is closely related to the topological Andr\'e--Quillen homology of $B$. 
We present an explicit construction which to every 1--dimensional 
and commutative formal group law $F$ over $B$ associates a morphism 
of ring spectra $F_*\co H\mZ\to DB$ from the
Eilenberg--MacLane ring spectrum of the integers.
We show that formal group laws account for all such ring spectrum maps, 
and we identify the space of ring spectrum maps between $H\mZ$ and $DB$. 
That description involves formal group law data and
the homotopy units of the ring spectrum $DB$.
\end{abstract}
\asciiabstract{%
We explain a new relationship between formal group laws and ring
spectra in stable homotopy theory. We study a ring spectrum denoted DB
which depends on a commutative ring B and is closely related to the
topological Andre-Quillen homology of B.  We present an explicit
construction which to every 1-dimensional and commutative formal group
law F over B associates a morphism of ring spectra F_*: HZ --> DB from
the Eilenberg-MacLane ring spectrum of the integers.  We show that
formal group laws account for all such ring spectrum maps, and we
identify the space of ring spectrum maps between HZ and DB.  That
description involves formal group law data and the homotopy units of
the ring spectrum DB.}

\asciikeywords{Ring spectrum, formal group law, Andre-Quillen homology}
\maketitlepage

\section{Introduction}

In this paper we explain a new relationship between formal group laws 
and ring spectra in stable homotopy theory. 
We use formal group laws to construct  maps of ring spectra 
and describe spaces of ring spectrum maps in terms of formal group data.

Our main object of study is a ring spectrum denoted $DB$ 
which functorially depends on a commutative ring $B$. 
The ring spectrum $DB$ controls the stable homotopy theory 
of augmented commutative $B$--algebras. Its homotopy groups $\pi_{\ast} DB$ 
are the Cartan--Bousfield--Dwyer algebra of stable homotopy operations 
of commutative simplicial $B$--algebras 
\cite{Cartan, Bous:operations, Dwyer:oper}. 
The homotopy groups $\pi_{\ast} DB$ are also isomorphic to the 
{\em $\Gamma$--homology}~\cite{Robinson-Whitehouse:GammaHomology}, 
relative to $B$, of the polynomial algebra $B[x]$ on one generator,
and to the 
{\em topological Andr\'e--Quillen homology}~\cite{Basterra:AQHomology}
of the associated Eilenberg--MacLane spectra. Both $\Gamma$--homology
and topological Andr\'e--Quillen homology arise in obstruction theory 
for $E_{\infty}$ ring spectrum structures~\cite{Robinson:Gamma-Einfty}.
We elaborate more on the precise relationship in Remark~\ref{DB and TAQ} below,
and we also refer to the survey paper~\cite{Basterra-Richter} 
by Basterra and Richter.

We present an explicit construction 
which to every 1--dimensional and commutative formal group law $F$ over $B$ 
associates a homomorphism of ring spectra 
\[ F_{\ast} \,\co\, H\mZ \ \to \ DB \]
from the Eilenberg--MacLane ring spectrum of the integers to $DB$. 
We prove that in this way formal group laws account 
for all ring spectrum maps, ie, we show:

\medskip

{\bf Theorem}\qua {\sl The construction which sends a 1--dimensional,
commutative formal group law $F$ to the ring spectrum map $F_{\ast}$
induces a natural bijection between the set of strict isomorphism
classes of formal group laws over $B$ and the set of homotopy classes
of ring spectrum maps from $H\mZ$ to $DB$.}

\medskip

This theorem is a corollary of the identification of 
the space of all ring spectrum maps between the Eilenberg--MacLane 
ring spectrum $H\mZ$ and $DB$. That description involves the 
{\em homotopy units} $(DB)^{\times}$ of the ring spectrum $DB$. 
The group-like simplicial monoid $(DB)^{\times}$ of homotopy units 
acts by conjugation on the space Ring$(H\mZ,DB)$ of ring spectrum maps. 
The construction $F_{\ast}$ of the first theorem extends to a natural weak map
from the classifying space of the groupoid of formal group laws 
and strict isomorphisms $\mathcal{FGL}\rscript{str}(B)$ 
to the homotopy orbit space 
$\Rg(H\mZ,DB)/\mbox{conj.}$ of the conjugation action 
of the identity component $(DB)^{\times}_1$ of the units 
on the space of ring spectrum maps. In Theorem \ref{main} we show

\medskip

{\bf Theorem}\qua {\sl The weak map 
\[ \mathcal{FGL}\rscript{\em str}(B) \ \to \ \Rg(H\mZ,DB)/\mbox{\em conj.} \]
from the classifying space of the groupoid of formal group laws 
and strict isomorphisms to the homotopy orbit space of the conjugation action 
is a weak homotopy equivalence.}

\medskip

Another corollary of this main theorem is that the ring spectrum $DB$ 
is not equivalent to the Eilenberg--MacLane ring spectrum of 
any differential graded algebra (or, equivalently, of any simplicial ring);
see Corollary \ref{DB is not simplicial} for the precise statement. 
This should be compared to the fact that the spectrum underlying $DB$ 
(ie, ignoring the multiplication)
is stably equivalent to the smash product $H\mZ\sm HB$ 
(see Theorem \ref{DB summary} (b)), in particular it is equivalent
to a product of Eilenberg--MacLane spectra. 
Corollary \ref{DB is not simplicial} says that the multiplicative structure 
of $DB$ is considerably more complicated.

\medskip
\noindent
{\bf Prerequisites}\qua We freely use the language and standard results
from the theory of model categories; the original source for this material is  
\cite{Quillen:HA}, a more modern introduction can be found in 
\cite{Dwyer-Spalinski}, and the ultimate reference 
is currently \cite{Hovey:book}.
Our notion of ring spectrum is that of a {\em Gamma-ring} 
(see \cite[2.13]{Lydakis:Smash} or \cite[Def.\ 2.1]{Sch:Gamma}). 
Gamma-rings are based on a symmetric monoidal smash product for 
$\Gamma$--spaces with good homotopical properties 
\cite{Segal:categories, BF, Lydakis:Smash}. 
The foundational material about the homotopy theory of Gamma-rings 
and their modules can be found in \cite{Sch:Gamma}; a summary is also
given in Section \ref{general review}.
The results of this paper can be translated into other frameworks
for ring spectra by the general comparison procedures described in
\cite{MMSS, Sch:comparison}.
We also need a few basic facts from the theory of
formal group laws, which in this paper 
(with the exception of Section \ref{loose ends}) are always 
1--dimensional and commutative; 
all we need is contained in \cite{Lazard:formels} or in
Chapter III, \S 1 of \cite{Froehlich}.

\bigskip
\noindent
{\large\bf  Outline of the paper}

\medskip

{\bf Section \ref{general review}}\qua 
We review some general facts about $\Gamma$--spaces and Gamma-rings. 
We recall the assembly map (\ref{assembly map})
from the smash product to the composition product of $\Gamma$--spaces
which is used several times in this paper.

\medskip
\noindent
{\bf Section \ref{DB review}}
We review the Gamma-ring $DB$ and summarize some of its properties 
in Theorem \ref{DB summary}.
Construction  \ref{Gamma-ring map} associates 
a homomorphism of ring spectra $F_*\co H\mZ\to DB$ 
to every formal group law $F$.
We state in Theorem \ref{homotopy classes} how this accounts 
for all homotopy classes of Gamma-ring maps.

\medskip
\noindent
{\bf Section \ref{conjugation}}\qua 
For every Gamma-ring $R$ we construct a natural conjugation action 
of the simplicial monoid of homotopy units $R^{\times}$ (\ref{unit definition})
on $R$ through Gamma-ring homomorphisms.

\medskip
\noindent
{\bf Section \ref{comparison map}}\qua We construct a weak map 
\[ \xymatrix{ \mathcal{FGL}\rscript{str}(B) & 
\widetilde{\mathcal{FGL}}\rscript{str}(B)  \ar[l]_-{\sim}   
\ar[r]^-{\kappa} &  \Rg(H\mZ,DB)/\mbox{conj.} }\]
from the classifying space of the groupoid of formal group laws 
and strict isomorphisms to the homotopy orbit space of the conjugation action
of Section \ref{conjugation}. 
The left map is a weak equivalence by construction.
Theorem \ref{main}, which is the main theorem of this paper, 
says that the right map $\kappa$ is a weak equivalence.

\medskip
\noindent
{\bf Section \ref{filtration}}\qua We use a filtration of the Gamma-ring $DB$,
coming from powers of the augmentation ideal,
to reduce the proof of Theorem \ref{main} to showing that a truncated version 
\[  \kappa_k \,\co\,\widetilde{{\mathcal B}ud}^k(B) \ 
\to \   \Rg(H\mZ,D_kB)/\mbox{conj.} \]
of the map $\kappa$ of Section  \ref{comparison map}
is a weak equivalence for all $k\geq 1$ (see Theorem \ref{k-bud comparison}). 
Here $\widetilde{{\mathcal B}ud}^k(B)$ is (weakly equivalent to)
the classifying space  of the groupoid of $k$--buds (or $k$--jets) 
of formal group laws
and $D_kB$ is the ``quotient'' Gamma-ring of $DB$ by the ``ideal''
coming from $(k+1)$st powers of the augmentation ideal. 

\medskip
\noindent
{\bf Section \ref{singular extensions}}\qua 
We exploit that two successive stages in the filtration
of $DB$ form a ``singular extension'' of ring spectra
\[ (B\tensor {\mathcal S}^k)^! \ \to \ D_kB \ \to \ D_{k-1}B  \]
where $(B\tensor {\mathcal S}^k)^!$ is a ``square zero ideal'' 
coming from the $k$-th symmetric power functor. 
This allows us to reduce the problem  to showing that a certain map
\[ \kappa_{B\tensor {\mathcal S}^k} \,\co\,\widetilde{\Z}(B\tensor {\mathcal S}^k) \ \to \ 
\der(H\mZ,B\tensor {\mathcal S}^k)/\mbox{conj.} \]
to the homotopy orbit space of the derivations of $H\mZ$ with coefficients in 
the symmetric power functor is a weak equivalence.
Here $\widetilde{\Z}(B\tensor {\mathcal S}^k)$ is (weakly equivalent to) 
the classifying space of the groupoid 
of symmetric 2--cocycles of degree $k$ over $B$.

\medskip
\noindent
{\bf Section \ref{2-cocycles}}\qua 
For later use we define a map 
\[ \lambda_G \,\co\,\Z(G) \ \to \ 
\map_{\mathcal{GR}}(H\mZ\rscript{c},H\mZ\times G_{st}^!)_{hG_{st}(\mZ)}\] 
for every functor $G$ from the category of finitely generated free abelian 
groups to the category of all abelian groups; the important case is
when $G$ is a symmetric power functor $B\tensor {\mathcal S}^k$.
Here $\Z(G)$ is the classifying space of the groupoid 
of symmetric 2--cocycles of the functor $G$ (\ref{def-cocycles}),
$\map_{\mathcal{GR}}$ denotes the simplicial mapping space of Gamma-rings,
$G_{st}$ is the Dold--Puppe stabilization of $G$ (\ref{DP and Q}) 
and the Gamma-ring $H\mZ\times G_{st}^!$ is the split singular extension
of $H\mZ$ by the bimodule $G_{st}^!$ \eqref{split extension}.

\medskip
\noindent
{\bf Section \ref{simplification}}\qua 
We compare the map  $\kappa_{B\tensor {\mathcal S}^k}$ 
of Section \ref{singular extensions}
with the map $\lambda_{B\tensor {\mathcal S}^k}$ of Section \ref{2-cocycles}
by means of a commutative square
\[(\ref{comparison square}) \hspace*{1cm}
\xymatrix@C=20mm{ \widetilde{\Z}(B\tensor {\mathcal S}^k)\ar_{\sim}[d] 
\ar^-{\kappa_{B\tensor {\mathcal S}^k}}[r] & 
\der(H\mZ,B\tensor {\mathcal S}^k)/\mbox{conj.} \ar^{\sim}[d] \\
\Z(B\tensor {\mathcal S}^k) \ar_-{\lambda_{B\tensor {\mathcal S}^k}}[r] & 
\map_{\mathcal{GR}}(H\mZ\rscript{c},H\mZ\times (B\tensor {\mathcal S}^k)_{st}^!)_{h(B\tensor {\mathcal S}^k)_{st}(\mZ)}
}\hspace*{1cm}\]
in which the vertical maps are weak equivalences.
Hence instead of showing that $\kappa_{B\tensor {\mathcal S}^k}$ is a weak equivalence
we may show that $\lambda_{B\tensor {\mathcal S}^k}$ is.

\medskip
\noindent
{\bf Section \ref{universal equivalence}}\qua 
By the results of the previous sections, the proof of the main theorem
is reduced to an identification of the space of Gamma-ring maps
\[ \map_{\mathcal{GR}}(H\mZ\rscript{c},
H\mZ\times (B\tensor {\mathcal S}^k)_{st}^!) \]
(or more precisely a certain homotopy orbit space thereof)
with the classifying space of symmetric 2--cocycles.
In this section we reinterpret the above mapping space in terms of
the category $s\F$ of simplicial functors from the category of 
finitely generated free abelian groups to the category of abelian groups.
We note that the construction which sends $G\in s\F$ to 
the split extension $H\mZ\times G^!$ (\ref{split extension}) has a left adjoint
\[ \Lc \,\co\,\GR/H\mZ \ \to \ \ s\F \] 
from the category of Gamma-rings over $H\mZ$ to the category 
of simplicial functors. Moreover, the two functors
form a Quillen adjoint pair between model categories. In order to identify 
the above mapping space, we evaluate the left adjoint $\Lc$
on the Gamma-ring $H\mZ\rscript{c}$, 
the cofibrant replacement of $H\mZ$. We denote by $J$ the functor 
which supports the universal symmetric 2--cocycle (\ref{universal cocycle}).
The main result of this section,
Theorem \ref{thm-adjoint of universal}, states that the map
\[ \Lc(H\mZ\rscript{c}) \ \to \ J \] 
which is adjoint to the ``universal derivation'' (\ref{con-derivation})
\[  (1,d_u) \,\co\,H\mZ\rscript{c} \ \to \ H\mZ\times J^! \] 
is a stable equivalence of simplicial functors.
This implies that for any reduced functor $G$ the homotopy groups of the space
\[ \map_{\mathcal{GR}}(H\mZ\rscript{c},H\mZ\times G_{st}^!) 
\ \iso \  \map_{s\F}(\Lc(H\mZ\rscript{c}),G_{st}) \]
are isomorphic to the hyper-cohomology groups $\hExt_{\F}^*(J,QG)$,
for $*\leq 0$,
in the abelian category $\F$ of reduced functors from the category of 
finitely generated abelian groups to the category of abelian groups. 
The chain complex $QG$ is MacLane's cubical construction for the functor $G$
(\ref{DP and Q}).

\medskip
\noindent
{\bf Section \ref{der and hyper}}\qua 
We prove a homological criterion,
Theorem \ref{homological criterion}, in terms of the functor $G$ 
for when the map $\lambda_G$ defined in Section \ref{2-cocycles}
is a weak equivalence. 
Loosely speaking, the criterion requires that ``MacLane cohomology equals
topological Hochschild cohomology'' for the functor $G$,
compare Remark \ref{MacLane vs TH}.
The precise meaning of this is that for 
all integers $m\leq 2$ the map
\[ \Ext^m_{\F}(I,G) \ \to \  \hExt^m_{\F}(I,QG) \]
be an isomorphism. 
The map $G\to QG$ is a model for Dold--Puppe stabilization
and it is initial, in the derived category of $\F$, among maps
from $G$ to complexes whose homology functors are additive. 
$\hExt^m_{\F}(I,-)$ denotes hyper-Ext groups of the functor $I$
with coefficients in a chain complex of functors. 
In Example \ref{counterexample} we also give a functor 
for which the criterion fails.

\medskip
\noindent
{\bf Section \ref{sec-Ext}}\qua 
In Theorem \ref{finite functors} we verify
the homological criterion of Theorem \ref{homological criterion} for the
symmetric power functors $G=B\tensor {\mathcal S}^k$. 
This finishes the proof of our main theorem:
Theorem \ref{homological criterion} shows that the map
$\lambda_{B\tensor {\mathcal S}^k}$ is a weak equivalence, hence by 
the commutative square (\ref{comparison square}) the map 
$\kappa_{B\tensor {\mathcal S}^k}$ is a weak equivalence.
The commutative diagram (\ref{singular fibre sequences}) of fibre sequences
shows inductively that the maps 
$\kappa_k\co \widetilde{{\mathcal B}ud}^k(B) \to \Rg(H\mZ,N_k)/\mbox{conj.}$ 
of Theorem \ref{k-bud comparison} are weak equivalences.
Hence the (weak) map 
$\kappa\co \mathcal{FGL}\rscript{str}(B) \to \Rg(H\mZ,DB)/\mbox{conj.}$ 
of Theorem \ref{main} is a weak equivalence.

\medskip
\noindent
{\bf Section \ref{loose ends}}\qua 
In the last section we give an application of the main theorem
as well as an outlook towards possible generalizations and future directions.
Variations of the main construction are possible, and some of them
are described in this final section. 
For example, formal groups 
can be replaced by formal $A$--modules where $A$ is any ring 
(not necessarily commutative). Such formal $A$--modules give rise 
to Gamma-ring maps from the Eilenberg--MacLane Gamma-ring of $A$ 
into $DB$. Furthermore when considering higher dimensional formal group laws, 
the natural target of the construction is a matrix Gamma-ring over $DB$ 
of the corresponding dimension.

\medskip
\noindent
{\bf Acknowledgments}\qua The author would like to thank 
Jeff Smith for suggesting the method used in Section \ref{conjugation} 
to obtain the conjugation action; moreover, many ideas that appear
in Section \ref{universal equivalence} arose in a
joint project with Smith on derivations of ring spectra.

\section{Review of $\Gamma$--spaces and Gamma-rings} \label{general review}

In this section we review some general facts about  
$\Gamma$--spaces and Gamma-rings, and we fix notation and terminology.
None of this material is new, but we present it in a form 
which is convenient for this paper.
We also prove certain properties of the assembly map 
(\ref{assembly map}) 
from the smash product to the composition product of $\Gamma$--spaces
which are used several times in this paper.

\begin{numbered paragraph} \label{basic definition}
{\bf $\Gamma$--spaces}\qua 
The category of $\Gamma$--spaces was introduced by Segal 
\cite{Segal:categories}, who showed that 
it has a homotopy category equivalent to the stable homotopy category of connective spectra. Bousfield and Friedlander \cite{BF} 
considered a bigger category of $\Gamma$--spaces in which the ones introduced 
by Segal appeared as the {\it special} $\Gamma$--spaces (\ref{special}). 
Their category admits a closed simplicial model category structure 
with a notion of stable weak equivalences giving rise again 
to the homotopy category of connective spectra. 
Then Lydakis \cite{Lydakis:Smash} introduced internal function objects 
and a symmetric monoidal smash product with good homotopical properties.

The category $\Gamma^{\text{op}}$  has one object 
$n^+=\{0,1,\dots,n\}$ for every non-negative integer $n$, 
and morphisms are the maps of sets which send 0 to 0. 
$\Gamma^{\text{op}}$ is equivalent to the opposite 
of Segal's category $\Gamma$ \cite{Segal:categories}. 
A $\Gamma$--space is a covariant functor from $\Gamma^{\text{op}}$ 
to the category of simplicial sets taking $0^+$ to a one point simplicial set.
A morphism of $\Gamma$--spaces is a natural transformation of functors. 
We denote by $\mathbb S$ the $\Gamma$--space which takes $n^+$ 
to $n^+$, considered as a constant simplicial set. 
If $X$ is a $\Gamma$--space and $K$ a pointed simplicial set, 
a new $\Gamma$--space $X\sm\, K$ is defined by setting 
$(X\sm\, K)(n^+) = X(n^+)\sm\, K$.

A $\Gamma$--space $X$ can be prolonged, by direct limit, 
to a functor from the category of finite pointed sets to the
category of (not necessarily finite) pointed sets. 
By degreewise evaluation and formation of the diagonal 
of the resulting bisimplicial sets, it can furthermore be promoted 
to a functor from the category of pointed simplicial sets to itself 
\cite[\S 4]{BF}. 
The extended functor preserves weak equivalences of simplicial sets 
\cite[Prop.\ 4.9]{BF} and is automatically simplicial, ie, 
it comes with coherent natural maps 
$K\, \sm\, X(L) \to X(K\, \sm\, L)$.
We will not distinguish notationally between the prolonged functor 
and the original $\Gamma$--space. 

The homotopy groups of a $\Gamma$--space $X$ are defined as
\[ \pi_n X \ = \ \mbox{colim}_i \, \pi_{n+i}|X(S^i)| \ , \]
where the colimit is formed using the maps
\[ S^1 \, \sm\, X(S^n) \ \to \  X(S^1\, \sm\, S^n)\ . \] 
A map of $\Gamma$--spaces  is a {\em stable equivalence} 
if it induces isomorphisms on homotopy groups.
Since the functor given by a prolonged $\Gamma$--space preserves connectivity
\cite[4.10]{BF}, the homotopy groups of a $\Gamma$--space are
 always trivial in negative dimensions. 
\end{numbered paragraph}

\begin{numbered paragraph} \label{smash product} {\bf Smash products}\qua 
In \cite[Thm.\ 2.2]{Lydakis:Smash}, Lydakis defines a smash product 
for $\Gamma$--spaces which is characterized by the universal property that 
$\Gamma$--space maps \mbox{$X\, \sm \, Y \to Z$} 
are in bijective correspondence with maps
\[ X(k^+) \, \sm\, Y(l^+) \quad \to \quad Z(k^+\sm\, l^+) \]
which are natural in both variables. By \cite[Thm.\ 2.18]{Lydakis:Smash}, 
the smash product of $\Gamma$--spaces is associative and 
commutative with unit $\mathbb S$, up to coherent natural isomorphism. 
There are also internal homomorphism $\Gamma$--spaces 
\cite[2.1]{Lydakis:Smash},
adjoint to the smash-product, so that the category of $\Gamma$--spaces forms 
a symmetric monoidal closed category.
\end{numbered paragraph}

\begin{numbered paragraph} \label{special} 
{\bf Special $\Gamma$--spaces}\qua  
A $\Gamma$--space $X$ is called {\em special\/} if the map 
$X(k^+{\scriptstyle \vee}\, l^+)$\break $\to X(k^+)\times X(l^+)$ 
induced by the projections from $k^+ \, {\scriptstyle \vee} \, l^+$ to $k^+$ 
and $l^+$ is a weak equivalence for all $k$ and $l$. 
In this case, the weak map
\[ \xymatrix{  X(1^+)\times X(1^+) &  
X(2^+)  \ar[l]_-{\sim} \ar[r]^-{X(\nabla)} & X(1^+) } \]
induces an abelian monoid structure on $\pi_0\,(X(1^+))$. 
Here $\nabla\co 2^+\to 1^+$ is the fold map defined by 
$\nabla(1)=1=\nabla(2)$. $X$ is called {\em very special} if it is special 
and the monoid $\pi_0\,(X(1^+))$ is a group. By Segal's theorem 
(\cite[Prop.\ 1.4]{Segal:categories}, see also \cite[Thm.\ 4.2]{BF}), 
the spectrum associated to a very special $\Gamma$--space $X$ 
is an $\Omega$--spectrum in the sense that the maps 
$|X(S^n)|\to \Omega |X(S^{n+1})|$ adjoint to the spectrum structure 
maps are homotopy equivalences.
In particular, the homotopy groups of a very special $\Gamma$--space $X$ 
are naturally isomorphic to the homotopy groups 
of the simplicial set $X(1^+)$.
\end{numbered paragraph}

\begin{numbered paragraph} \label{model structures} 
{\bf Model structures}\qua 
Bousfield and Friedlander introduce 
two model category structures for $\Gamma$--spaces called the {\em strict} 
and the {\em stable} model categories \cite[3.5, 5.2]{BF}. 
It will be more convenient for our purposes to work with slightly 
different model category structures, which we called the Quillen- 
or Q--model category structures in \cite{Sch:Gamma} and \cite{Sch:stable}. 
The strict and stable Q--structures have the same weak equivalences, 
hence the same homotopy categories, as the corresponding 
Bousfield--Friedlander model category structures.
In this paper we never consider the Bousfield--Friedlander
model structures, so we drop the decoration `Q' for the other
model structure. 

We call a map of $\Gamma$--spaces a {\em strict fibration} 
(resp.\ a {\em strict equivalence}) if it is a Kan fibration 
(resp.\ weak equivalence) of simplicial sets when evaluated 
at every $n^+\!\in\Gamma^{\text{op}}$. 
{\em Cofibrations} are defined as the maps having 
the left lifting property with respect to all strict acyclic fibrations. 
The cofibrations can be characterized as the injective maps 
with projective cokernel, see \cite[Lemma A3 (b)]{Sch:Gamma}
for the precise statement. 
By \cite[II.4 Thm.\ 4]{Quillen:HA}, the strict equivalences, 
strict fibrations and cofibrations make the category 
of $\Gamma$--spaces into a closed simplicial model category.

More important is the {\em stable} model category structure. 
This one is obtained by localizing the strict model 
structure with respect to the stable equivalences. 
We call a map of $\Gamma$--spaces a {\em stable fibration}
if it has the right lifting property with respect to the cofibrations
which are also stable equivalences.
By \cite[Thm.\ 1.5]{Sch:Gamma}, the stable equivalences, stable
fibrations and cofibrations make the category 
of $\Gamma$--spaces into a closed simplicial model category. 
A $\Gamma$--space $X$ is stably fibrant 
if and only if it is very special and $X(n^+)$ 
is fibrant as a simplicial set for all $n^+\in\Gamma\rscript{op}$.

A $\Gamma$--space $X$ defines a spectrum $X(S)$ 
(in the sense of  \cite[Def.\ 2.1]{BF})
whose $n$-th term is the value of the prolonged $\Gamma$--space at $S^n$.
For example, the $\Gamma$--space $\mathbb S$ becomes isomorphic 
to the identity functor of the category of pointed simplicial sets 
after prolongation. So the associated spectrum is the sphere spectrum. 
The functor that sends a $\Gamma$--space $X$ to the spectrum $X(S)$ 
has a right adjoint \cite[Lemma 4.6]{BF}, 
and these two functors form a Quillen pair. 
One of the main theorems of \cite{BF} says that this Quillen pair 
induces an equivalence between the homotopy category of $\Gamma$--spaces, 
taken with respect to the stable equivalences, and
the stable homotopy category of connective spectra (see \cite[Thm.\ 5.8]{BF}).
We do not use this result here, but it is the main motivation
for the study of $\Gamma$--spaces. 
\end{numbered paragraph}

\begin{numbered paragraph} \label{assembly map}{\bf The assembly map}\qua 
Given two $\Gamma$--spaces $X$ and $Y$, there is a natural map 
$X\,\sm\, Y \to X\circ Y$ from the smash product 
to the composition product \cite[2.12]{Lydakis:Smash}, \cite[1.8]{Sch:stable}.
The formal and homotopical properties of this {\em assembly map} 
are of importance to this paper.
Since $\Gamma$--spaces prolong to functors defined 
on the category of pointed simplicial sets, they can be composed. 
Explicitly, for $\Gamma$--spaces $X$ and $Y$, we set 
$(X\circ Y)(n^+)=X(Y(n^+))$. This composition $\circ$ is a monoidal 
(though not symmetric monoidal) product on the category of $\Gamma$--spaces. 
The unit is the same as for the smash product, it is the $\Gamma$--space 
$\mathbb S$ which as a functor 
is the inclusion of $\Gamma\rscript{op}$ into all pointed simplicial sets.

The assembly map is obtained as follows. 
Prolonged $\Gamma$--spaces are naturally simplicial functors \cite[\S3]{BF}, 
which means that there are natural coherent maps 
$X(K)\, \sm\, L \to X(K\sm\, L)$. This simplicial structure gives maps
\[ X(n^+)\, \sm\, Y(m^+) \ \to  \  
X(n^+\sm\, Y(m^+))\ \to \ X(Y(n^+\sm\, m^+)) \]
natural in both variables. From this the assembly map 
$X\,\sm\, Y \to X\circ Y$ is obtained by the universal property 
of the smash product of $\Gamma$--spaces. 
The assembly map is associative and unital, 
$\mathbb S$ being the unit for both $\sm$ and $\circ$. 
In technical terms: the identity functor on the category of 
$\Gamma$--spaces becomes a lax monoidal functor from ($\mathcal GS,\sm)$ 
to $(\mathcal GS,\circ)$. 
The homotopical properties of smash and composition product and
of the assembly map are summarized in the following theorem,
which is due to Lydakis \cite{Lydakis:Smash}.
\end{numbered paragraph}

\begin{theorem} \label{assembly properties}
\renewcommand{\labelenumi}{\rm(\alph{enumi})}
\begin{enumerate}
\item The composition product of $\Gamma$--spaces preserves stable equivalences
in each of its variables.
\item The smash product with a cofibrant $\Gamma$--space 
preserves stable equivalences.
\item Let $X$ and $Y$ be $\Gamma$--spaces, one of which is cofibrant.
Then the assembly map 
\[ X\,\sm\, Y \ \to \  X\circ Y \] 
is a stable equivalence.
\end{enumerate}
\end{theorem}
\begin{proof}
Parts (b) and (c) are  \cite[Prop.\ 5.12]{Lydakis:Smash} and
 \cite[Prop.\ 5.23]{Lydakis:Smash} respectively.

For any $\Gamma$--space $F$ the structure map 
$S^1\sm F(S^n)\to F(S^1\sm S^n)$ is
$(2n+1)$--connected \cite[Prop.\ 5.21]{Lydakis:Smash}.
So the map $\pi_{*}\, \Omega^n |F(S^n)| \to \pi_* \, F$ 
is an isomorphism for $\ast< n$.
Hence if $F\to F'$ is a stable equivalence of $\Gamma$--spaces, 
then the map $F(S^n)\to F'(S^n)$ is $2n$--connected.
If $X$ is any $\Gamma$--space, then the prolonged functor preserves
connectivity  \cite[Prop.\ 5.20]{Lydakis:Smash}, 
so the map $X(F(S^n))\to X(F'(S^n))$ is also
$2n$--connected. Thus the map $X\circ F\to X\circ F'$ is a stable equivalence.

It remains to show that the map $F\circ X\to F'\circ X$ is also a
stable equivalence. By the previous paragraph we may assume that $X$ is
cofibrant, and then the claim follows from parts (b) and (c).
\end{proof}

\begin{numbered paragraph} \label{Gamma-rings and modules}
{\bf Gamma-rings and their modules}\qua 
Our notion of ring spectrum is that of a {\em Gamma-ring}. 
Gamma-rings are the monoids in the symmetric monoidal category of 
$\Gamma$--spaces with respect to the smash product and 
they are special cases of `Functors with Smash Product' 
(FSPs, cf.\ \cite[1.1]{Boek:THH}, \cite[2.2]{Pira-Wald}). 
One can describe Gamma-rings as `FSPs defined on finite sets'. 
From a Gamma-ring one obtains an FSP or {\em symmetric ring spectrum} 
(\cite{HSS}, \cite[2.1.11]{Shipley:THH}) 
by prolongation and, in the second case, evaluation on spheres.
A more detailed discussion of the homotopy theory of Gamma-rings 
can be found in \cite{Sch:Gamma}.
Homotopy theoretic results about the Gamma-rings
can be translated into other frameworks
for ring spectra by the general comparison procedures of \cite{MMSS}.

Explicitly, a Gamma-ring is a $\Gamma$--space $R$ equipped with maps
\[  {\mathbb S}\ \to \ R \quad \mbox{and} \quad R\,\sm\, R \ \to\ R \ , \]
called the unit and multiplication map, which satisfy certain 
associativity and unit conditions (see \cite[VII.3]{MacL:categories}). 
A morphism of Gamma-rings is a map of $\Gamma$--spaces commuting 
with the multiplication and unit maps. If $R$ is a Gamma-ring, 
a {\em left 
$R$--module} is a $\Gamma$--space 
$N$ together with an action map $R\,\sm\, N\to N$ 
satisfying associativity and 
unit conditions (see again \cite[VII.4]{MacL:categories}).
A morphism of left $R$--modules is a map of $\Gamma$--spaces commuting 
with the action of $R$. We denote the category of left $R$--modules by $R$--mod.

One similarly defines right modules and bimodules. 
Because of the universal property of the smash product of $\Gamma$--spaces 
(\ref{smash product}), Gamma-rings are in bijective correspondence 
with lax monoidal functors from the category $\Gamma\rscript{op}$ 
to the category of pointed simplicial sets (both under smash product). 
\end{numbered paragraph}

\begin{numbered paragraph} \label{Gamma-ring examples} {\bf Examples}\qua 
The unit $\mathbb S$ of the smash product is a Gamma-ring in a unique way. 
The category of $\mathbb S$--modules is isomorphic to the category of 
$\Gamma$--spaces. 
Other standard examples of Gamma-rings are monoid rings 
over the sphere Gamma-ring $\mathbb S$ and Eilenberg--MacLane models 
of classical rings. If $M$ is a simplicial monoid, 
we define a $\Gamma$--space ${\mathbb S}[M]$ by
\[ {\mathbb S}[M]\,(n^+) \ = \ M^+  \sm \, n^+  \]
(so ${\mathbb S}[M]$ is isomorphic, as a $\Gamma$--space, 
to ${\mathbb S}\,\sm\, M^+$, see (\ref{basic definition})\,). 
There is a unit map ${\mathbb S} \to {\mathbb S}[M]$ 
induced by the unit of $M$ and a multiplication map 
${\mathbb S}[M]\, \sm \, {\mathbb S}[M] \to {\mathbb S}[M]$ 
induced by the multiplication of $M$ which turn ${\mathbb S}[M]$ 
into a Gamma-ring. This construction of the monoid ring 
over $\mathbb S$ is left adjoint to the functor which takes a 
Gamma-ring $R$ to the simplicial monoid $R(1^+)$.

If $A$ is an ordinary ring, then the Eilenberg--MacLane 
$\Gamma$--space $HA$ is given by the functor which takes a pointed set 
$K$ to the reduced free $A$--module $\widetilde A[K]$ generated by $K$.  
The unit map
\[ \eta \,\co\,K \ \to \ \widetilde A[K] \]
is the inclusion of generators. The multiplication map
\[ \mu \,\co\,\widetilde A[K]\, \sm\, \widetilde A[L] \ \to 
\ \widetilde A[K\sm\, L] \]
takes a smash product 
$(\sum_{k\in K}a_k\cdot k)\,\sm\,(\sum_{l\in L}b_l\cdot l)$ 
to the element $\sum a_kb_l\cdot (k\,\sm \, l)$.
For later reference we note that the multiplication $\mu\co HA\sm HA\to HA$
factors as the composition
\begin{equation} \label{HZ product factor} 
\xymatrix@C=16mm{ HA\sm HA \ar[r]^-{\text{assembly}} &
HA\circ HA  \ar[r]^-{\text{eval.}} &  HA } 
\end{equation}  
of the assembly map (\ref{assembly map}) and the evaluation map.

More examples of Gamma-rings arise from algebraic theories 
and as endomorphism Gamma-rings, see \cite[4.5, 4.6]{Sch:stable}. 
The Gamma-ring $DB$ (\ref{DB definition}) is such an example.
\end{numbered paragraph}

The modules over a fixed Gamma-ring and the category of all Gamma-rings 
form simplicial model categories, created by the forgetful
functor to $\Gamma$--spaces \cite[Thm.\ 2.2 and Thm.\ 2.5]{Sch:Gamma}. 
In these model structures a map is a weak equivalence (resp.\ fibration) 
if it is a stable equivalence (resp.\ stable fibration) 
as a map of $\Gamma$--spaces.
For a ring $A$, the Eilenberg--MacLane functor $H$ is the right adjoint
of a Quillen equivalence between the model category 
of $HA$--modules and the model category of simplicial $A$--modules 
\cite[Thm.\ 4.4]{Sch:Gamma}.


\section{The Gamma-ring $DB$ and formal~group~laws} 
\label{DB review}

In this section we recall the definition and some basic properties 
of our main object of study, the Gamma-ring $DB$, 
for $B$ a fixed commutative ring. This Gamma-ring parameterizes 
the stable homotopy theory of augmented commutative simplicial $B$--algebras. 
In Construction \ref{Gamma-ring map} we explain
how a formal group law $F$ over $B$ gives rise to a homomorphism 
of ring spectra $F_*\co H\mZ\to DB$. 
The rest of this paper is then devoted to studying 
the homotopical significance of that construction.

By ``parameterizing the stable homotopy theory'' 
we mean that there is a Quillen-equivalence between 
the model category of $DB$--modules and the model category of spectra 
of simplicial commutative $B$--algebras,
see Theorem~\ref{DB summary}~(d). 
Commutative simplicial algebras have been the object of much study 
\cite{Quillen:comrings, Dwyer:oper, Goerss:Hilton-Milnor, 
Goerss:AQ-homology, Sch:modelspec, Turner:looking_glass}. 
The homology theory arising as the derived functor of abelianization 
in this category is known as Andr\'e--Quillen homology for commutative rings.
The stable homotopy category of simplicial commutative $B$--algebras 
is equivalent to the homotopy category of $DB$--modules. 
The homotopy groups of $DB$ are the coefficients 
of the universal homology theory for commutative algebras.

\begin{numbered paragraph} \label{DB definition} 
{\bf The Gamma-ring $DB$}\qua
The $\Gamma$--space underlying the Gamma-ring $DB$ takes a pointed set $K$ 
to the augmentation ideal of the power series ring generated by $K$, 
considered as a constant simplicial set:
\[ DB(K) \ = \ \mbox{kernel}\,(\,\widetilde{B}[\![K]\!] \to 
\widetilde{B}[\![\ast]\!] = B\, ) \ . \]
The tilde over $\widetilde{B}[\![K]\!]$ indicates that the power series 
generator corresponding to the basepoint of $K$ has been set equal to 0;
thus $\widetilde{B}[\![\ast]\!]$ reduces to the coefficient ring.

The product which makes $DB$ into a Gamma-ring comes from 
substitution of power series. To define the multiplication map
$\mu\co DB\sm DB\to DB$ we need to describe a natural associative map
\[ \mu \,\co\,DB(K)\, \sm\, DB(L) \ \to \ DB(K\,\sm\, L) \]
for pointed sets $K$ and $L$. An element of $DB(K)$ is represented 
by a power series $f$ in variables $K$ without constant term. 
Similarly, an element of $DB(L)$ is represented by a power series $g$ 
in variables $L$. The multiplication map $\mu$ takes the smash product 
\[       f\, \sm\, g \ \in \  DB(K)\, \sm\, DB(L) \]
to the power series $\mu(f\sm\, g)$ in the variables 
$K\sm\, L$ defined by
\[ \mu(f(k_1,\dots,k_m)\,\sm\, g(l_1,\dots,l_n)) \quad = \hspace*{4cm} \]
\[ \hspace*{2cm}
f(\, g(k_1\sm\, l_1,\dots,k_1\sm\, l_n),\, \dots\, ,g(k_m\sm\, l_1,\dots,k_m\sm\, l_n)\,)
 \ . \]
The unit map $\eta\co \mS\to DB$, ie, natural transformation
\[ \eta\, \co\,  K \to DB(K) \ , \]
sends an element of $K$ to the generator it represents. 
The multiplication and unit transformations $\mu$ and $\eta$ are associative 
and unital because substitution of power series is, so $DB$ is in fact a 
Gamma-ring.
\end{numbered paragraph}

Some important properties of $DB$ are summarized in the following theorem.
Most of these results are compiled 
from \cite{Sch:stable}. In \cite[7.9]{Sch:stable}, we use the
notation $DB$ for a slightly different Gamma-ring, namely the Gamma-ring
associated to the algebraic theory of augmented commutative $B$--algebras.
Let $DB\rscript{pol}$ denote the sub-$\Gamma$--space of $DB$ whose
value at a pointed set $K$ consists of all {\em polynomials} in $K$, ie,
the power series in $DB(K)$ with only finitely many non-zero coefficients.
The sub-$\Gamma$--space $DB\rscript{pol}$ is closed under the multiplication
of $DB$, and the unit map $\eta\co\mS\to DB$ factors through $DB\rscript{pol}$.
Hence $DB\rscript{pol}$ is a Gamma-ring and the inclusion 
$DB\rscript{pol}\to DB$ is a homomorphism. 
$DB\rscript{pol}$ is exactly the Gamma-ring which is denoted $DB$ 
in \cite[7.9]{Sch:stable}; there is no homotopical
difference between the two Gamma-rings since the inclusion 
$DB\rscript{pol}\to DB$ is a stable equivalence, 
see Theorem \ref{DB summary} (a) below.
However, in this paper it is more convenient to work with the power series 
model, so we give it the simpler name.

\begin{theorem} \label{DB summary}
\renewcommand{\labelenumi}{\rm(\alph{enumi})}
\begin{enumerate}
\item The inclusion $DB\rscript{\em pol}\to DB$ 
is a stable equivalence of Gamma-rings. 
\item As a $\Gamma$--space, $DB$ is stably equivalent to the derived
smash product of the Eilenberg--MacLane $\Gamma$--spaces $H\mZ$ and $HB$;
in particular, the stable homotopy groups of $DB$ are additively isomorphic
to the spectrum homology of the  Eilenberg--MacLane spectrum of $B$.
\item The graded ring of homotopy groups of $DB$ is isomorphic 
to the ring of stable homotopy operations of commutative augmented 
simplicial $B$--algebras.
\item There is a Quillen adjoint functor pair between the model category of 
$DB$--modules and the model category ${\mathcal S}p(B\mbox{\em -alg})$ of
spectra of  commutative augmented simplicial $B$--algebras.
The adjoint pair passes to an equivalence between the homotopy category 
of $DB$--modules and the homotopy category of connective spectra 
of commutative augmented simplicial $B$--algebras.
\end{enumerate}
\end{theorem}
\begin{proof}
(a)\qua Let $(B\tensor {\mathcal S}^k)^!$ denote the $\Gamma$--space 
defined by
\[ (B\tensor {\mathcal S}^k)^!(K) \ = \ 
B \tensor (\widetilde{\mZ}[K]^{\tensor k}/\Sigma_k) \ , \] 
the tensor product of $B$ with the $k$-th symmetric power 
of the reduced free abelian group on $K$.
An isomorphic description of $(B\tensor {\mathcal S}^k)^!(K)$ is as
the free reduced $B$--module on the $k$-th symmetric power of $K$,
\[ (B\tensor {\mathcal S}^k)^!(K) \ \iso \ \widetilde{B}[K^{\sm k}/\Sigma_k] \ . \] 
The underlying $\Gamma$--space of $DB$ 
is the infinite product of the symmetric power functors  $(B\tensor {\mathcal S}^k)^!$
for all $k\geq 1$.
The polynomial model $DB\rscript{pol}$ is the {\em weak} product 
of these symmetric power functors.
Since the stable homotopy groups of the $\Gamma$--space $(B\tensor {\mathcal S}^k)^!$
are trivial up to dimension $2k-3$ \cite[12.3]{Dold-Puppe}, 
the weak product and the product are stably equivalent.

(b)\qua We let $SP$ denote the $\Gamma$--space which
takes a pointed set $K$ to its infinite symmetric product,
ie, the free abelian monoid generated by $K$ with basepoint as identity
element.  
By the Dold--Thom theorem, the group completion map $SP\to H\mZ$
is a stable equivalence of $\Gamma$--spaces.
We choose a cofibrant replacement $HB\rscript{c}$ of $HB$ as a $\Gamma$--space
and obtain a chain of homomorphisms of $\Gamma$--spaces
\begin{eqnarray*}
 HB\rscript{c} \sm H\mZ & \varl{1cm}{\sim} &  HB\rscript{c} \sm SP \ 
\varr{1cm}{\sim} \  HB\rscript{c} \circ SP \\
& \varr{1cm}{\sim} &
  HB \circ SP \ \iso \ DB\rscript{pol} \ \varr{1cm}{\sim} \ DB \ . 
\end{eqnarray*}
The first map is a stable equivalence since smashing with a cofibrant 
$\Gamma$--space preserves stable equivalences \cite[5.12]{Lydakis:Smash}.
The second map is the assembly map (\ref{assembly map}), which is 
a stable equivalence by \cite[5.23]{Lydakis:Smash}.
The third map is a stable equivalence since the composition product
of $\Gamma$--spaces preserves stable equivalences in both variables
(Theorem \ref{assembly properties} (a)).
The reduced free $B$--module generated by $SP(K)$ 
is isomorphic to the polynomials without constant term generated by $K$,
subject to the basepoint relation. This gives as isomorphism 
of $\Gamma$--spaces between  $HB \circ SP$ and $DB\rscript{pol}$.
The  last map is a stable equivalence by part (a).

Part (c) and (d) are special cases of \cite[4.11]{Sch:stable}.
and  \cite[4.4]{Sch:stable} respectively; see also \cite[7.9]{Sch:stable}.
\end{proof}

The ring of stable homotopy operations of commutative augmented 
simplicial $B$--algebras
--- ie, the graded ring of homotopy groups of $DB$ ---
is sometimes called the stable Cartan--Bousfield--Dwyer-algebra since 
these authors calculated the unstable operations for $B=\mF_p$, 
see \cite{Cartan, Bous:operations, Dwyer:oper}. 
An explicit calculation of $\pi_*\, D\mF_p$ can be found 
as Theorems 12.3 (for $p=2$) and 12.6 (for $p$ odd)
in Bousfield's unpublished paper \cite{Bous:operations}. 
The fact that the ring of  stable homotopy operations is generally not
commutative shows that $DB$ is {\em not} stably equivalent  
to the derived smash product $HB\sm^L H\mZ$ as a Gamma-ring (unless
$B$ is a $\mQ$--algebra). 

\begin{numbered paragraph} \label{DB and TAQ}
{\bf Relation to topological Andr\'e-Quillen homology and $\Gamma$-homology}\qua
There are isomorphisms of graded abelian groups
\begin{equation}\label{DB, HGamma and TAQ} 
\pi_{\ast} DB \ \iso \ H\Gamma_*(B[x]|B;B) \ \iso \ TAQ_*(HB[x]|HB;HB) \ .
\end{equation}
Here $H\Gamma_*(B[x]|B;B)$ is the {\em $\Gamma$--homology} in the sense of 
Robinson and Whitehouse~\cite{Robinson-Whitehouse:GammaHomology},
of the polynomial algebra $B[x]$, considered as
an augmented $B$--algebra; moreover, $TAQ_*(HB[x]|HB;HB)$ denotes the
{\em topological Andr\'e--\break Quillen homology}~\cite{Basterra:AQHomology}
of the Eilenberg--MacLane ring spectrum $HB[x]$ relative to $HB$.
Both $\Gamma$--homology and topological Andr\'e--Quillen homology groups
are studied because they carry obstructions to the existence 
of $E_{\infty}$ ring spectrum structures~\cite{Robinson:Gamma-Einfty}
(it is a coincidence that the symbol $\Gamma$ occurs both in $\Gamma$--homology
and as the category $\Gamma\rscript{op}$).

The first isomorphism in \eqref{DB, HGamma and TAQ} comes about as follows.
By \ref{DB summary} (a) above, the Gamma-ring $DB$ 
has a stably equivalent polynomial model $DB\rscript{pol}$;
as a $\Gamma$--space, $DB\rscript{pol}$ is isomorphic to the functor
which assigns to the object $n^+$ of $\Gamma\rscript{op}$ the $B$--module
$B[x]^{\tensor_B n}$. By a theorem of Pirashvili and Richter
\cite[Thm.~1]{Pira-Richter}, the homotopy groups 
$\pi_{\ast}DB$ are thus isomorphic to the $\Gamma$--homology
$H\Gamma_*(B[x]|B;B)$ of the polynomial algebra $B[x]$ relative to $B$. 
The second isomorphism in \eqref{DB, HGamma and TAQ} is due to Basterra 
and McCarthy \cite{Basterra-McCarthy}, who show
that for Eilenberg--MacLane spectra of classical rings, 
topological Andr\'e--Quillen homology coincides with $\Gamma$--homology. 
The survey article~\cite{Basterra-Richter} by Basterra and Richter
discusses all these identifications in more detail.

The isomorphisms  \eqref{DB, HGamma and TAQ}
do not mention the multiplicative structure of $\pi_{\ast} DB$.
In Sections 5.1 and 7.9 of \cite{Sch:stable} we associate to
any augmented commutative $B$--algebra $A$ a $DB$--module 
$\Sigma^{\infty}_BA$ which models the suspension spectrum of $A$ as
an augmented commutative $B$--algebra. The underlying $\Gamma$--space of
$\Sigma^{\infty}_BA$  sends $n^+\in\Gamma\rscript{op}$ to the $B$--module
$A^{\tensor_B n}$.
So  \cite[Thm.~1]{Pira-Richter} yields an isomorphism
of graded abelian groups
\[ \pi_{\ast}\left(\Sigma^{\infty}_BA \right)\ \iso \ H\Gamma_*(A|B;B) \]
(at least if $A$ is flat as a $B$--module). 
The $DB$--action on $\Sigma^{\infty}_BA$ gives more structure 
to $\Gamma$--homology, since the left hand side above has
a natural action of the graded ring $\pi_{\ast}DB\iso H\Gamma_*(B[x]|B;B)$.
We see moreover that the homotopical object underlying $\Gamma$--homology 
is not just a chain complex, but a $DB$--module spectrum.
\end{numbered paragraph}

The main objective of this paper is the study of the homotopy type of $DB$ 
as a ring spectrum and of a close relationship to formal group theory. 
More precisely 
we will describe the space of Gamma-ring maps from $H\mZ$, 
the Eilenberg--MacLane Gamma-ring of the integers (\ref{Gamma-ring examples}), 
to $DB$ in terms of formal group law data. Unless stated otherwise, 
formal group laws will always be {\em 1--dimensional and commutative}.
To see how non-commutative and higher-dimensional formal group laws fit 
into our context see Remarks \ref{non-commutative} and
\ref{higher dimensional}.

\begin{construction} \label{Gamma-ring map}
Suppose that $F$ is a 1--dimensional and commutative formal group law over 
the commutative ring $B$. In other words, $F$ is a power series in 
two variables $x$ and $y$ with coefficients in $B$ which satisfies
\begin{eqnarray*}
F(x,0)  & = & x \ = \ F(0,x) \ , \\
F(x,y) & = & F(y,x) \mbox{\quad and} \\
F(F(x,y),z) &  = & F(x,F(y,z)) \ . 
\end{eqnarray*}
We define a map $F_{\ast}\co H\mZ\to DB$ 
of Gamma-rings. For every pointed set $K$ we have to specify a map
\[ F_{\ast}(K) \,\co\, \widetilde\mZ[K] = H\mZ(K) \ \to \ 
 DB(K) \ \subseteq \ \widetilde{B}[\![K]\!] \]
which is natural in $K$ and respects the multiplication and unit maps. 
The map $F_{\ast}(K)$ simply takes a sum $\sum_{k\in K} a_k\cdot k$ 
of generators of the free abelian group on $K$ to the formal sum
\[ {\sum_{k\in K}}^F [a_k]_F(k) \]
with respect to $F$, of the same elements viewed as generators of the 
power series ring. Here $[n]_F$ denotes the $n$--series of the formal group law
$F$ for every integer $n\in\mZ$. We omit the verification that the map 
$F_{\ast}$ indeed commutes with the multiplication and unit map. 
For example, on the level of underlying monoids this means that the map 
\[ [-]_F \,\co\,  \mZ  \, \iso  \, H\mZ(1^+) \ \to 
\  DB(1^+) \, = \,  x\cdot B[\![x]\!] \]
is a homomorphism from the multiplicative monoid of the integers to 
the monoid of power series without constant term under substitution,
ie, it boils down to the relation 
\[ [n]_F([m]_F(x)) \ = \ [n\cdot m]_F(x) \]
for $n,m\in\mZ$.
\end{construction}

\begin{remark} We offer two additional ways of looking at the above
construction of the Gamma-ring map $F_*$; the three definitions correspond
to looking at a commutative 1--dimensional formal group law as either
\begin{enumerate}
\item[(a)] a power series $F(x,y)$ in two variables with certain properties,
\item[(b)] an abelian cogroup structure, in the category of complete, augmented
commutative $B$--algebras, on the power series ring $B[\![x]\!]$,
\item[(c)] or a morphism of algebraic theories 
from the theory of abelian groups to the theory of complete, augmented, 
commutative $B$--algebras.
\end{enumerate}
The first point of view leads to the explicit formula for the Gamma-ring map
$F_*$ that was just given in Construction \ref{Gamma-ring map}.

The interpretation (b) of a formal group law exhibits the Gamma-ring map
$F_*$ as a special case of a more general construction 
associated to an abelian cogroup object. Indeed, to every object $X$ in 
a pointed category with coproducts $\C$ one can associate an {\em endomorphism
Gamma-ring} $\End_{\C}(X)$
(see \cite[4.6]{Sch:stable} or \ref{coordinate free}).
Then every abelian cogroup structure on the object $X$ leads to a map of
Gamma-rings $H\mZ\to \End_{\C}(X)$; we refer to  \ref{coordinate free}
for more details. 

Perspective (c) leads to a compact description 
in the language of algebraic theories \cite[3.3.1]{Borceux:2}.
Specifying a 1--dimensional, commutative formal group law over $B$ is
the same as specifying a morphism of algebraic theories \cite[3.7.1]{Borceux:2}
from the theory of abelian groups to the theory of complete,
augmented, commutative $B$--algebras. The construction \cite[4.5]{Sch:Gamma}
which associates to a pointed algebraic theory $T$ its stable Gamma-ring 
$T^s$ is functorial for morphisms of algebraic theories.
Now  $H\mZ$ is the Gamma-ring associated to the theory of abelian groups
and $DB$  is the Gamma-ring associated to the theory of
complete, augmented, commutative $B$--algebras.
Hence a formal group law $F$ defines a morphism of algebraic theories,
thus a morphisms of associated Gamma-rings.
\end{remark}

\begin{example} The {\em additive formal group law} is given by
\[ F^{a}(x,y) \ = \ x + y \ . \]
The associated   Gamma-ring map $F^{a}_*$ is the composite 
\[ \xymatrix{ H\mZ \ar[r] & HB \ar[r]^-{\text{incl.}} & DB } \]
of the unique Gamma-ring map $H\mZ\to HB$ with the `inclusion' of $HB$
into $DB$ as the linear power series.
Conversely, this is the only case in which the
Gamma-ring map $F_*$ factors over the inclusion $HB\to DB$ 
on the point-set level:
the power series $F(x,y)$ can be recovered from the Gamma-ring map
$F_*$ as the image of $x+y\in \mZ[x,y]$ under
\[ F_* \,\co\, \mZ[x,y] \ \iso \ H\mZ(\{x,y\}^+) \ \to 
\ DB(\{x,y\}^+) \ \subseteq B[\![x,y]\!] \ ; \] 
so if $F_*$ factors over $HB$, then $F(x,y)$ has only linear terms,
so that necessarily $F(x,y)=x+y$.
 
Another formal group law which exists over any ring $B$ 
is the {\em multiplicative} one given by
\[ F^{m}(x,y) \ = \ x + y + xy \ . \]
The multiplicative formal group law can be used to
express the Gamma-ring $DB$ additively as the smash product of
the Eilenberg--MacLane Gamma-rings for $B$ and $\mZ$;
indeed the composite map
\[\xymatrix@C=10mm{ 
HB \sm^L H\mZ \ar[rr]^-{\text{incl.}\ \sm^L\  F^{m}_*} && 
DB \sm^L DB  \ar[r]^-{\mu} & DB }\]
is in the same homotopy class as the stable equivalence of 
Theorem \ref{DB summary} (b).
\end{example}

The homotopical significance of the point-set level construction of 
\ref{Gamma-ring map} is summarized in the following 

\begin{theorem} \label{homotopy classes} 
Construction \ref{Gamma-ring map} which to a 1--dimensional,
commutative formal group law $F$ associates 
the Gamma-ring map $F_{\ast}\co H\mZ\to DB$ induces a natural bijection
\[\xymatrix{ \mbox{\em FGL}(B)/\mbox{\em strict isomorphism}\ \ar[r]^-{\iso}
& \ [H\mZ,DB]_{\text{\em Gamma-rings}}} \]
between the strict isomorphism classes of formal group laws over $B$ and 
the set of maps from $H\mZ$ to $DB$ in the homotopy category of Gamma-rings. 
\end{theorem}

Theorem \ref{homotopy classes} is the $\pi_0$--part of a 
space level statement relating formal groups to
Gamma-ring maps between $H\mZ$ and $DB$ in Theorem \ref{main}. 
In the simplicial model category 
of Gamma-rings every pair of objects has a homomorphism space 
(ie, simplicial set).
As usual with model categories, in order to give 
the morphism space a homotopy invariant meaning, the source object has 
to be replaced by a weakly equivalent cofibrant one, 
and the target object has to be replaced by a weakly equivalent fibrant one.
Then the components of the (derived) space of Gamma-ring maps
are the morphisms  in the homotopy category of Gamma-rings. 
Theorem \ref{main} below identifies the derived homomorphism space 
of Gamma-ring maps from $H\mZ$ to $DB$. The answer is given
in terms of formal group law data and the simplicial monoid of 
{\em homotopy units} of $DB$.

\begin{numbered paragraph} \label{unit definition} {\bf Homotopy units}\qua
Let $R$ be a Gamma-ring and  $R^{\f}$ a stably fibrant replacement 
of $R$ in the model category structure of \cite[2.5]{Sch:Gamma}.
As for any Gamma-ring, the underlying space $R^{\f}(1^+)$ 
is a simplicial monoid with product induced by the 
multiplication map $R^{\f}\sm \,R^{\f}\to R^{\f}$. 
Moreover, $R^{\f}(1^+)$ is a model for the infinite loop space of the
spectrum represented by $R$.
We define the {\em homotopy units} $R^{\times}$ 
as the union of the invertible components 
of the simplicial monoid $R^{\f}(1^+)$. 
So $R^{\times}$ is a group-like simplicial monoid which is independent
up to weak equivalence of the choice of fibrant replacement.
Moreover there are natural isomorphisms of homotopy groups 
\[ \pi_0 \, R^{\times} \ \iso  \ \mbox{units\,}(\pi_0 \, R) 
\mbox{\quad and \quad}
\pi_i \, R^{\times} \ \iso  \ \pi_i \, R \mbox{\quad for $i\geq 1$.} \]
For any classical ring $A$, the units of $A$ act by conjugation on $A$ 
and hence on the set of ring homomorphisms from any other ring to $A$. 
In Section \ref{conjugation} we make sense of the analogous conjugation 
action for Gamma-rings. 
For any Gamma-ring $R$ we construct, after change of models, 
a strict action of the homotopy units  $R^{\times}$ on $R$ by Gamma-ring
homomorphisms.
More precisely we introduce a simplicial group $UR^{\times}$, 
weakly equivalent to 
$R^{\times}$, and a fibrant Gamma-ring, stably equivalent to $R$, 
and on which $UR^{\times}$ acts by conjugation.
Below we only use the action of the simplicial subgroup  $UR^{\times}_1$,
the connected component of the unit element. 
We can consider the homotopy orbit space 
of the conjugation action of $UR^{\times}_1$ on the simplicial set
$\Rg(H\mZ,R)=\map_{\GR}(H\mZ\rscript{c},R^{\f})$. 
We denote that homotopy orbit space by $\Rg(H\mZ,R)/\mbox{conj.}$

Construction \ref{Gamma-ring map} associates to every formal group law over 
$B$ a Gamma-ring map from the Eilenberg--MacLane Gamma-ring $H\mZ$ to $DB$. 
This map gives rise to a point in the space of Gamma-ring maps 
$\Rg(H\mZ,DB)$. In Section \ref{comparison map} we extend this to 
a natural weak map from the classifying space $\mathcal{FGL}\rscript{str}(B)$ 
of the groupoid of formal group laws and strict isomorphisms to the homotopy 
orbits space of the conjugation action, ie, we construct a diagram
of simplicial sets
\[ \xymatrix{  \mathcal{FGL}\rscript{str}(B) &
\widetilde{\mathcal{FGL}}\rscript{str}(B)  \ar[l]_-{\sim} 
\ar[r]^-{\kappa} &  \Rg(H\mZ,DB)/\mbox{conj.}}  \]
in which the first map is a weak equivalence. 
The main result of this paper, Theorem \ref{main},
says that the map $\kappa\co\widetilde{\mathcal{FGL}}\rscript{str}(B)\to
\Rg(H\mZ,DB)/\mbox{conj.}$ is a weak equivalence.
Theorem \ref{homotopy classes} is the bijection induced on
path components by the weak equivalence $\kappa$.
\end{numbered paragraph}

\section{The conjugation action} \label{conjugation}

In this section we construct the conjugation action of the homotopy units 
of a Gamma-ring on the Gamma-ring. In view of the 
application to $DB$ we need a construction relative 
to a group which maps to the multiplicative monoid of the Gamma-ring.
In the example of $DB$ that group is 
the group of invertible power series on one generator over $B$. 

\begin{construction} \label{con-action of group}
We consider a Gamma-ring $R$ together with a simplicial group $G$ and
a homomorphism of simplicial monoids $\psi\co G\to R(1^+)$.
The conjugation action of $G$ on $R$ is described by a map of 
simplicial monoids
\[ c \,\co\,  G \ \to \ \map_{\mathcal{GR}}(R,R) \ . \]
Here for $g\in G$ the conjugation map $c(g)\co R\to R$ is defined 
at a pointed set $K$ as the composite
\begin{eqnarray*} R(K) \ 
\varr{1.8cm}{g\,\sm\,\text{id}\,\sm\, g^{-1}} \ 
G_+\sm\, R(K) \,\sm \, G_+  
& \varr{1.5cm}{\psi\, \sm\,\text{id}\,\sm\, \psi} &
R(1^+)\,\sm\, R(K) \, \sm\, R(1^+) \\ 
& \varr{1.5cm}{\text{mult.}} & R(1^+\sm K\sm 1^+)\ \iso \ R(K) \ . 
\end{eqnarray*}
We omit the formal verification that the conjugation map is in fact 
a homomorphism of Gamma-rings, that the definition extends 
to higher dimensional simplices of $G$, and that the formula 
$c(g\cdot g')=c(g)\circ c(g')$ holds. 
The monoid map $c$ can now be used 
to let the group $G$ act on the space of Gamma-ring maps 
from any Gamma-ring $S$ to $R$ via 
\[ G \ \times \ \map_{\mathcal{GR}}(S,R)\ \to \ 
\map_{\mathcal{GR}}(S,R) \  , \quad  (g,f) \ \mapsto \ c(g)\circ f \ . \]
\end{construction}

The goal of this section is to extend the conjugation action 
from the given group $G$ to the homotopy units \eqref{unit definition} of $R$. 
The problem is that Construction~\ref{con-action of group}
makes use of strict inverses, whereas the homotopy units 
$R^{\times}$ are only a group-like simplicial monoid.
One way to solve this would be to find 
a stable equivalence of Gamma-rings from $R$ 
to some stably fibrant Gamma-ring $R^{\f}$ 
which has the property that every element in an invertible component 
of $R^{\f}(1^+)$ has a point-set inverse. 
But it seems unlikely that this can be done in general, 
and we use a different approach. 

\begin{construction} \label{conjugation construction}
As above we consider a Gamma-ring $R$ together with a simplicial group $G$ and
a homomorphism of simplicial monoids $\psi\co G\to R(1^+)$.
If $R\to R^{\f}$ is a stably fibrant 
replacement of $R$ in the model category of Gamma-rings then $R^{\times}$ 
was defined in \ref{unit definition} as the simplicial 
monoid of invertible components in $R^{\f}(1^+)$. 
The image of the simplicial group $G$ under the map 
$R(1^+)\to R^{\f}(1^+)$ is contained in the invertible components, 
which provides a morphism of simplicial monoids $G\to R^{\times}$. 
We explained in \ref{con-action of group} how the group $G$ acts 
on spaces of Gamma-ring maps into $R$, and we now want to extend 
this to a conjugation action of the homotopy units $R^{\times}$.

We start by factoring the homomorphism $G\to R^{\times}$ 
in the model category of simplicial monoids 
(a special case of \cite[II.4 Thm.\ 4]{Quillen:HA}) 
as a cofibration followed by an acyclic fibration
\[\xymatrix{ G \ \ar@{>->}[r] & \ cR^{\times} \ \ar@{->>}^{\sim}[r] & \ 
R^{\times}. }\]
We denote by $UR^{\times}$ the algebraic group completion 
of the simplicial monoid $cR^{\times}$. So $UR^{\times}$ is obtained 
from $cR^{\times}$ by formally adjoining inverses in every simplicial 
dimension. Since the map $G\to cR^{\times}$ is a cofibration of 
simplicial monoids,  the monoid $cR^{\times}$ is a retract of a simplicial 
monoid which is dimensionwise a free product of a group and 
a free monoid \cite[II p.\ 4.11 Rem.\ 4]{Quillen:HA}. 
By the following Lemma the group completion map 
$cR^{\times}\to UR^{\times}$ is thus a weak equivalence.
\end{construction}

\begin{lemma} \label{good monoid proof} 
Let $M$ be a simplicial monoid which in every dimension 
is a free product of a group and a free monoid,
and such that $\pi_0 M$ is a group. 
Then the group completion map $M\to UM$ is a weak equivalence.
\end{lemma}
\begin{proof} We call a simplicial monoid $N$ {\em good} 
if the group completion map 
$N\to UN$ induces a weak equivalence $BN\to BUN$ of classifying spaces. 
For any simplicial monoid $N$ whose components form a group, 
$BN$ is a delooping of $N$, ie, the map $|N|\to\Omega|BN|$ 
is a weak equivalence; this follows for example by applying \cite[B.4]{BF}
to the sequence of bisimplicial sets $N\to E_{\bullet}N \to B_{\bullet}N$. 
So it suffices to show that any monoid 
as in the statement of the lemma is good.

Suppose $M$ and $N$ are good, discrete monoids. We claim that then
the free product $M\ast N$ is also good. 
By \cite[Lemma 4]{McDuff:monoids}, the canonical map of simplicial sets 
$BM\,\Wedge\, BN\to B(M\ast N)$ is a weak equivalence.
Since $U(M\ast N)\iso UM\ast UN$ (the coproduct in the category of groups 
coincides with the free product of underlying monoids), 
the map $BUM\,\Wedge\, BUN\to BU(M\ast N)$ is also a weak equivalence
and the claim follows. 

Every group, viewed as a constant simplicial monoid, is good.
The classifying spaces of the free monoid on one generator and of the free
group on one generator are both weakly equivalent to a circle. 
So a free monoid on one generator is good.
By the above and by direct limit, a free product of a group with
a free monoid is good.

The classifying space $BN$ of a simplicial monoid $N$ is the diagonal
of the bisimplicial set given by the classifying spaces $BN_m$
of the individual monoids $N_m$ in the various simplicial dimensions.
Hence if $N$ is a simplicial monoid such that the discrete monoid $N_m$ 
is good for all $m\geq 0$, then $N$ itself is good. This proves the lemma.
\end{proof}

We take the adjoint $\mS[cR^{\times}]\to R^{\f}$ of the monoid map 
$cR^{\times}\to R^{\times}\to R^{\f}(1^+)$ 
(where $\mS[cR^{\times}]$ is the  monoid Gamma-ring 
(\ref{Gamma-ring examples})) and factor it 
in the model category of Gamma-rings as a cofibration 
followed by an acyclic fibration
\[\xymatrix{ \mS[cR^{\times}] \ \ar@{>->}[r] & \ R_1 \ar@{->>}^{\sim}[r] &
\ R^{\f} \ . }\]
We then define another Gamma-ring $R_2$ as the pushout, 
in the category of Gamma-rings, of the diagram
\[\xymatrix{\ \mS[cR^{\times}]\ \ar@{>->}[r] \ar_{\sim}[d] & R_1 \ar[d] \\
\ \mS[UR^{\times}]\  \ar[r] & \ R_2 \ . }\]

For every pointed simplicial set $K$ the $\Gamma$--space $\mS\sm K$ is 
cofibrant; in particular the underlying $\Gamma$--spaces of
$\mS[cR^{\times}]$ and $\mS[UR^{\times}]$ are cofibrant.
By \cite[4.1 (3)]{SS:monoidal} the underlying  $\Gamma$--space of
$R_1$ is also cofibrant, so by the following lemma the map 
$R_1\to R_2$ is a stable equivalence of Gamma-rings.

\begin{lemma} \label{flat pushout proof} 
Consider a diagram of Gamma-rings
\[ X \ \varl{1cm}{\sim} \ Y \ \to \ Z \]
in which the left map is a stable equivalence, the right map
is a cofibration and
all three Gamma-rings are cofibrant as $\Gamma$--spaces.
Then the map from $Z$ to the pushout of the diagram is also 
a stable equivalence.
\end{lemma}
\begin{proof} We denote the pushout of the diagram by $P$.
We first consider the situation where the map $Y\to Z$ 
is obtained by cobase change from a generating cofibration. 
In other words we assume that there exists a cofibration $K\to L$ 
of $\Gamma$--spaces such that $Z$ is the pushout of the diagram of Gamma-rings
\[ Y \ \varl{1cm}{} \ T(K) \ \to \ T(L) \]
where $T$ denotes the tensor algebra functor. This special kind of pushout 
in the category of Gamma-rings is analyzed in the proof of 
\cite[Lemma 6.2]{SS:monoidal}. The pushout $Z$ is then the colimit 
of a sequence of cofibrations of $\Gamma$--spaces
\[ Y=Z_0 \ \to \ Z_1 \ \to \ \cdots \ \to Z_n \ 
\to \ \cdots \]
for which the subquotient $Z_n/Z_{n-1}$ is isomorphic 
to $(L/K)^{\sm n}\sm\, Y^{\sm (n+1)}$. 
In the same way the composite pushout $P$ is the colimit of a sequence 
of cofibrations of $\Gamma$--spaces with subquotients isomorphic to 
$(L/K)^{\sm n}\sm\, X^{\sm (n+1)}$. The smash product of $\Gamma$--spaces 
preserves stable equivalences between cofibrant objects 
\cite[Thm.\ 5.12]{Lydakis:Smash}, 
so the induced maps on the subquotients of the filtrations 
for $Z$ and $P$ are all stable equivalences. So the maps induced 
on all finite stages, and finally the map $Z\to P$ on colimits 
are also stable equivalences.

By induction the lemma thus holds whenever the cofibration $Y\to Z$ 
is the composite of finitely many maps obtained by cobase changes 
from generating cofibrations. Since homotopy groups of $\Gamma$--spaces 
commute with transfinite compositions over cofibrations, 
the lemma holds whenever $Y\to Z$ is such a transfinite composition 
of cobase changes of generating cofibrations. Finally, if the lemma holds 
for a cofibration $Y\to Z$, then it also holds for any retract. 
But all cofibrations of Gamma-rings are obtained by a sequence 
of these constructions by the small object argument.
\end{proof}

We finally let $R_3$ be a stably fibrant replacement of the Gamma-ring 
$R_2$. We can display all the relevant objects thus constructed 
in a commutative diagram of Gamma-rings
\[\xymatrix@C=15mm{\mS[G]\ \ar[rr] \ar[d] & & R \ar@{>->}^{\sim}[d] \\
\mS[cR^{\times}]\ \ar_{\sim}[d]  \ar@{>->}[r] & 
R_1 \ar@{->>}^{\sim}[r] \ar^{\sim}[d] & R^{\f} \\
\mS[UR^{\times}]\ \ar[r] & R_3 }\]
In this diagram the Gamma-rings $R^{\f}, R_1$ and $R_3$ are fibrant, 
and the maps between them are stable equivalences. 
Furthermore, the induced map from the simplicial group $UR^{\times}$ 
to the invertible components of the underlying monoid of $R_3$ 
is a weak equivalence. 
As described in \ref{con-action of group},
the simplicial group $UR^{\times}$ acts by conjugation on $R_3$ 
via homomorphisms of Gamma-rings, and this action extends the action of $G$.

\begin{remark}
Construction \ref{conjugation construction} can be made functorial 
in the triple $(R,G,\psi: G\to R(1^+))$
since the factorizations in the model categories of simplicial monoids
and of Gamma-rings can be made functorial.
\end{remark}

\section{The comparison map} \label{comparison map}

We now use the conjugation action of the previous section  
in the case of the Gamma-ring $DB$. 
We obtain a model for the homotopy invariant space of Gamma-ring maps 
from $H\mZ$ to $DB$ on which the homotopy units of $DB$ act by conjugation. 
So we can form the homotopy orbit space with respect to the conjugation action.
We then construct a weak map from the classifying space of the groupoid 
of formal group laws and strict isomorphisms to the homotopy orbit space 
of the conjugation action. 
In the remaining sections we show that that map is a weak equivalence.

In the later sections we need a version of the weak map 
from the classifying space $\mathcal{FGL}\rscript{str}(B)$ 
to the homotopy orbit space 
for other Gamma-rings. So we set up the construction of the (weak) map 
in a slightly more general context.

\begin{construction} \label{construction weak map} 
Again we consider a Gamma-ring $R$, a simplicial group $G$ and 
a homomorphism of simplicial monoids $\psi\co G\to R(1^+)$.
We now make the additional assumption that {\em the image of $G$ lands
in the unit component  of $R$}, ie, the composite
map
\[ G \ \varr{1cm}{\psi} \ R(1^+) \ \to \ \pi_0 \, R \]
is constant with value $1\in \pi_0 R$. 
This assumption is not very important, but it slightly simplifies
certain arguments later.
Moreover we are given a simplicial set $X$ with an action 
of the simplicial group $G$ 
and suppose we are are also given a $G$--equivariant map 
$X\to \map_{\mathcal{GR}}(H\mZ,R)$, where $G$ acts on the mapping space 
by conjugation (\ref{con-action of group}). 
In our main example, $R$ will be the Gamma-ring $DB$ 
and $G$ will be the discrete group of power series 
in one variable with leading term $x$ and with multiplication
give by substitution. 
Furthermore, $X$ will be the set of formal group laws over $B$, 
and the map FGL$(B)\to  \map_{\mathcal{GR}}(H\mZ,DB)$ takes $F$ to $F_{\ast}$.

In Construction \ref{conjugation construction} we produced a
simplicial group $UR^{\times}$, a homomorphism of simplicial groups 
$G\to UR^{\times}$ and a commutative diagram of Gamma-rings
\[\xymatrix{ {\mathbb S}[G] \ar[d] & {\mathbb S}[G] \ar@{=}[l] \ar[d] \ar[r] &
{\mathbb S}[UR^{\times}] \ar[d] \\
R^{\f} & R_1 \ar_{\sim}[r] \ar@{->>}^{\sim}[l] & R_3 }\]
in which the lower horizontal maps are stable equivalences 
between stably fibrant models of $R$. Furthermore, the induced map 
from the simplicial group $UR^{\times}$ to the invertible components 
of the underlying monoid of $R_3$ is a weak equivalence. 
The simplicial group $G$ acts by conjugation on $R$, $R_1$ and $R_3$,
hence on the spaces of Gamma-ring maps from any other Gamma-ring 
into $R$, $R_1$ and $R_3$. 
If $S$ is a cofibrant Gamma-ring, then the lower horizontal maps 
induce weak equivalences
\[ \map_{\GR}(S,R^{\f}) \ \varl{1cm}{\sim} \ 
\map_{\GR}(S,R_1) \ \varr{1cm}{\sim} 
\map_{\GR}(S,R_3) \]
which are $G$--equivariant. Furthermore the action of $G$ on
$\map_{\GR}(S,R_3)$ extends to an action 
of the group $UR^{\times}$.
So we have extended, up to weak equivalence, the action of the group $G$ 
to the action of a simplicial group weakly equivalent to the homotopy units 
of $R$. 

If we now choose a cofibrant replacement $H\mZ\rscript{c}\to H\mZ$ 
in the model category of Gamma-rings, then the space 
$\map_{\mathcal{GR}}(H\mZ\rscript{c},R_3)$ is a model 
for the homotopy invariant space of Gamma-ring maps. 
Furthermore, the simplicial group $UR^{\times}$ acts on this space 
by conjugation, hence so does its subgroup $UR^{\times}_1$,
the connected component of the identity element. 
We abbreviate the space $\map_{\mathcal{GR}}(H\mZ\rscript{c},R_3)$ to 
$\Rg(H\mZ,R)$ and write $\Rg(H\mZ,R)/\mbox{conj.}$ 
for the homotopy orbit space of the conjugation action of 
the connected simplicial group $UR^{\times}_1$,
\[ \Rg(H\mZ,R)/\mbox{conj.} \ = \  
\map_{\mathcal{GR}}(H\mZ\rscript{c},R_3)_{hUR^{\times}_1} \ . \]
We now use the $G$--space $X$ and the equivariant map 
$X\to \map_{\mathcal{GR}}(H\mZ,R)$ to construct a weak map 
from the homotopy orbit space $X_{hG}$ to the homotopy orbit space 
$\Rg(H\mZ,R)/\mbox{conj.}$. 
We let $\widetilde{X}$ denote the pullback of the diagram
\[\xymatrix{ & \map_{\mathcal GR}(H\mZ\rscript{c},R_1) \ar@{->>}^{\sim}[d] \\
X \ar[r] & \map_{\mathcal GR}(H\mZ\rscript{c},R^{\f}) }\]
Since all spaces in the diagram have an action by the group 
$G$ and all maps are equivariant, the group $G$ 
acts on the space  $\widetilde{X}$. Since the map $R_1\to R^{\f}$ 
is an acyclic fibration of Gamma-rings, the induced map 
on homomorphism spaces is an acyclic fibration, 
hence so is the map $\widetilde{X}\to X$.
The stable equivalence of Gamma-rings $R_1\to R_3$ induces a 
$G$--equivariant map of homomorphism spaces 
$\map_{\mathcal{GR}}(H\mZ\rscript{c},R_1)\to
\map_{\mathcal{GR}}(H\mZ\rscript{c},R_3)$. 

The conjugation action of $G$ on the target space extends 
to an action of the simplicial group $UR^{\times}$. 
By our assumption on the homomorphism $\psi\co G\to R(1^+)$ the image
of $G$ lands in the identity component $UR^{\times}_1$.
We denote by $\kappa$ the map induced on homotopy orbit spaces
\begin{eqnarray*} \kappa \,\co\,\widetilde{X}_{hG} & \varr{1.3cm}{} & 
\map_{\mathcal{GR}}(H\mZ\rscript{c},R_1)_{hG} \\ 
& \varr{1.3cm}{} & 
\map_{\mathcal{GR}}(H\mZ\rscript{c},R_3)_{hUR^{\times}_1} \ = \  
\Rg(H\mZ,R)/\mbox{conj.}
\end{eqnarray*}
So altogether we have obtained a weak map of homotopy orbit spaces
\[ X_{hG} \ \varl{1cm}{\sim} \ 
 \widetilde{X}_{hG} \ \varr{1cm}{\kappa} \Rg(H\mZ,R)/\mbox{conj.} \]
\end{construction}

Now we return to the main example and apply 
Construction \ref{construction weak map} to the Gamma-ring $DB$. 
In this case $G$ is the discrete group $\Phi(B)$ of power series $\varphi(x)$
in one variable over $B$ with leading term $x$, 
with composition (substitution) of power series  as the group structure. 
This group acts on the set of formal group laws over $B$ via
\[ F^{\varphi}(x,y) \ = \ 
\varphi(F(\varphi^{\mbox{-}1}(x),\varphi^{\mbox{-}1}(y))) \ .  \] 
In fact, $F^{\varphi}$ is defined so that $\varphi\co F\to F^{\varphi}$ 
is a strict isomorphism of formal group laws. 
The homomorphism $\Phi(B)\to DB(1^+)\iso x\cdot B[\![x]\!]$ is the
inclusion.
Because of the equality
\[ (F^{\varphi})_* \ = \ \varphi \cdot F_* \cdot \varphi^{\mbox{-}1} \]
as maps of Gamma-rings $H\mZ\to DB$, the assignment $F\mapsto F_{\ast}$ 
from the set FGL$(B)$ of formal group laws to the set of Gamma-ring 
maps from $H\mZ$ to $DB$ is $\Phi(B)$--equivariant. 
The homotopy orbit space FGL$(B)_{h\Phi(B)}$ is isomorphic 
to the classifying space of the groupoid of formal group laws 
and strict isomorphisms, which we denote  by $\mathcal{FGL}\rscript{str}(B)$.
So Construction \ref{construction weak map} yields maps
\[\xymatrix{ \mathcal{FGL}\rscript{str}(B) = 
\mbox{FGL}\rscript{str}(B)_{h\Phi(B)} 
& \widetilde{\mbox{FGL}\rscript{str}(B)}_{h\Phi(B)} \ar[l]_-{\sim}  
\ar[d]^-{\kappa} \\ 
& \Rg(H\mZ,DB)/\mbox{conj.} }\]
where the upper map is a weak equivalence.
We use the notation $\widetilde{\mathcal{FGL}}\rscript{str}(B)$ 
for the homotopy orbit space 
$\widetilde{{\mbox{FGL}}\rscript{str}(B)}_{h\Phi(B)}$. 
The following theorem is our main result:

\begin{theorem} \label{main} The map 
\[ \kappa \,\co\,\widetilde{\mathcal{FGL}}\rscript{\em str}(B)\ \to\ 
\Rg(H\mZ,DB)/\mbox{\em conj.} \]
is a weak equivalence.
\end{theorem}

The proof of Theorem \ref{main} occupies the rest of this paper.
Since the homotopy orbit space construction defining
$\Rg(H\mZ,DB)/\mbox{conj.}$ involves a connected simplicial group,
the quotient map 
\[ \Rg(H\mZ,DB)\ \to \  \Rg(H\mZ,DB)/\mbox{conj.} \]
induces a bijection of path components. 
The set $[H\mZ,DB]_{\text{Gamma-rings}}$ is canonically isomorphic
to the components of the mapping space $\Rg(H\mZ,DB)$,
so Theorem \ref{homotopy classes} is just the bijection of path components
induced by the weak equivalence $\kappa$ of Theorem \ref{main}.

In addition to the homotopy classes of Gamma-ring maps, Theorem \ref{main} 
allows us to identify the higher homotopy groups of the space  $\Rg(H\mZ,DB)$ 
of Gamma-ring maps.  
Since the group  $\pi_1(DB)^{\times} \iso \pi_1 DB$ is trivial by
\ref{DB summary} (b), the simplicial group $(DB)^{\times}_1$ is
1--connected, so the quotient map 
$\Rg(H\mZ,DB)\to \Rg(H\mZ,DB)/\mbox{conj.}$
induces an equivalence of fundamental groupoids.
Together with Theorem \ref{main} this  implies that the fundamental groupoid
of the space $\Rg(H\mZ,DB)$ is equivalent to the groupoid
$\mathcal{FGL}\rscript{str}(B)$.
In particular this yields isomorphisms
\[ \pi_1 \Rg(H\mZ,DB;F_{\ast}) \ \iso \ \Aut\rscript{strict}(F)  \] 
between the fundamental group at the basepoint $F_*\in\Rg(H\mZ,DB)$ 
and the strict automorphism group of $F$.
Since the homotopy orbit space of the conjugation action is weakly equivalent 
to the classifying space of a groupoid, its homotopy groups
are trivial above dimension 1; so for every $F$ the action map
\[ U(DB)^{\times}_1 \ \to  \  \Rg(H\mZ,DB) \ ,
\quad u \, \longmapsto \ u \cdot F_{\ast} \cdot u^{-1} \]
induces isomorphisms of homotopy groups 
$\pi_n DB \iso \pi_n \, \Rg(H\mZ,DB;F_{\ast})$ for $n\geq 2$.

\begin{example} For algebras over the rational numbers, the map
$\kappa$ of Theorem \ref{main} is trivially an equivalence.
Indeed, if $B$ is a $\mQ$--algebra, and $F$ a 1--dimensional and commutative
formal group law over $B$, then there is a strict isomorphism 
(called the {\em logarithm} of $F$) between $F$ 
and the additive formal group law \cite[III.1 Cor.\ 1]{Froehlich}. 
Moreover, $F$ has no non-trivial strict automorphisms, 
so the classifying space of the groupoid $\mathcal{FGL}\rscript{str}(B)$
is weakly contractible.

On the other hand, the Gamma-ring $DB$ is now stably equivalent to the
Eilenberg--MacLane Gamma-ring $HB$ by Theorem \ref{DB summary} (b).
So both the space $\Rg(H\mZ,DB)$ and the unit component of $(DB)^{\times}$ 
are weakly contractible, hence so is
the homotopy orbit space  $\Rg(H\mZ,DB)/\mbox{conj.}$
\end{example}

\begin{remark} 
Instead of taking homotopy orbits with respect to the connected simplicial
group $U(DB)^{\times}_1$ one can divide out the conjugation action of
the entire homotopy units $U(DB)^{\times}$ on the space 
$\Rg(H\mZ,DB)$. The resulting orbit space receives a (weak) map from the 
groupoid of formal group laws and all (ie, not necessarily strict)
isomorphisms.
The same proof as for Theorem \ref{main} shows that that map
\[ \widetilde{\mathcal{FGL}}(B)\ \to\ 
\Rg(H\mZ,DB)_{hU(DB)^{\times}} \]
is a weak equivalence.
\end{remark}

\section{A filtration of $DB$} \label{filtration}

The Gamma-ring $DB$ has a natural filtration arising from powers 
of the augmentation ideal of the power series rings. 
There are truncated versions $D_kB$ 
of the Gamma-ring $DB$ and analogues of the map $\kappa$ 
of Theorem \ref{main} for every $k$. 
In this section we reduce the proof of Theorem \ref{main} 
to the analogous statement about the stages of the filtration, 
see Theorem \ref{k-bud comparison}.

For $k\geq 1$ we denote by $D_kB$ the truncated version 
of the Gamma-ring $DB$ obtained by dividing out all power series in the 
$(k+1)$-st power of the augmentation ideal. So as a $\Gamma$--space,
\[ D_kB(K) \, = \,  \mbox{kernel}\,(\,\widetilde{B}[\![K]\!]/I^{k+1} \to 
\widetilde{B}[\![\ast]\!] = B \, ) \]
where $I=K\cdot\widetilde{B}[\![K]\!]$ is the ideal 
of $\widetilde{B}[\![K]\!]$ consisting of power series without
constant term. 
The unit map again comes from the inclusion of generators 
and the multiplication is induced by substitution of (truncated) power series,
similar to the definition for $DB$. 
In other words: $D_kB$ has a unique Gamma-ring structure for which 
the natural projection map $DB\to D_kB$ is a homomorphism of Gamma-rings. 
Note that $D_1B$ is isomorphic to the Eilenberg--MacLane 
Gamma-ring $HB$. There are maps of Gamma-rings 
$DB\to D_kB$ and $D_kB\to D_{k-1}B$ induced by truncation of power series.

Now we apply Construction \ref{construction weak map} 
to the Gamma-ring $D_kB$. 
The group we work relative to is $\Phi_k(B)$, 
the quotient of the group $\Phi(B)$ of power series over $B$ with
leading term $x$ by the normal subgroup of power series which are congruent 
to the power series $x$ modulo $x^{k+1}$. 
The group $\Phi_k(B)$ injects into the monoid $D_kB(1^+)$, 
so Construction \ref{conjugation construction} 
provides a diagram of Gamma-rings
\[\xymatrix{ {\mathbb S}[\Phi_k(B)] \ar[d] & 
{\mathbb S}[\Phi_k(B)] \ar@{=}[l] \ar[d] \ar[r] & 
{\mathbb S}[U(D_kB)^{\times}] \ar[d] \\
(D_kB)^{\f} & (D_kB)_1 \ar_{\sim}[r] \ar@{->>}^{\sim}[l] & (D_kB)_3 }\]
in which the lower horizontal maps are stable equivalences 
between fibrant Gamma-rings. Furthermore the induced map 
from the simplicial group\break $U(D_kB)^{\times}$ to the invertible components 
of the underlying monoid of $(D_kB)_3$ is a weak equivalence. 
Hence the simplicial group $U(D_kB)^{\times}$ acts by conjugation on the space 
\[ \Rg(H\mZ,D_kB) \, = \, \map_{\mathcal GR}(H\mZ\rscript{c}, (D_kB)_3) \]
extending the action of the group $\Phi_k(B)$. 
As in~\eqref{construction weak map} we denote by
$$\Rg(H\mZ,D_kB)/\mbox{conj.}$$ 
the homotopy orbit space of  
$\Rg(H\mZ,D_kB)$ by the conjugation action of $U(D_kB)^{\times}_1$,
the identity component of $U(D_kB)^{\times}$.

Since the Constructions \ref{conjugation construction} 
and \ref{construction weak map} can be made functorial, 
truncations induce compatible maps
\begin{eqnarray*} \Rg(H\mZ,DB) & \to & \Rg(H\mZ,D_kB) \\
\mbox{and \qquad} \Rg(H\mZ,D_kB)  & \to & \Rg(H\mZ,D_{k-1}B)  
\end{eqnarray*}
and similarly for the orbit spaces by the conjugation actions.

\begin{lemma} \label{ring holim} The maps 
\[ \Rg(H\mZ,DB) \ \varr{1cm} \ \holim_k \ \Rg(H\mZ,D_kB) \] 
and
\[ \Rg(H\mZ,DB)/\mbox{\em conj.} \ \varr{1cm} \ 
\holim_k \ (\Rg(H\mZ,D_kB)/\mbox{\em conj.}) \]
induced by truncation are weak equivalences.
\end{lemma}
\begin{proof} We apply the homotopy limit construction \cite{Bous-Kan} 
objectwise to Gamma-rings to obtain 
a construction of homotopy limits for Gamma-rings. 
As we explained in the proof of Theorem \ref{DB summary} (a),
the homotopy fibre of the projection $DB\to D_kB$ is the product 
of certain $\Gamma$--spaces $(B\tensor {\mathcal S}^m)^!$ for $m>k$, each of which
is $(2k-1)$--connected. Hence the map
\[ (DB)_3 \ \to \ \holim_k \ (D_k B)_3 \] 
is a stable equivalence of Gamma-rings.
So after taking homomorphism spaces from $H\mZ\rscript{c}$ 
we obtain a weak equivalence of spaces
\begin{eqnarray*}  \map_{\mathcal{GR}}(H\mZ\rscript{c},(DB)_3) &
\varr{1cm}{\sim} & \map_{\mathcal{GR}}(H\mZ\rscript{c},\holim_k 
\, (D_kB)_3) \\
& &  \iso \ \holim_k \, \Rg(H\mZ,D_kB) \  
\end{eqnarray*}
which proves the first statement. 
Similarly, the induced map of homotopy units
\[ (DB)^{\times}  \to \ \holim_k \ (D_k B)^{\times}  \] 
is a weak equivalence, and so also the homotopy orbit space
$\Rg(H\mZ,DB)/\mbox{conj.}$ is the homotopy inverse limit of the
truncated versions.
\end{proof}

On the formal group law side, the classifying space 
$\mathcal{FGL}\rscript{str}(B)$ 
can also be expressed as a homotopy limit of suitable truncated versions. 
We denote by ${\B}ud^k_B$ the classifying space of the groupoid of $k$--buds 
\cite[Def.\ 2.1]{Lazard:formels} (also called $k$--jets) 
of formal group laws over $B$ and $k$--buds of strict isomorphisms. 
Truncation induces maps
\[ \mathcal{FGL}\rscript{str}(B) \ \to {\B}ud^k_B 
\mbox{\qquad and \qquad}  {\B}ud^k_B \ \to {\B}ud^{k-1}_B \ . \]
The classifying spaces $\mathcal{FGL}\rscript{str}(B)$ and ${\B}ud^k_B$
are isomorphic to the homotopy orbit spaces $\mbox{FGL}(B)_{h\Phi(B)}$ and 
$(\mbox{Bud}^k_B)_{h\Phi_k(B)}$ respectively. 
Since the group $\Phi(B)$ is the inverse limit of the groups 
$\Phi_k(B)$ and the set of formal group laws is the inverse limit 
of the sets of $k$--buds, we have 

\begin{lemma} \label{FGL holim} The map
\[  \mathcal{FGL}\rscript{str}(B) \ \to \ 
\mbox{\em holim}_k \ {\B}ud^k_B \]
induced by truncation on classifying spaces is a weak equivalence.
\end{lemma}
\begin{proof} This an instance of a general fact 
about homotopy orbits of groups acting on sets, alias translation categories. 
Suppose $\{G_k\to G_{k-1}\}_{k\geq 1}$ is a sequence of surjective
group homomorphisms,
$\{X_k\to X_{k-1}\}_{k\geq 1}$ a tower of sets, and suppose that $G_k$
acts on $X_k$ in such a way that the map $X_k\to X_{k-1}$ is $G_k$--equivariant.
Then the inverse limit $G=\lim_k G_k$ of groups acts on the inverse limit
$X=\lim_k X_k$ of sets and the canonical map
\[ X_{hG} \ \to \ \holim_k \ (X_k)_{hG_k} \]
is a weak equivalence. This follows from the homotopy fibre sequences
\[ X_k \ \to \ (X_k)_{hG_k} \ \to \ BG_k \]
by passage to homotopy inverse limit together with the fact that
the natural maps
\[ X = {\lim}_k \ X_k  \ \to \ \holim_k \ X_k 
\mbox{\qquad and \qquad}
 BG \ \to \ \holim_k \ BG_k  \]
are weak equivalences (we note that the path components of the homotopy inverse
limit of the classifying spaces $BG_k$ are in bijective correspondence with 
the set  $\lim^1_kG_k$; since we consider surjective group homomorphisms,
this $\lim^1$--term is trivial and the map from $BG$ to the homotopy inverse
limit is indeed a weak equivalence).
In our example $X_k$ is the set of $k$--buds
of formal group laws and $G_k$ is the group $\Phi_k(B)$ of $k$--buds
of power series conjugating the $k$--buds of formal group laws.
\end{proof}

Every $k$--bud of formal group law $F$ gives rise to a map of Gamma-rings 
$H\mZ\to D_kB$ in the same fashion as genuine formal group laws give maps 
to $DB$ (\ref{Gamma-ring map}). 
The group $\Phi_k(B)$ of truncated invertible power series 
conjugates the set of $k$--buds and the map 
Bud$_B^k\to \map_{\mathcal{GR}}(H\mZ,D_kB)$ is  $\Phi_k(B)$--equivariant. 
If we carry out Construction \ref{construction weak map} with $R=D_kB$ 
and with the group $\Phi_k(B)$ conjugating the $k$--buds 
of formal group laws, we obtain a map of homotopy orbit spaces
\begin{eqnarray*} \kappa_k \,\co\, \widetilde{{\B}ud}^k_B \, = \, 
(\widetilde{\mbox{Bud}_B^k})_{h\Phi_k(B)} & \to & 
\map_{\mathcal{GR}}(H\mZ\rscript{c},(D_kB)_3)_{hU(D_kB)_1^{\times}} \\
& = & \Rg(H\mZ,D_kB)/\mbox{conj.}  
\end{eqnarray*}
The constructions are natural, so we end up with a commutative diagram
\begin{equation} \label{homotopy limit square}
\xymatrix@C=20mm{ \widetilde{\mathcal{FGL}}\rscript{str}(B) \ar_{\sim}[d] 
\ar^-{\kappa}[r] & \Rg(H\mZ,DB)/\mbox{conj.} \ar^{\sim}[d] \\
\holim_k \ \widetilde{{\B}ud}^k_B \ar_-{\kappa_k}[r] & 
\holim_k \ \Rg(H\mZ,D_kB)/\mbox{conj.} } \end{equation}
in which the vertical maps are weak equivalences 
by Lemmas \ref{ring holim} and  \ref{FGL holim}.
So Theorem \ref{main} follows once we have shown

\begin{theorem} \label{k-bud comparison} For all $k\geq 1$ the map
\[ \kappa_k \,\co\,   \widetilde{{\B}ud}^k_B \ \to \ 
 \Rg(H\mZ,D_kB)/\mbox{\em conj.}  \]
is a weak equivalence.
\end{theorem}

\begin{remark} \label{constant unless prime power}
If $k$ is not a prime power, then the functor $B\tensor {\mathcal S}^k$
is a retract of a diagonalizable functor (Lemma \ref{Dold-Puppe lemma}), 
and so the $\Gamma$--space underlying $(B\tensor {\mathcal S}^k)^!$ 
is stably contractible (\ref{DP and Q} (e)). 
Hence the reduction map $D_kB\to D_{k-1}B$ is a stable equivalence 
of Gamma-rings and the induced map on homotopy units 
$(D_kB)^{\times}\to (D_{k-1}B)^{\times}$
is a weak equivalence of simplicial monoids.
Thus the map of homotopy orbit spaces
\[ \Rg(H\mZ,D_kB)/\mbox{conj.}  \ \varr{1cm}{\text{red.}} 
\ \Rg(H\mZ,D_{k-1}B)/\mbox{conj.} \]
is a weak equivalence.
By Remark \ref{Lazard's thm} below, the reduction functor
${{\B}ud}^k_B\to{{\B}ud}^{k-1}_B$ is an equivalence of categories.
Hence in the inductive step nothing happens unless $k=p^h$ is a prime power. 
However it is convenient to make this case distinction only
at the very end (see Step 1 in the proof of Theorem \ref{deduction theorem}).
\end{remark}

\section{Some singular extensions of Gamma-rings}
\label{singular extensions}

In this section we start the inductive proof of Theorem \ref{k-bud comparison}.
We exploit that the truncation maps $D_kB\to D_{k-1}B$ 
are ``singular extensions'' of Gamma-rings.
This lets us reduce the problem to a comparison
of the derivation space of the Gamma-ring $H\mZ$ 
with coefficients in the ``kernel'' $(B\tensor {\mathcal S}^k)^!$ 
of the extension to the groupoid of symmetric 2--cocycles.

The case $k=1$ of Theorem \ref{k-bud comparison} is straightforward. 
There is only one 1--bud of formal group law, and the only
1--bud of strict automorphism is the identity.
On the other hand, $D_1B$ is isomorphic 
to the Eilenberg--MacLane Gamma-ring $HB$, so both the space
$\Rg(H\mZ,D_1B)$ and the unit component  $(D_1B)^{\times}_1$ 
of the homotopy units are weakly contractible. 
Hence source and target of the map
\[ \kappa_1 \,\co\,   \widetilde{{\B}ud}^1_B \ \to \ 
 \Rg(H\mZ,D_1B)/\mbox{conj.}  \]
are weakly contractible and $\kappa_1$ is a weak equivalence.

\begin{numbered paragraph} \label{B cocycles} 
{\bf Symmetric 2--cocycles}\qua
For the inductive step we recall how the difference
between $k$--buds  and $(k-1)$--buds of formal group laws 
is controlled by symmetric 2--cocycles. 
A {\em symmetric 2--cocycle} of degree $k$ with values in $B$ 
is a homogeneous polynomial $c(x,y)\in B[x,y]$ of degree $k$ 
which satisfies the relations
\[ c(x,y) \ = \ c(y,x) \mbox{\quad and \quad} 
c(x,y) \ + \ c(x+y,z)  \ = \  c(x,y+z) \ + \ c(y,z)\ . \]
If $F$ is any $k$--bud of formal group law over $B$ 
and $c$ is a $k$--homogenous 2--cocycle, then the truncated power series 
$F(x,y)+c(x,y)$ is another $k$--bud of formal group law with the same 
$(k-1)$--bud as $F$. Conversely, if $F$ and $F'$ are two $k$--buds with the same
reduction modulo the $k$--th powers of the augmentation ideal, then $c=F-F'$ 
is a $k$--homogenous 2--cocycle.
The proof of this is straightforward, compare
\cite[Sec.\ II]{Lazard:formels} or \cite[III.1 Lemma 2]{Froehlich}.

We define $\Z(B\tensor {\mathcal S}^k)$, the {\em groupoid of symmetric 2--cocycles} 
of degree $k$,
as the category whose objects are the  symmetric 2--cocycles of degree $k$ 
over $B$. The set of morphisms from a cocycle $c$ to a cocycle $c'$
consists of those $b\in B$ satisfying 
$c'=c+ b\cdot\left[  x^k + y^k - (x+y)^k \right]$; composition is given
by addition in $B$.
\end{numbered paragraph}

Suppose $F$ is a $k$--bud of formal group law.
Then we can define a functor 
\[ F+- \,\co\,  \Z(B\tensor {\mathcal S}^k) \ \to \ {{\B}ud}^k_B  \]
on objects by $(F+-)(c) = F+c$ and on morphisms by 
$(F+-)(b) = x +b\cdot x^k$. The relation
\[ F^{x+bx^k}(y,\bar y) \ \equiv \ 
F(y,\bar y) + b\cdot\left[ y^k + \bar y^k -  (y+\bar y)^k\right] 
\mod (y,\bar y)^{k+1} \] 
shows that if $b\co c\to c'$ is a morphism of cocycles, 
then the power series $x+b\cdot x^k$ is
indeed an isomorphism from $F+c$ to $F+c'$.

\begin{lemma} \label{FGL fibre}
For every $k$--bud of formal group law $F$, the functor
\[ F+- \,\co\,\Z(B\tensor {\mathcal S}^k) \ \to \ {{\B}ud}^k_B \]
induces a weak equivalence from the classifying space of
the groupoid \mbox{$\Z(B\tensor {\mathcal S}^k)$} to the homotopy fibre
of the truncation map ${{\B}ud}^k_B \to {{\B}ud}^{k-1}_B$ 
over the basepoint $F$.
\end{lemma}
\begin{proof} This  again is an instance of a general fact 
about homotopy orbits of groups acting on sets, alias translation categories. 
Suppose $G\to \bar G$ is an epimorphism of groups with kernel $K$. 
Moreover, let $X$ be a $G$--set, $\bar X$ a $\bar G$--set and
$\pi\co X\to \bar X$ a $G$--equivariant map. Then for every point
$x\in \bar X$ the kernel $K$ acts on the preimage $\pi^{-1}(x)$. 
In this situation the sequence of homotopy orbit spaces
\[ \pi^{-1}(x)_{hK} \ \to \ X_{hG} \ \to 
\ \bar X_{h\bar G} \] 
is a homotopy fibre sequence over the point $x$.

In the situation at hand the epimorphism is the truncation
$\Phi_k(B)\to \Phi_{k-1}(B)$, whose kernel is isomorphic to the
additive group of $B$ via $b\longmapsto x+b\cdot x^k$.
The groups conjugate the $k$--buds respectively $(k-1)$--buds 
of formal group laws.
For every choice of $k$--bud of formal group law $F$, 
the map $F+-$ is an isomorphism, equivariant for 
$B\iso\mbox{kernel}\co \Phi_k(B)\to \Phi_{k-1}(B)$, from the
symmetric 2--cocycles to the $k$--buds which have the same $(k-1)$--bud as $F$.
Hence the lemma follows.
\end{proof}

\begin{remark} \label{Lazard's thm} 
The symmetric 2--cocycles have been identified by Lazard 
\cite[II Lemme 3]{Lazard:formels}, see also \cite[III.1 Thm.\ 1a]{Froehlich}. 
There is a universal integral symmetric 2--cocycle $c_k\in \mZ[x,y]$ 
of degree $k$ given by
\[ c_k(x,y) \ = \ \frac{1}{d_k} \,  \left[  x^k + y^k - (x+y)^k \right] \]
and the degree $k$ symmetric 2--cocycles over $B$ are precisely the multiples 
$b\cdot c_k$ for $b\in B$. 
Here $d_k$ is the greatest common divisor of the binomial coefficients  
$k \choose i$ for $1\leq i\leq k-1$, which evaluates to
\[ d_k \ = \ \left\lbrace \begin{array}{cl} p & \mbox{if } k=p^h 
\mbox{ for a prime $p$ and $h\geq 1$} \\
1 & \mbox{else.} \end{array} \right. \]

Hence the classifying space of $\Z(B\tensor {\mathcal S}^k)$ can be identified as follows:
if $k$ is not a prime power, then the classifying space is weakly contractible
and the reduction functor
${{\B}ud}^k_B\to{{\B}ud}^{k-1}_B$ is an equivalence of categories.
If $k=p^h$ for some prime $p$, then the groupoid
$\Z(B\tensor {\mathcal S}^{p^h})$ is isomorphic to the translation category of the
action of $B$ on itself given by
\[ (b,x) \ \longmapsto \ p\cdot b \, + \, x \ .\]
Hence the components of the classifying space of $\Z(B\tensor {\mathcal S}^{p^h})$
are in bijective correspondence with the set  $B/pB$ and the fundamental
group at each basepoint is isomorphic to the group of those $b\in B$
such that $pb=0$.
However, in the rest of this paper we will not use this explicit 
knowledge about the symmetric 2--cocycles. 
\end{remark}

Now we identify the difference between the spaces of Gamma-ring maps 
from $H\mZ$ to $D_kB$ and to $D_{k-1}B$, 
ie, we study the homotopy fibers of the reduction map
\[ \Rg(H\mZ,D_kB)/\mbox{conj.} \ \to\ \Rg(H\mZ,D_{k-1}B)/\mbox{conj.} \]
For this purpose we consider the Gamma-ring 
$H\mathbb Z\times(B\tensor {\mathcal S}^k)^!$. 
We present the construction in a more general context, since we need
it later.

\begin{numbered paragraph} \label{split extension}
{\bf Split singular extensions}\qua
Let $G$ be a functor from the category of finitely generated
free abelian groups to the category of all abelian groups,
and suppose that $G$ is {\em reduced} in the sense that $G(0)\iso 0$.
We define a Gamma-ring $H\mZ\times G^!$, the {\em split singular extension}
of $H\mZ$ by the bimodule $G^!$.

First there is an $H\mZ$--bimodule $G^!$ associated to the functor $G$;
the notation is taken from \cite[Ex.\ 2.6]{Pira-Wald}, 
where the construction first appeared.
As a $\Gamma$--space, $G^!$ is the composite 
\[ \Gamma\rscript{op} \ \varr{1cm}{\widetilde{\mZ}} \ 
\mbox{(f.\ g.\ free ab.\ groups)} \ \varr{1cm}{G} \ \Ab \   
\varr{1cm}{\Phi} \ \mbox{(pt.\ sets)} \] 
of the reduced free functor $\widetilde\mZ$,
followed by $G$ and the forgetful functor $\Phi$ 
from abelian groups to pointed sets.
The $\Gamma$--space $G^!$ has the structure of an $H\mZ$--bimodule via
the composite
\[ \xymatrix@C=15mm{ H\mZ \sm G^! \sm  H\mZ  \ar[r]^-{\text{assembly}} &  
H\mZ \circ G^! \circ H\mZ \ar[r] & G^! \ .} \]
The first map is an instance of the assembly map (\ref{assembly map}) 
from the smash product to the composition product of $\Gamma$--spaces;
the second map is induced by evaluation maps (this uses that
the original functor $G$ was defined for and takes values in abelian groups).

The product $H\mZ \times G^!$  becomes a Gamma-ring as follows.
The unit map of  $H\mZ \times G^!$ is the composition of the unit
$\eta\co \mS\to H\mZ$ with the inclusion  $H\mZ\to H\mZ \times G^!$.
The multiplication map is the composite
\begin{eqnarray*} (H\mZ \times G^!) \, \sm \, (H\mZ \times G^!) 
& \varr{2cm}{} &  (H\mZ \sm H\mZ) \times(H\mZ \sm G^!) \times
(G^!\sm H\mZ) \\
& \varr{2cm}{(\mu_{H\mZ},\ l+r)} &  \qquad  H\mZ \, \times \, G^! 
\end{eqnarray*} 
where $l$ and $r$ denote the left respectively right action of $H\mZ$ on $G^!$.
\end{numbered paragraph}

In particular we can apply Construction \ref{split extension} to the
functor $B\tensor {\mathcal S}^k$ which takes 
a finitely generated free abelian group
$A$ to the tensor product of $B$ with the $k$-th symmetric power of $A$.
If $F$ is any $k$--bud of formal group law, then $F_*\co H\mZ\to D_kB$ is
a morphism of Gamma-rings.
Homogenous polynomials of degree $k$ naturally inject into
the quotient of power series by terms of degree $k+1$, which gives
a map of $\Gamma$--spaces $\mbox{Incl}\co (B\tensor {\mathcal S}^k)^!\to D_kB$. 
Moreover, their pointwise sum in $D_kB$ 
\[ F_{\ast} + \mbox{Incl.} \,\co\, H\mZ \times(B\tensor {\mathcal S}^k)^! \ 
\to \ D_kB \]
is again a map of Gamma-rings.

\begin{lemma} \label{fibre square} The commutative square
\[\xymatrix@C=20mm{ H\mZ \times (B\tensor {\mathcal S}^k)^! \ar_{\text{\em proj.}}[d]
\ar^-{F_{\ast} + \text{\em Incl.}}[r] & D_kB \ar[d] \\
H\mathbb Z \ar_-{F_{\ast}}[r] & D_{k-1}B }\]
is a homotopy fibre square of Gamma-rings. 
\end{lemma}
\begin{proof} It suffices to show that the underlying square 
of $\Gamma$--spaces is homotopy cartesian. 
As a $\Gamma$--space $DB$ splits as a product 
\[ D_kB \ \iso  \ D_{k-1}B\, \times \, (B\tensor {\mathcal S}^k)^!  \]
and under this isomorphism the map $F_{\ast} + \text{Incl.}$
becomes the map 
\[ F_*\times \mbox{Id} \,\co\,H\mZ \times (B\tensor {\mathcal S}^k)^!
\ \to \  D_{k-1}B \, \times \, (B\tensor {\mathcal S}^k)^! \ , \]
so the claim follows.
\end{proof}

The additive group of the ring $B$ includes into the underlying monoid 
of the Gamma-ring $H\mZ \times(B\tensor {\mathcal S}^k)^!$ via the map which sends 
$b\in B$ to the polynomial $x+b\cdot x^k$, considered as an element 
of $(H\mZ \times(B\tensor {\mathcal S}^k)^!)(1^+)\iso \mZ[x]\times B[x^k]$ 
(where $x$ is an indeterminate corresponding to the non-basepoint element 
of $1^+$). We can now apply Construction \ref{conjugation construction} 
to the Gamma-ring $H\mZ \times(B\tensor {\mathcal S}^k)^!$ relative to the additive group
of $B$. 
The construction produces a commutative diagram of Gamma-rings
\[\xymatrix@C=11mm{ \mS[B] \ar[d] & \mS[B] \ar@{=}[l] \ar[d] \ar[r] & 
\mS[U(H\mZ\times(B\tensor {\mathcal S}^k)^!)^{\times}] \ar[d] \\
(H\mZ \times(B\tensor {\mathcal S}^k)^!)^{\f} & 
(H\mZ \times(B\tensor {\mathcal S}^k)^!)_1 
\ar_{\sim}[r] \ar@{->>}^{\sim}[l] & 
(H\mZ \times(B\tensor {\mathcal S}^k)^!)_3 }\]
in which the lower horizontal maps are stable equivalences between 
fibrant Gamma-rings. The induced map from the simplicial group 
$U(H\mZ\times(B\tensor {\mathcal S}^k)^!)^{\times}$ 
to the invertible components of the underlying monoid of 
$(H\mZ \times(B\tensor {\mathcal S}^k)^!)_3$ is a weak equivalence. 

We denote by 
\[ \der(H\mZ,B\tensor {\mathcal S}^k) \ = \ 
\map_{\GR}(H\mZ\rscript{c},(H\mZ \times(B\tensor {\mathcal S}^k)^!)_3) \]
the space of Gamma-ring maps 
from $H\mZ\rscript{c}$ to $(H\mZ \times(B\tensor {\mathcal S}^k)^!)_3$ 
and refer to this space as the space of {\em derivations} of $H\mZ$ 
with coefficients in $B\tensor {\mathcal S}^k$.
The group  $U(H\mZ\times (B\tensor {\mathcal S}^k)^!)^{\times}$ acts by conjugation 
on the space of derivations and we denote by
$\der(H\mZ,B\tensor {\mathcal S}^k)/\mbox{conj.}$ the homotopy orbit space 
of the conjugation action of the identity component 
$U(H\mZ\times (B\tensor {\mathcal S}^k)^!)^{\times}_1$.
Since the square of monoids 
\[\xymatrix@R=15mm{ B \ar[r] \ar_{\begin{array}{c} \scriptstyle b\ \longmapsto \\ \scriptstyle  x\, +\, b\cdot x^k\end{array}}[d] &
(H\mZ\times(B\tensor {\mathcal S}^k)^!)(1^+) \ar^{F_{\ast} + \text{Incl.}}[d] \\
\Phi_k(B) \ar[r] &  D_kB(1^+) }\]
commutes and the constructions of Sections \ref{conjugation} are functorial 
in Gamma-rings equipped with a map from a group to the underlying monoid, 
the map 
\[ F_{\ast} + \mbox{Incl.} \,\co\,  H\mZ \times (B\tensor {\mathcal S}^k)^! 
\ \to \  D_kB \] 
induces maps between the respective derived spaces of Gamma-ring maps
from $H\mZ$ and their homotopy orbit spaces.

\begin{lemma} \label{ring fibre} The two maps 
\begin{eqnarray*}\der(H\mZ,B\tensor {\mathcal S}^k) \quad &
\varr{1.6cm}{F_{\ast} + \text{\em Incl.}} & \Rg(H\mZ,D_kB)
\qquad \mbox{and} \\
\der(H\mZ,B\tensor {\mathcal S}^k)/\mbox{\em conj.} &
\varr{1.6cm}{F_{\ast} + \text{\em Incl.}} & 
\Rg(H\mZ,D_kB)/\mbox{\em conj.}  
\end{eqnarray*}
induce weak equivalences between the derivation space, respectively
its homotopy orbit space, and the respective homotopy fibres of the
truncation maps 
\begin{eqnarray*} \Rg(H\mZ,D_kB)\quad & \to & \Rg(H\mZ,D_{k-1}B) 
\mbox{\qquad and \qquad} \\
\Rg(H\mZ,D_kB)/\mbox{\em conj.} & \to & 
\Rg(H\mZ,D_{k-1}B)/\mbox{\em conj.} 
\end{eqnarray*}
over the basepoint $F_*$.
\end{lemma} 
\begin{proof} The first statement is a direct consequence of the fact 
that the square of Lemma \ref{fibre square} is homotopy cartesian. 
The second follows from the first since the sequence of simplicial groups
\[ U(H\mZ\times (B\tensor {\mathcal S}^k)^!)^{\times}_1 \ \varr{1cm} \ 
U(D_kB)_1^{\times} \ \to \ U(D_{k-1}B)_1^{\times} \]
is also a homotopy fibre sequence, again because of Lemma \ref{fibre square}.
\end{proof}

We denote by $\widetilde{\mbox{Z}}^2_s(B\tensor {\mathcal S}^k)$ 
the pullback of the diagram
\[\xymatrix{ & 
\map_{\mathcal{GR}}(H\mZ\rscript{c},
(H\mZ \times (B\tensor {\mathcal S}^k)^!)_1)
\ar@{->>}^{\sim}[d] \\
\mbox{Z}^2_s(B\tensor {\mathcal S}^k) \ar[r] &
\map_{\mathcal{GR}}(H\mZ\rscript{c},
(H\mZ \times (B\tensor {\mathcal S}^k)^!)^{\f}) }\] 
All maps are equivariant with respect to the action of the additive group 
of $B$, so this group acts on  
$\widetilde{\mbox{Z}}^2_s(B\tensor {\mathcal S}^k)$. 
Furthermore the weak equivalence 
\[ \map_{\mathcal{GR}}(H\mZ\rscript{c},
(H\mZ \times (B\tensor {\mathcal S}^k))_1) \ 
\varr{1cm}{\sim} \
\map_{\mathcal{GR}}(H\mZ\rscript{c},
(H\mZ \times (B\tensor {\mathcal S}^k))_3) \]
is $B$--equivariant, and on the target the action extends to an action 
by the simplicial group $U(H\mZ\times (B\tensor {\mathcal S}^k))^{\times}_1$. 
So we get an induced map on homotopy orbits
\begin{align*} \kappa_{B\tensor {\mathcal S}^k} \,\co\, 
\widetilde{\mbox{Z}}^2_s(B\tensor {\mathcal S}^k)_{hB}\ \to\ 
&\map_{\mathcal{GR}}(H\mZ\rscript{c},(H\mZ \times (B\tensor {\mathcal S}^k)^!)_1)_{hB} 
\ \to  \\
\map&_{\mathcal{GR}}(H\mZ\rscript{c},
(H\mZ \times (B\tensor {\mathcal S}^k)^!)_3)_{h U(H\mZ\times (B\tensor {\mathcal S}^k)^!)^{\times}_1}  \\
= \ &\der(H\mZ,B\tensor {\mathcal S}^k)/\mbox{conj.} \nonumber
\end{align*}
Note that the homotopy orbit space of the action of $B$ on the set 
Z$^2_s(B\tensor {\mathcal S}^k)$ of symmetric 2--cocycles is isomorphic 
to the classifying space of the groupoid $\Z(B\tensor {\mathcal S}^k)$ (\ref{B cocycles});
hence we use the notation $\widetilde{\Z}(B\tensor {\mathcal S}^k)$
for the weakly equivalent homotopy orbit space
$\widetilde{\mbox{Z}}^2_s(B\tensor {\mathcal S}^k)_{hB}$.

Now we can reduce the inductive step of Theorem \ref{k-bud comparison}
to a statement about the map  $\kappa_{B\tensor {\mathcal S}^k}$.
We assume inductively that the map  
\[ \kappa_{k-1} \,\co\,   \widetilde{{\B}ud}^{k-1}_B \ \to \ 
 \Rg(H\mZ,D_{k-1}B)/\mbox{conj.} \]
is a weak equivalence.
This guarantees in particular that the $(k-1)$--buds 
of formal group laws account for all components of the target space. 
Since every $(k-1)$--bud of formal group law extends to a $k$--bud, 
the reduction map 
\[  \Rg(H\mZ,D_kB)/\mbox{conj.} \ \to \ 
\Rg(H\mZ,D_{k-1}B)/\mbox{conj.}  \]
is surjective on components.
If we fix a $k$--bud of $F$ of formal group law,
then the diagram
\begin{equation} \label{singular fibre sequences}
\xymatrix@C=20mm{ \widetilde{\Z}(B\tensor {\mathcal S}^k) \ar_{ \kappa_{B\tensor {\mathcal S}^k}}[d] 
\ar^-{F+-}[r] & \widetilde{{{\B}ud}^k_B} \ar^{\kappa_k}[d] \\
\der(H\mZ,B\tensor {\mathcal S}^k)/\mbox{conj.} \ar_-{F_{\ast} + \text{Incl.}}[r] & 
\Rg(H\mZ,D_kB)/\mbox{conj.} } \end{equation}
is commutative. By Lemmas \ref{FGL fibre} and \ref{ring fibre},
the horizontal maps identify the respective homotopy fibres of the
truncation maps over the basepoints $F$ and $F_*$.  
Since formal group laws account for all components of the target space,
we have thus reduced the inductive step of the proof
of Theorem \ref{k-bud comparison}, and hence of the main theorem, to showing 

\begin{theorem} \label{thm-kappa S^k}
The map
\[  \kappa_{B\tensor {\mathcal S}^k} \,\co\, \widetilde{\Z}(B\tensor {\mathcal S}^k) \
\to \ \der(H\mZ,B\tensor {\mathcal S}^k)/\mbox{\em conj.} \]
is a weak equivalence for all commutative rings $B$ and all $k\geq 1$.
\end{theorem}

The remaining sections are spent  verifying that the map 
$\kappa_{B\tensor {\mathcal S}^k}$ is indeed a weak equivalence. 
If $k$ is not a prime power, then source and target
of the map $ \kappa_{B\tensor {\mathcal S}^k}$ are weakly contractible,
compare Remark \ref{constant unless prime power}. 
However there is no need to make this case distinction until
the very end (see Step 1 in the proof of Theorem \ref{deduction theorem}).

\section{Symmetric 2--cocycles and derivations} \label{2-cocycles}

In this section we provide some general constructions which will
be needed in the sequel. 
We consider functors $G$ from finitely generated free abelian groups 
to all abelian groups which are {\em reduced} in the sense that $G(0)\iso 0$.
For such functors we discuss symmetric 2--cocycles and show how these
lead to Gamma-ring maps into the split extension 
$H\mZ\times G^!$~\eqref{split extension}.
We also recall the Dold--Puppe stabilization $G_{st}$ 
of the functor $G$~\eqref{DP and Q}.
The work of this section is summarized in a certain map (\ref{lambda})
\[  \lambda_G \,\co\,\Z(G) \  \to
\map_{\mathcal{GR}}(H\mZ\rscript{c},H\mZ\times G_{st}^!)_{hG_{st}(\mZ)} \]
with source the groupoid $\Z(G)$ of symmetric 2--cocycles of $G$.
In the case where $G$ is the symmetric power functor $B\tensor {\mathcal S}^k$, 
we identify the map $\lambda_{B\tensor {\mathcal S}^k}$ 
(up to weak equivalence) with the map $\kappa_{B\tensor {\mathcal S}^k}$ 
of Theorem \ref{thm-kappa S^k} in Section \ref{simplification}.

\begin{numbered paragraph} \label{def-cocycles}
{\bf Cocycles of a functor}\qua
We let  $\F$ denote the abelian category of reduced functors 
from finitely generated free abelian groups to all abelian groups.
Suppose  $G\in \F$ is a such functor.
Then a {\em symmetric 2--cocycle} with values in $G$ 
is an element $c\in G(\mZ\oplus\mZ)$ which is 
\renewcommand{\labelenumi}{\rm(\alph{enumi})}
\begin{enumerate}
\item fixed by the involution of $G(\mZ\oplus\mZ)$
induced by the interchange of summands and
\item in the kernel of the map
\[ G\left(\begin{array}{ccc} 1 & 0 & 0 \\ 0 & 1 & 0 \end{array}\right) 
- G\left(\begin{array}{ccc} 1 & 0 & 0 \\ 0 & 1 & 1 \end{array}\right)
+ G\left(\begin{array}{ccc} 1 & 1 & 0 \\ 0 & 0 & 1 \end{array}\right)
- G\left(\begin{array}{ccc} 0 & 1 & 0 \\ 0 &0 & 1  \end{array}\right) \ : \]
\[ G(\mZ\oplus \mZ) \ \to \  G(\mZ\oplus \mZ\oplus \mZ) \ . \] 
\end{enumerate}
We denote by Z$^2_s(G)$ the group of symmetric 2--cocycles in $G$. 

A group homomorphism
\[ \theta \,\co\,G(\mZ) \ \to \  \mbox{Z}^2_s(G) \subset 
G(\mZ\oplus\mZ) \]
is defined by 
\[ \theta \ = \ G(1,0)- G(1,1) + G(0,1)  \ . \]
The cocycle $\theta(a)$ associated to $a\in G(\mZ)$ is sometimes referred to as
the {\em principal} cocycle of $a$.
We denote by  $\Z(G)$ the translation category of the action of $G(\mZ)$ 
on the set $\mbox{Z}^2_s(G)$ of symmetric 2--cocycles given by 
$(a,c)\longmapsto \theta(a)+c$.
More precisely,  $\Z(G)$ is the groupoid whose objects are
the symmetric 2--cocycles Z$^2_s(G)$ and where
the set of morphisms from a cocycle $c$ to a cocycle $c'$ consists of
those elements $a\in G(\mZ)$ such that $c'=\theta(a)+c$.
Composition in $\Z(G)$ is given by addition in the group $G(\mZ)$.
\end{numbered paragraph}

\begin{example}
If $G=B\tensor {\mathcal S}^k$ is the symmetric power functor, then
$G(\mZ\oplus \mZ)$ is the group of homogenous polynomials of degree $k$
in two variables over the ring $B$. 
For $c\in (B\tensor {\mathcal S}^k)(\mZ\oplus \mZ)$, the cocycle condition~(a)
translates into $c(x,y)=c(y,x)$ and condition~(b) translates into the equation
$$ c(x,y) + c(x+y,z) \ = \ c(x,y+z) + c(y,z) \ . $$
Hence the symmetric 2--cocycles with values in $B\tensor {\mathcal S}^k$ 
coincide with the homogenous 2--cocycles of degree $k$ as defined 
in \ref{B cocycles}. 
Moreover, the group $(B\tensor {\mathcal S}^k)(\mZ)$ is isomorphic
to the additive group of $B$ and the map 
$$ \theta\,\co\,B\iso (B\tensor {\mathcal S}^k)(\mZ) 
\ \to \ \mbox{Z}^2_s(B\tensor {\mathcal S}^k)$$
sends  $b\in B$ to the cocycle $b\cdot\left[ x^k + y^k - (x+y)^k \right]$.
Hence the cocycle groupoid for the functor 
$B\tensor {\mathcal S}^k$ as defined in \ref{def-cocycles} 
coincides with $\Z(B\tensor {\mathcal S}^k)$ as defined in \ref{B cocycles}.
\end{example}

\begin{numbered paragraph} \label{universal cocycle}
{\bf The universal 2--cocycle}\qua
We let $I\in\F$ denote the inclusion functor and $P\in \F$ the functor
which takes an abelian group $A$ to the reduced free abelian group
generated by the underlying pointed set of $A$.
The Yoneda isomorphism $\Hom_{\F}(P,G)\iso G(\mZ)$ shows that
$P$ is a projective object of $\F$.
Evaluation gives a natural epimorphism 
\[  P(A) \ = \ \widetilde{\mZ}[A] \ \to \ A = I(A) \ , \]
ie, an epimorphism $\epsilon\co P\to I$ in the category $\F$. 
We let $J$ denote the kernel; the relevance for us is that $J$ represents the
symmetric 2--cocycle functor.
The element 
\begin{equation} \label{eqn-universal cocycle}  
c_u \ =\ [1,0]-[1,1]+[0,1] \ \in \ \widetilde\mZ[\mZ\oplus\mZ] 
= P(\mZ\oplus\mZ) 
\end{equation}
is in the kernel of the evaluation map, so it is an element 
of $J(\mZ\oplus\mZ)$.
As an element of $P(\mZ\oplus\mZ)$, $c_u$ is the principal cocycle 
associated to $[1]\in P(\mZ)$. Hence $c_u$ is a symmetric 2--cocycle 
of $J$ (but it is not principal for $J$ 
because the element $[1]\in P(\mZ)$ does not belong to $J(\mZ)$). 
\end{numbered paragraph}

\begin{lemma} \label{lem-universal cocycle}
The  symmetric 2--cocycle  $c_u$ with values in the
functor $J$ is universal in the sense that the map 
\[ \Hom_{\F}(J,G) \ \to \ \mbox{\em Z}^2_s(G) \ ,
\quad f \longmapsto f(c_u) \]
is an isomorphism for all functors $G\in\F$.
\end{lemma}
\begin{proof} Both  $\Hom_{\F}(J,G)$ and $\mbox{Z}^2_s(G)$ are
additive and left exact in the functor $G$. Hence it suffices to
check the claim for a set of injective cogenerators of the category $\F$.
If $A$ is a finitely generated free abelian group we define a functor
$\I_A$ by the formula
\[ \I_A(M) \ = \ \map_*(\Hom(M,A),\mQ/\mZ) \ . \] 
Here $\Hom(M,A)$ denotes the set of group homomorphisms from $M$ to $A$
and `$\map_*$' refers to the group of set-theoretic maps from 
$\Hom(M,A)$ to $\mQ/\mZ$ preserving 0,
with group structure by pointwise addition.
The Yoneda isomorphism
\[ \Hom_{\F}(G,\I_A) \ \iso \ \Hom(G(A),\mQ/\mZ) \] 
implies that $\I_A$ is injective and that the functors $\I_A$ form a collection
of injective cogenerators as $A$ varies.

It remains to verify that for all finitely generated free abelian groups $A$
evaluation at the cocycle $c_u\in\mbox{Z}^2_s(J)$ is an isomorphism
from $\Hom_{\F}(J,\I_A)$ to $\mbox{Z}^2_s(\I_A)$.
We claim that the group $\mbox{Z}^2_s(\I_A)$ can be identified
with the quotient of the group $\map_*(A,\mQ/\mZ)$ of pointed set maps from 
$A$ to $\mQ/\mZ$ by the subgroup $\Hom(A,\mQ/\mZ)$ of additive maps.
We use the natural basis of $\mZ\oplus\mZ$ to identify 
$\I_A(\mZ\oplus\mZ)=\map_*(\Hom(\mZ\oplus\mZ,A),\mQ/\mZ)$ with
the group of set-theoretic maps $g\co A\oplus A\to\mQ/\mZ$ 
satisfying $g(0,0)=0$.
Under this identification the cocycle conditions for elements of
$\I_A(\mZ\oplus\mZ)$ translate into the conditions
\[ g(x,y) \ = \ g(y,x) \mbox{\quad and \quad} 
g(x,y)  \ + \ g(x+y,z)  \ = \ g(x,y+z) \ + \ g(y,z)  \]
for all $x,y,z\in A$ on the function $g\in\map_*(A\oplus A,\mQ/\mZ)$.
In other words: $g$ is a factor set of an abelian group extension of $A$
by $\mQ/\mZ$.
Since $A$ is free abelian, every such extension splits, so $g$ is principal, 
ie, there exists a function $h\co A\to\mQ/\mZ$ satisfying
$g(x,y) \ = \ h(x) - h(x+y) + h(y)$ and $h(0)=0$. 
Moreover, $h$ is uniquely determined by this up to an additive function. 

All this means that the bottom row in the commutative diagram of abelian groups
\[\xymatrix{ 0 \ar[r] &  \Hom_{\F}(I,\I_A) \ar_{\iso}[d] \ar[r] & 
\Hom_{\F}(P,\I_A) \ar_{\iso}[d] \ar[r] &  
\Hom_{\F}(J,\I_A) \ar^{f\mapsto f(c_u)}[d] \ar[r] & 0 \\
0 \ar[r] & \Hom(A,\mQ/\mZ) \ar[r] & \map_*(A,\mQ/\mZ) \ar_-{\vartheta}[r] &
\mbox{Z}^2_s(\I_A) \ar[r] & 0 }\]
is exact where $\vartheta$ is defined by $\vartheta(h)(x,y)=h(x)-h(x+y)+h(y)$. 
The upper row is exact since $\I_A$ is an injective object and
$J$ is the kernel of the epimorphism $\epsilon\co P\to I$.
The left and middle vertical maps are special cases of the Yoneda
isomorphism $\Hom_{\F}(G,\I_A)\iso \Hom(G(A),\mQ/\mZ)$, so the right
vertical map is also an isomorphism and $c_u\in\mbox{Z}^2_s(J)$ 
is indeed a universal symmetric 2--cocycle.
\end{proof}

\begin{numbered paragraph} \label{con-derivation}
{\bf The universal derivation}\qua 
We saw in Lemma \ref{lem-universal cocycle} 
that the functor $J=\mbox{kernel}(\epsilon\co P\to I)$ 
supports a universal symmetric 2--cocycle. 
Now we construct a universal derivation,
ie, a certain homomorphism of Gamma-rings
\[ H\mZ\ \to \  H\mZ\times J^!  \] 
into the split extension (\ref{split extension}) of $H\mZ$ by $J^!$.
The first component of this map is the identity map of $H\mZ$.
To describe the second component we consider the map of
$\Gamma$--spaces
\[ \eta\circ 1 \ - \ 1\circ \eta \,\co\,H\mZ \ \to \ H\mZ\circ H\mZ \] 
where $\eta\co \mS\to H\mZ$ is the unit map of $H\mZ$ given by inclusion of
generators into the free abelian groups, and where we use the
identifications $\mS\circ H\mZ\iso H\mZ\iso H\mZ\circ \mS$.
We note that $H\mZ\circ H\mZ=P^!$ and that the composite of 
$\eta\circ 1 - 1\circ \eta$ with the evaluation map 
$\epsilon^!\co H\mZ\circ H\mZ=P^!\to I^!=H\mZ$ is trivial.
So we can define $d_u\co H\mZ\to J^!$ as the unique morphism of $\Gamma$--spaces
whose composite with the inclusion $J^!\subseteq P^!$ is
the difference $\eta\circ 1 - 1\circ \eta$.
We refer to $d_u$ as the {\em universal derivation}.
Now we claim that the map
\[ (1,d_u) \,\co\, H\mZ\ \to \  H\mZ\times J^! \] 
is a homomorphism of Gamma-rings. The only thing to verify is the 
multiplicativity, and it suffices to do this after composition with
the injective  Gamma-ring map
$H\mZ\times J^!\to H\mZ\times P^!$ induced by the inclusion $J\to P$.
By the definition of the product of 
$H\mZ\times P^!=H\mZ\times (H\mZ\circ H\mZ)$ in \ref{split extension}
this boils down to verifying the commutativity of the  diagram
\[\xymatrix@C=15mm{ H\mZ\circ H\mZ \ar^{(1\circ(\eta\circ 1 - 1\circ\eta)\ ,\
(\eta\circ 1 - 1\circ \eta)\circ 1)}[rr] \ar_{\mu}[d] & &
H\mZ^{\circ 3} \times H\mZ^{\circ 3} \ar^{(\mu\circ 1)\times(1\circ\mu)}[d] \\
H\mZ  \ar_-{\qquad \eta\circ 1 - 1\circ \eta\quad }[dr] & & 
(H\mZ\circ H\mZ)\times (H\mZ\circ H\mZ) \ar^{+}[dl] \\
& H\mZ\circ H\mZ }\]
where $\mu\co H\mZ\circ H\mZ\to H\mZ$ is the ``multiplication'' induced
by the evaluation map $P\to I$. 
The commutativity of the diagram in turn follows from the identities
\[ (\mu\circ 1)(1\circ \eta\circ 1) \ = \  1_{H\mZ\circ H\mZ}
\mbox{\qquad and \qquad} 
(\mu\circ 1)(1\circ 1\circ \eta) \ =  \ (1\circ\eta)\mu \]
(here juxtaposition means composition of $\Gamma$--space maps)
and their analogues for $\mu\circ 1$ replaced by $1\circ \mu$.
\end{numbered paragraph}

Hence for a functor $G\in\F$ we can define a map
\begin{equation} \label{prelim lambda} \mbox{Z}^2_s(G) \ \iso \ \Hom_{\F}(J,G) 
\ \to \ \map_{\mathcal{GR}}(H\mZ,H\mZ\times G^!) 
\end{equation}
by sending a morphism $f\co J\to G$ to the composite Gamma-ring map
\[ H\mZ \ \varr{1.5cm}{d_u} \  H\mZ\times J^! 
\ \varr{1.5cm}{1 \times f^!} \  H\mZ\times G^!\ . \] 
If $c\in \mbox{Z}^2_s(G)$ is a cocycle represented by $f_c\co J\to G$,
then we refer to the above Gamma-ring map as the 
{\em derivation} associated to the 2--cocycle $c$.

The group $G(\mZ)$ maps to the underlying monoid 
$(H\mZ \times G^!)(1^+)\iso \mZ\times G(\mZ)$
of the split extension via $a\longmapsto (1,a)$, 
hence it acts on the Gamma-ring $H\mZ \times G^!$ by conjugation
(\ref{con-action of group}).
The group $G(\mZ)$ also acts on the symmetric 2--cocycles via
the homomorphism  $\theta\co  G(\mZ) \to \mbox{Z}^2_s(G)$ (\ref{def-cocycles}).

\begin{lemma} \label{G(Z) equivariance} The map
\[ \mbox{\em Z}^2_s(G) \ \iso \ \Hom_{\F}(J,G) \ \to 
\ \map_{\mathcal{GR}}(H\mZ,H\mZ\times G^!) \] 
which sends a morphism  $f\,\co J\to G$ to the Gamma-ring map
$(1,f^!\circ d_u)\,\co$\break $H\mZ\to H\mZ\times G^!$
is equivariant for the action of $G(\mZ)$.
\end{lemma}
\begin{proof} Let $a\in G(\mZ)$ be an element and $f\co J\to G$ 
a morphism of functors. Under the Yoneda isomorphism 
$G(\mZ)\iso\Hom_{\F}(P,G)$ and the isomorphism  
$\mbox{Z}^2_s(G)\iso \Hom_{\F}(J,G)$ of Lemma \ref{lem-universal cocycle},
the map  $\theta\co  G(\mZ) \to \mbox{Z}^2_s(G)$ corresponds to the
map $\Hom_{\F}(P,G)\to \Hom_{\F}(J,G)$ induced by the inclusion 
$\iota\co J \to P$.

So we have to verify the equality
\[ (1,d_{f+a\circ \iota}) \ = \ (1,a) \cdot (1,d_f) \cdot (1,a)^{-1} \] 
as maps $H\mZ\to H\mZ\times G^!$. Only the second component matters.
If we substitute definitions \ref{con-action of group} of the
conjugation action and \ref{split extension} of the multiplication
in $H\mZ\times G^!$ we see that the second component of right hand side 
is a sum
\[ \mbox{proj}_{G^!} \circ \left[ (1,a) \cdot (1,d_f) \cdot (1,a)^{-1} 
\right] \ = \ r_{\mbox{-}a} \ + \ d_f \ + \ l_a  \ . \]
Here $l_a\co H\mZ\to G^!$ is left multiplication with $a\in G^!(1^+)\iso G(\mZ)$,
ie, the composition
\[ H\mZ \iso \mS\sm H\mZ \ \varr{1cm}{a\sm 1} \  G^!\sm H\mZ
\varr{1.5cm}{\text{action}} \ G^! \ , \]
and $r_{\mbox{-}a}$ is right multiplication with $-a\in G^!(1^+)$.

Since $d_{f+a\circ \iota}=d_f +d_{a\circ \iota}$ it suffices to show that 
\[ d_{a\circ \iota} \ = \ r_{\mbox{-}a} \ + \ l_a  \] 
as maps  $H\mZ\to G^!$. By naturality it is enough to check the universal
example, ie, to take $G=P$ and $a=[1]\in\widetilde{\mZ}[\mZ]=P(\mZ)$,
which corresponds to the identity of $P$ under the Yoneda isomorphism.
By definition 
\[ d_\iota \ = \ \eta\circ 1- 1\circ \eta \,\co\,H\mZ\ \to \  
P^! = H\mZ\circ H\mZ \ . \]
So the claim follows since  $l_{[1]} = \eta\circ 1$ and 
$r_{\mbox{-}[1]} = - 1\circ \eta$.
\end{proof}

\begin{numbered paragraph} \label{DP and Q}
{\bf Dold--Puppe stabilization and MacLane's $Q$--construction}\qua
We recall the Dold--Puppe stabilization 
of a functor $G\in\F$ (compare \cite[8.3]{Dold-Puppe}).
We work with a specific model
for the  Dold--Puppe stabilization which uses MacLane's
{\em cubical construction} \cite[Sec.\ 12]{EM:cubical}.
In the original paper of Eilenberg and MacLane the cubical construction
was defined for the free functor $P=\widetilde{\mZ}$, but
the definition extends to arbitrary reduced functors in $\F$,
see \cite[Sec. 4]{Pira:approximation} or \cite[6.2]{Johnson-McCarthy}
for the definition.
A convenient reference for the relationship between Dold--Puppe 
stabilization and  MacLane's cubical construction 
is \cite{Johnson-McCarthy}. 
The cubical construction $QG$ of a functor
$G\in\F$ is a chain complex of functors, concentrated 
in non-negative dimensions, with the following properties:
\renewcommand{\labelenumi}{\rm(\alph{enumi})}
\begin{enumerate}
\item $QG$ is {\em homotopy-additive} 
in the sense that for every pair of finitely generated free abelian groups
$A$ and $A'$ the canonical map 
$$ QG(A)\oplus QG(A')\ \to \ QG(A\oplus A') $$ 
is a quasi-isomorphism.
\item there is a natural isomorphism
\[ H_* \, QG(\mZ)  \ \iso \  L\rscript{st}_* \, G \ \iso \  \pi_* \, G^! \ , \]
ie, the homology groups of the complex $QG(\mZ)$
are isomorphic to the Dold--Puppe stable derived functors of $G$ and
to the stable homotopy groups of the $\Gamma$--space $G^!$. 
\item In dimension zero, $(QG)_0=G$ and in positive dimensions $QG$ is a finite
sum of higher order cross-effects (see \cite[Sec.9]{EM:K(pi)-II} or
\cite[Sec. 7]{Johnson-McCarthy}) of~$G$.
\item As a functor of $G$, the assignment $G\longmapsto QG$ is 
additive, exact, and commutes with limits and colimits.
\item Suppose the functor $G$ is {\em diagonalizable}, ie, there exists
a functor 
\[ \bar G \,\co\,(\mbox{f.g. free ab. groups})\times
(\mbox{f.g. free ab. groups}) \ \to \ \Ab \]
of two variables satisfying $\bar G(A,0)\iso 0\iso \bar G(0,A)$
and a natural isomorphism $G(A)\iso \bar G(A,A)$.
Then the complex $QG$ is acyclic.
\end{enumerate}

Property (a) is proved in \cite[4.2]{Pira:approximation} and 
\cite[6.3]{Johnson-McCarthy}.
The first isomorphism of part (b) follows from \cite[4.1]{Pira:approximation}
or  \cite[7.5]{Johnson-McCarthy};  essentially by definition
the stable homotopy groups of the $\Gamma$--space $G^!$ are the
Dold--Puppe stable derived functors of $G$.
Part (c) is  proved in \cite[6.3]{Johnson-McCarthy}. 
Property (d) follows from (c) since taking cross-effects 
is exact and commutes with limits and colimits.
If $G$ is diagonalizable, then \cite[5.20]{Dold-Puppe} shows
that the stable derived functors of $G$ are trivial.
So by part (b) the complex $QG(\mZ)$ is acyclic and by part (a)
$QG$ is acyclic as a complex of functors.

Properties (a) and (b) already characterize  $QG$ up to a chain of natural
quasi-isomorphisms; this is because on the level of homotopy
categories, $G\longmapsto QG$ is left adjoint to the inclusion of
the subcategory of homotopy additive complexes of functors.

We denote by $G_{st}$ the simplicial functor which corresponds to 
$QG$ under the Dold--Kan equivalence between simplicial objects
and non-negative chain complexes in the abelian category $\F$.
So $G_{st}$ is defined by the property that its normalized chain complex
is isomorphic to $QG$. By property (c), the functor of zero-simplices
of $G_{st}$ is $G$, which gives a map $G\to G_{st}$ which induces
isomorphism of stable homotopy groups upon passage to associated
$\Gamma$--spaces $G^!\to G_{st}^!$.
\end{numbered paragraph}

Construction \ref{split extension} which associates the Gamma-ring
$H\mZ\times G^!$ to a functor $G\in\F$ makes perfect 
sense for {\em simplicial functors}, ie, simplicial objects in the
abelian category $\F$. Moreover, the stabilization map $G\to G_{st}$
induces a stable equivalence of Gamma-rings 
$H\mZ\times G^!\to H\mZ\times G_{st}^!$.
Combining the map (\ref{prelim lambda}) with this Gamma-ring homomorphism
and the approximation map $H\mZ\rscript{c}\to H\mZ$
gives a map
\begin{eqnarray*}
 \mbox{Z}^2_s(G) \ \iso \ \Hom_{\F}(J,G) & 
\varr{2.5cm}{f\ \mapsto \ (1,\ f^!\circ d_u)} & 
\map_{\mathcal{GR}}(H\mZ,H\mZ\times G^!) \\
& \varr{2.5cm}{} &
\map_{\mathcal{GR}}(H\mZ\rscript{c},H\mZ\times G_{st}^!) \ .
\end{eqnarray*}
By Lemma \ref{G(Z) equivariance} this map is equivariant for the action
of $G(\mZ)$ by translation and conjugation respectively.
Hence passing to homotopy orbits yields a map
\begin{equation} \label{lambda}  \lambda_G \,\co\,\Z(G) \ \iso \
\mbox{Z}^2_s(G)_{hG(\mZ)} \  \to \
\map_{\mathcal{GR}}(H\mZ\rscript{c},H\mZ\times G_{st}^!)_{hG_{st}(\mZ)} \ . 
\end{equation}
Here we identified the classifying space of the groupoid $\Z(G)$ 
of symmetric 2--cocycles (\ref{def-cocycles}) 
with the homotopy orbit space of the action of $G(\mZ)$ 
on the set $\mbox{Z}^2_s(G)$.

\section{A change of models} \label{simplification}

So far we have reduced the proof of our main theorem, Theorem \ref{main},
to the verification that the map 
\[  \kappa_{B\tensor {\mathcal S}^k} \,\co\, \widetilde{\Z}(B\tensor {\mathcal S}^k) \
\to \ \der(H\mZ,B\tensor {\mathcal S}^k)/\mbox{conj.} \]
constructed in Section \ref{singular extensions} 
is a weak equivalence for all $k\geq 1$.
In this section we bring the map $\kappa_{B\tensor {\mathcal S}^k}$
into a more manageable form by constructing a commutative square
\begin{equation} \label{comparison square} 
\xymatrix@C=20mm{ \widetilde{\Z}(B\tensor {\mathcal S}^k) \ar_{\sim}[d] 
\ar^-{\kappa_{B\tensor {\mathcal S}^k}}[r] & 
\der(H\mZ,B\tensor {\mathcal S}^k)/\mbox{conj.} \ar^{\sim}[d] \\
\Z(B\tensor {\mathcal S}^k) \ar_-{\lambda_{B\tensor {\mathcal S}^k}}[r] & 
\map_{\mathcal{GR}}(H\mZ\rscript{c},H\mZ\times (B\tensor {\mathcal S}^k)_{st}^!)_{h(B\tensor {\mathcal S}^k)_{st}(\mZ)} }\end{equation}
in which the vertical maps are weak equivalences.
The lower horizontal map $\lambda_{B\tensor {\mathcal S}^k}$ 
is an instance of (\ref{lambda}).
This then leaves us with the task to verify that
$\lambda_{B\tensor {\mathcal S}^k}$ is a weak equivalence for all $k\geq 1$.

The construction of the square (\ref{comparison square}) is technical.
The idea is that by properties \ref{DP and Q} (a) and (b) 
of the cubical construction,
$(B\tensor {\mathcal S}^k)^!_{st}$ is a stably fibrant model 
of the $\Gamma$--space $(B\tensor {\mathcal S}^k)^!$.
Hence the map 
$H\mZ\times (B\tensor {\mathcal S}^k)^!\to H\mZ\times (B\tensor {\mathcal S}^k)_{st}^!$ 
is a stable equivalence of Gamma-rings with fibrant target.
In particular, the space of Gamma-ring maps 
$\map_{\GR}(H\mZ\rscript{c}, H\mZ\times (B\tensor {\mathcal S}^k)_{st}^!)$ is a model
for the homotopy invariant morphism space, ie, it is weakly equivalent
to the derivation space $\der(H\mZ,B\tensor {\mathcal S}^k)$.
In order to identify these two derivation spaces we have to take a little
more care because we need to work relative to the classifying 
space of symmetric 2--cocycles. 

We reexamine  Construction \ref{conjugation construction}, applied to
the Gamma-ring $H\mZ \times (B\tensor {\mathcal S}^k)^!$.
This construction yields a commutative diagram of Gamma-rings
\[\xymatrix@C=11mm{ \mS[B] \ar[d] & {\mathbb S}[B] \ar[d] \ar@{=}[l] \ar[r] &
\mS[U(H\mZ\times (B\tensor {\mathcal S}^k)^!)^{\times}] \ar[d] \\
(H\mZ \times (B\tensor {\mathcal S}^k)^!)^{\f} & 
(H\mZ \times (B\tensor {\mathcal S}^k)^!)_1  \ar_{\sim}[r] 
\ar@{->>}^{\sim}[l] & (H\mZ \times (B\tensor {\mathcal S}^k)^!)_3 }\] 
in which the lower horizontal maps are stable equivalences 
between stably fibrant Gamma-rings. 
Here we work relative to the homomorphism 
$$ B\ \to \ 
H\mZ \times (B\tensor {\mathcal S}^k)^!(1^+)= \mZ[x]\times B[x^k] $$ 
which sends $b\in B$ to $x+b\cdot x^k$.
Furthermore, the induced map from the simplicial group 
$U(H\mZ\times (B\tensor {\mathcal S}^k)^!)^{\times}$ 
to the invertible components 
of the underlying monoid of $(H\mZ \times (B\tensor {\mathcal S}^k)^!)_3$ 
is a weak equivalence. 

Since the Gamma-ring  $H\mZ\times (B\tensor {\mathcal S}^k)_{st}^!$ is stably fibrant 
and the map 
\[ H\mZ\times (B\tensor {\mathcal S}^k)^!\ \to \ 
(H\mZ\times (B\tensor {\mathcal S}^k)^!)^{\f} \] 
in the functorial fibrant replacement is an acyclic cofibration, 
we can choose a Gamma-ring map from $(H\mZ\times (B\tensor {\mathcal S}^k)^!)^{\f}$ to 
$H\mZ\times (B\tensor {\mathcal S}^k)_{st}^!$ under $H\mZ\times (B\tensor {\mathcal S}^k)^!$. 
This map will automatically be a stable equivalence. 
We can perform constructions \ref{conjugation construction} and 
\ref{construction weak map} starting from either of these two fibrant 
replacements. The stable equivalence between them induces a weak
equivalence between the two homotopy orbit spaces of the conjugation action.
In other words, we can assume that the stably fibrant replacement of
the split extension $H\mZ\times (B\tensor {\mathcal S}^k)^!$ which was chosen
in the beginning is {\em equal} to 
$H\mZ\times (B\tensor {\mathcal S}^k)_{st}^!$.

The simplicial group $U(H\mZ\times (B\tensor {\mathcal S}^k)^!)^{\times}$ 
is defined according to the recipe (\ref{conjugation construction})
by a factorization 
in the model category of simplicial monoids
\[ B \ \to \ c(H\mZ\times (B\tensor {\mathcal S}^k)^!)^{\times}
\varr{1cm}{\sim} \ (H\mZ\times (B\tensor {\mathcal S}^k)^!)^{\times} 
\]
and then forming the algebraic group completion of 
$c(H\mZ\times (B\tensor {\mathcal S}^k)^!)^{\times}$.
The simplicial monoid $(H\mZ\times (B\tensor {\mathcal S}^k)^!)^{\times}$ is,
by definition, the union of the invertible components in
$(H\mZ\times (B\tensor {\mathcal S}^k)_{st}^!)(1^+) \iso 
\mZ\times (B\tensor {\mathcal S}^k)_{st}(\mZ)$. So the simplicial monoid
\[ (H\mZ\times (B\tensor {\mathcal S}^k)_{st}^!)^{\times} \ \iso \ 
\{\pm 1\}\times (B\tensor {\mathcal S}^k)_{st}(\mZ) \]
is already a simplicial group, hence 
there exists a unique extension to a homomorphism of simplicial groups 
$U(H\mZ\times (B\tensor {\mathcal S}^k)^!)^{\times} \to 
\{\pm 1\}\times  (B\tensor {\mathcal S}^k)_{st}(\mZ)$
which is necessarily a weak equivalence by Lemma \ref{good monoid proof}. 

Since $k\geq 2$, the stable derived functors of $B\tensor {\mathcal S}^k$ 
vanish in dimension 0 and 1 \cite[12.3]{Dold-Puppe}, so the simplicial 
abelian group $(B\tensor {\mathcal S}^k)_{st}(\mZ)$ is simply connected, and
the unit component of the monoid 
$(H\mZ\times (B\tensor {\mathcal S}^k)_{st}^!)(1^+)$ is equal to  
$(B\tensor {\mathcal S}^k)_{st}(\mZ)$.
Restriction to the unit components thus
gives a weak equivalence of connected simplicial groups
$U(H\mZ\times (B\tensor {\mathcal S}^k)^!)^{\times}_1 \to (B\tensor {\mathcal S}^k)_{st}(\mZ)$.

The next step in the construction of the conjugation action
(\ref{conjugation construction}) was the formation of the pushout 
$(H\mZ\times  (B\tensor {\mathcal S}^k)^!)_2$ of Gamma-rings:
\[\xymatrix{\mS[c(H\mZ\times (B\tensor {\mathcal S}^k)^!)^{\times}] \ar[r] \ar_{\sim}[d] &
(H\mZ\times  (B\tensor {\mathcal S}^k)^!)_1 \ar[d] \ar[r] & 
H\mZ\times  (B\tensor {\mathcal S}^k)_{st}^! \\
\mS[U(H\mZ\times (B\tensor {\mathcal S}^k)^!)^{\times}] \ar[r] & 
(H\mZ\times  (B\tensor {\mathcal S}^k)^!)_2 \ar@{-->}[ur] \ar@{>->}_{\sim}[d] \\
& (H\mZ\times (B\tensor {\mathcal S}^k)^!)_3 \ar@{-->}[uur] }\]
The simplicial group map  
$U(H\mZ\times (B\tensor {\mathcal S}^k)^!)^{\times} \to 
\{\pm 1\}\times (B\tensor {\mathcal S}^k)_{st}(\mZ)$
adjoins to a homomorphism of Gamma-rings from the spherical group ring 
$\mS[U(H\mZ\times (B\tensor {\mathcal S}^k)^!)^{\times}]$ to 
$H\mZ\times  (B\tensor {\mathcal S}^k)_{st}^!$. 
This in turn induces a map from the pushout 
$(H\mZ\times (B\tensor {\mathcal S}^k)^!)_2$ to  $H\mZ\times (B\tensor {\mathcal S}^k)_{st}^!$,
 represented as the upper dotted arrow in
the previous diagram. Since the approximation map  
$(H\mZ\times (B\tensor {\mathcal S}^k)^!)_2\to (H\mZ\times (B\tensor {\mathcal S}^k)^!)_3$ 
is an acyclic cofibration of Gamma-rings, we can finally choose 
an extension to a stable equivalence from  
$(H\mZ\times (B\tensor {\mathcal S}^k)^!)_3$ 
to  $H\mZ\times (B\tensor {\mathcal S}^k)_{st}^!$. 
Since this map was constructed relative to the group ring 
of the simplicial group  $U(H\mZ\times (B\tensor {\mathcal S}^k)^!)^{\times}$, 
it is equivariant 
with respect to the conjugation action of  
$U(H\mZ\times (B\tensor {\mathcal S}^k)^!)^{\times}_1$ 
on $(H\mZ\times (B\tensor {\mathcal S}^k)^!)_3$ and (through the map  
$U(H\mZ\times (B\tensor {\mathcal S}^k)^!)^{\times}_1 \to (B\tensor {\mathcal S}^k)_{st}(\mZ))$ 
on $H\mZ\times (B\tensor {\mathcal S}^k)_{st}^!$. 
By passage to homotopy orbits we get a weak equivalence
\[  \der(H\mZ,B\tensor {\mathcal S}^k)/\mbox{conj.} \ = \
\map_{\mathcal{GR}}(H\mZ\rscript{c},(H\mZ \times (B\tensor {\mathcal S}^k)^!)_3)_{hU(H\mZ\times (B\tensor {\mathcal S}^k)^!)^{\times}_1} \]
\[ \hspace*{4cm} \varr{1.5cm}{\sim} \ 
\map_{\mathcal{GR}}(H\mZ\rscript{c},H\mZ\times (B\tensor {\mathcal S}^k)_{st}^!)_{h(B\tensor {\mathcal S}^k)_{st}(\mZ)} \]
Moreover, the square (\ref{comparison square}) commutes.

\section{A useful adjunction} \label{universal equivalence}

By the results of the previous two sections, the proof of the main theorem
is reduced to an identification of the space of Gamma-ring maps
\[ \map_{\mathcal{GR}}(H\mZ\rscript{c},
H\mZ\times (B\tensor {\mathcal S}^k)_{st}^!) \]
(or more precisely a certain homotopy orbit space thereof)
with the classifying space of symmetric 2--cocycles.
In this section we use an adjunction to reinterpret the above mapping space 
in terms of the category $s\F$ of simplicial functors from the category of 
finitely generated free abelian groups to the category of abelian groups.

We note in Lemma \ref{quillen pair} that the construction 
(\ref{split extension}) of the split singular extension has a left adjoint
\[ \Lc \,\co\,\GR/H\mZ \ \to \  s\F \] 
from the category of  Gamma-rings over
$H\mZ$ to the category of simplicial functors. Moreover, the two functors
form a simplicial Quillen adjoint pair of model categories. 
So the mapping space we are interested in is isomorphic to the mapping space
\[ \map_{s\F}(\Lc(H\mZ\rscript{c}),(B\tensor {\mathcal S}^k)_{st}) \] 
in the category of simplicial functors.
To identify the mapping space in this adjoint form, 
we evaluate the left adjoint $\Lc$ on the cofibrant replacement
of the Eilenberg--MacLane Gamma-ring $H\mZ$. 

As in~\eqref{universal cocycle}, $J\in\F$ denotes the kernel 
of the evaluation map $\epsilon\co P\to I$.
By Lemma \ref{lem-universal cocycle} the functor $J$ supports 
the universal symmetric 2--cocycle.
The Gamma-ring map 
\[  H\mZ\rscript{c} \ \to \ H\mZ\times J^! \] 
which is the composite of  the approximation map 
$H\mZ\rscript{c}\to H\mZ$ and the universal derivation (\ref{con-derivation})
is adjoint to a map of simplicial functors
\begin{equation} \label{delta}  \delta \,\co\,\Lc(H\mZ\rscript{c}) 
\ \to \ J \ . \end{equation}
The main result of this section is

\begin{theorem}  \label{thm-adjoint of universal}
The map $\delta\co \Lc(H\mZ\rscript{c}) \to J$
which is adjoint to the universal derivation 
is a stable equivalence of simplicial functors.
\end{theorem}

\begin{remark}
For any functor $G\in\F$ the Dold--Puppe stabilization $G_{st}$ 
was defined so that its normalized chain
complex is the cubical construction $QG$ (\ref{DP and Q}).
So Theorem \ref{thm-adjoint of universal} 
and the Dold--Kan correspondence between simplicial
objects and chain complexes in the abelian category $\F$ imply
that the homotopy groups of the space
\[ \map_{\mathcal{GR}}(H\mZ\rscript{c},H\mZ\times G_{st}^!) 
\ \iso \  \map_{s\F}(\Lc(H\mZ\rscript{c}),G_{st}) \]
are isomorphic to the hyper-cohomology groups $\hExt_{\F}^*(J,QG)$,
for $*\leq 0$.
\end{remark}

The category $s\F$ of simplicial functors admits a stable model
structure, see \cite[6.13]{Sch:stable}. In this model structure,
a map $G\to \bar G$ is a weak equivalence or fibration if and only if
the associated map of $\Gamma$--spaces $G^!\to \bar G^!$ is a stable
equivalence or stable fibration respectively.
The stably fibrant objects are precisely the homotopy additive 
simplicial functors.

The split extension functor (\ref{split extension}) which sends $G\in s\F$ 
to $H\mZ\times G^!$, considered as a Gamma-ring over $H\mZ$, 
commutes with limits. 
Moreover, the category $s\F$ of simplicial functors is complete,
well-powered and has a set of cogenerators. 
So by Freyd's Adjoint Functor Theorem \cite[V.8 Cor.]{MacL:categories} 
there is a left adjoint $\Lc\co \GR/H\mZ\to s\F$.
The right adjoint $H\mZ\times (-)^!$ preserves stable equivalences 
and stable fibrations since in both categories these are defined on 
``underlying'' $\Gamma$--spaces. Hence we obtain

\begin{lemma} \label{quillen pair}
The functor which sends a simplicial functor $G\in s\F$ to 
$H\mZ\times G^!$, 
viewed as a Gamma-ring over $H\mZ$, is the right adjoint of
a Quillen adjoint pair between the category $s\F$ of simplicial functors,
endowed with the stable model structure, and the 
stable model category of Gamma-rings over $H\mZ$.
\end{lemma}

We prove Theorem \ref{thm-adjoint of universal} by constructing 
a commutative square of $H\mZ\rscript{c}$--bi\-modules
\begin{equation}\label{bar square}
\xymatrix@C=15mm{ \widetilde\B(H\mZ\rscript{c}) \ar^-{\eta_{H\mZ^c}}[r] 
\ar_{\text{ass.}}[d] & 
(B\Lc(H\mZ\rscript{c}))^! \ar^{(B\delta)^!}[d] \\
\widetilde\B(P)^! \ar_-{\widetilde u^!}[r] & (BJ)^! }\end{equation}
For a simplicial functor $G\in s\F$ we denote by 
$BG=\widetilde\mZ[S^1]\tensor G$
the (additive) bar construction, another simplicial functor. 
The bimodules $\widetilde\B(H\mZ\rscript{c})$ and $\widetilde\B(P)^!$ 
are ``multiplicative'' bar constructions defined below.
Three of the four objects are actually $H\mZ$--bimodules, which we
view as $H\mZ\rscript{c}$--bimodules via restriction of scalars.
In the square the left vertical and the two horizontal maps are 
stable equivalences
by Theorems \ref{Thm A}, \ref{Thm B} and \ref{Thm C} respectively.
Hence the map $(B\delta)^!\co (B\Lc(H\mZ\rscript{c}))^!\to (BJ)^!$
is a stable equivalence.
Since the assignment $G\longmapsto (BG)^!$ detects stable equivalences,
the map $\delta\co \Lc(H\mZ\rscript{c})\to J$ is indeed a stable equivalence
of simplicial functors.

\begin{numbered paragraph} \label{bar construction} 
{\bf A bar construction}\qua
The {\em reduced bar construction} is a functor
\[ \widetilde\B \,\co\,
\mathcal{GR}/H\mZ\rscript{c} \ \to \ 
H\mZ\rscript{c}\mbox{-mod-}H\mZ\rscript{c}  \]
from the category of Gamma-rings over $H\mZ\rscript{c}$ 
to the category of $H\mZ\rscript{c}$--bimodules. 
In this construction it is important that we start with Gamma-rings over
the cofibrant approximation $H\mZ\rscript{c}$, not just over $H\mZ$.
If we worked over $H\mZ$, the bar construction would have the
wrong homotopy type since the point set level smash product of $H\mZ$
with itself is not stably equivalent to $H\mZ\rscript{c}\sm H\mZ\rscript{c}$

If $Q$ is a Gamma-ring over $H\mZ\rscript{c}$, 
then the (unreduced) bar construction $\B(Q)$ 
is defined as the realization of a simplicial   
$H\mZ\rscript{c}$--bimodule which in simplicial dimension $n$ has the form
\[  H\mZ\rscript{c} \, \sm \, \underbrace{Q \, \sm \dots \, \sm \, Q}_{n}  \, 
\sm \, H\mZ\rscript{c} \ . \]
The simplicial structure maps are induced by the multiplication 
and unit map of $Q$ and the structure map $Q\to H\mZ\rscript{c}$. 
The inclusion of the 0-simplices induces a map 
\[ H\mZ\rscript{c}\,\sm\, H\mZ\rscript{c}\ \to \  \B(Q) \]
of $H\mZ\rscript{c}$--bimodules, and the reduced bar construction  
$\widetilde\B(Q)$ is the quotient of $\B(Q)$ by  
$H\mZ\rscript{c}\sm H\mZ\rscript{c}$.

For every simplicial functor $G$ there is a map 
\[ \tau_G \,\co\,\widetilde\B(H\mZ\rscript{c}\times G^!) 
\ \to \ (BG)^! \]
defined as the geometric realization of a map 
of simplicial $H\mZ\rscript{c}$--bimodules 
\[ n \ \mapsto \left[ \ H\mZ\rscript{c} \, \sm \, 
\underbrace{(H\mZ\rscript{c}\times G^!) \, \sm \dots \, \sm \, (H\mZ\rscript{c}\times G^!)}_{n}  \, 
\sm \, H\mZ\rscript{c} \ \to \ 
\underbrace{G^! \times \cdots \times G^!}_n  \ \right] \]
whose $i$-th projection, for $1\leq i\leq n$, is given symbolically by
\[ x_0 \, \sm \, (x_1,g_1) \, \sm \, \dots \, \sm \, 
(x_n,g_n) \, \sm \, x_{n+1} \quad \mapsto 
\quad x_0 \cdot x_1 \cdots x_{i-1} \cdot g_i \cdot x_{i+1} \cdots x_n 
\cdot x_{n+1} \ . \]

We define a map $\eta_Q\co  \widetilde\B(Q) \to  (B\Lc Q)^!$ as the composite 
\[  \widetilde\B(Q) \ \to \ 
\widetilde\B(H\mZ\rscript{c}\times(\Lc Q)^!)
\varr{1.3cm}{\tau_{\Lc Q}} \ (B\Lc Q)^! \ ; \]
the first map is induced by the Gamma-ring map
$Q\to H\mZ\rscript{c}\times(\Lc Q)^!$, 
which in turn comes from the structure map $Q\to H\mZ\rscript{c}$ 
and the adjunction unit $Q\to H\mZ\times(\Lc Q)^!$.
For the proof of Theorem \ref{thm-adjoint of universal} 
we are only interested in the special case $Q=H\mZ\rscript{c}$, 
but to establish that $\eta_{H\mZ^c}$
is a weak equivalence we will use a resolution argument which requires 
the general case of an arbitrary Gamma-ring over $H\mZ\rscript{c}$.
\end{numbered paragraph}

\begin{theorem} \label{Thm A} The map of $H\mZ\rscript{c}$--bimodules
\[ \eta_Q \,\co\,\widetilde\B(Q) \ \to \ (B\Lc(Q))^!   \] 
is a stable equivalence for every cofibrant Gamma-ring
$Q$ over $H\mZ\rscript{c}$.
\end{theorem}
\begin{proof}
We first assume that $Q$ is free, ie, that 
\[ Q\ = \ TX \ = \ \bigvee_{n\geq 0} \ X^{\sm n} \] 
is the tensor algebra generated by a cofibrant 
$\Gamma$--space $X$ over $H\mZ\rscript{c}$. Then the $H\mZ\rscript{c}$--bimodule 
$\widetilde\B(TX)$ can be analyzed 
through a combinatorial filtration as follows. 
For $p\geq 0$ we define $F_p\B$ as the realization of a simplicial 
sub-$H\mZ\rscript{c}$--bimodule of the bar construction 
$\B(TX)$. In simplicial degree $n$ we set
\begin{eqnarray*} (F_p\B)_n & = & 
H\mZ\rscript{c}\,\sm\,(\bigvee_{i_1+\cdots + i_n\leq p}
 X^{\sm i_1}\sm\dots\sm\, X^{\sm i_n}) \,\sm\, H\mZ\rscript{c} \\
& \subseteq & H\mZ\rscript{c}\, \sm \, (TX)^{\sm\ n}  
\sm \, H\mZ\rscript{c} \, = \, \B(TX)_n \ .
\end{eqnarray*}
The 0-th filtration is  $H\mZ\rscript{c}\,\sm\,H\mZ\rscript{c}$ 
and the subquotient $F_1\B/F_0\B$ is isomorphic to the suspension of 
$H\mZ\rscript{c}\,\sm\, X\,\sm\, H\mZ\rscript{c}$. 
To identify the subquotients of the filtration 
we use certain simplicial sets $D_p$. 
We define $D_p$ as the quotient of the standard simplicial 
$p$--simplex with the union of the first and last $(p-1)$--face collapsed 
to a basepoint. Then $D_1=S^1$ is the simplicial circle 
and $D_p$ is weakly contractible for $p\geq 2$.

We note that the $p$-th subquotient is the realization of a simplicial object 
which in dimension $n$ is of the form
\[ (F_p\B/F_{p-1}\B)_n \ = \ 
H\mZ\rscript{c}\,\sm\,(\bigvee_{i_1+\cdots + i_n = p} 
X^{\sm i_1}\sm\dots\sm\, X^{\sm i_n}) \,\sm \,H\mZ\rscript{c}  \ . \]
The map 
\[ H\mZ\rscript{c}\,\sm\, X^{\sm p} \,\sm \,H\mZ\rscript{c}
\ \to \ (F_p\B/F_{p-1}\B)_p \]
indexed by the $p$--tuple $(1,1,\dots,1)$ yields a map
\[ H\mZ\rscript{c}\,\sm\, X^{\sm p} \,\sm \,H\mZ\rscript{c}\,\sm\, \Delta^p
\ \to \ F_p\B/F_{p-1}\B \]
which factors over an isomorphism  between 
$H\mZ\rscript{c}\,\sm\, X^{\sm p}\sm \, H\mZ\rscript{c}\,\sm\, D_p$ and 
$F_p\B/F_{p-1}\B$.
Since $D_p$ is weakly contractible for $p\geq 2$, 
all the filtration subquotients  $F_p\B/F_{p-1}\B$ 
are stably contractible for $p\geq 2$. So the inclusion  
\[ H\mZ\rscript{c} \, \sm \, \Sigma X \, \sm \, 
H\mZ\rscript{c} \ \iso \ F_1\B/F_0\B \ \varr{1cm}{i} \ 
\B(TX)/F_0\B  \ = \  \widetilde\B(TX) \]
is a stable equivalence of $H\mZ\rscript{c}$--bimodules. 
To complete the verification that $\eta_{TX}$ is a stable equivalence
it remains to show that the composite
\[ \eta_{TX}\circ i \,\co\,H\mZ\rscript{c} \, \sm \, \Sigma X \, \sm \, 
H\mZ\rscript{c} \ \to \ (B\Lc(TX))^! \]
is a stable equivalence. 

We can rewrite the target $(B\Lc(TX))^!$ in a more familiar form.
Let $\Phi$ denote the forgetful functor from the category of
finitely generated free abelian groups to the category of pointed sets,
and let $\widetilde{\mZ}$ denote the reduced free functor from  
the category of pointed sets to the category of all abelian groups.
By composing the $\Gamma$--space $X$ with these two functor we obtain
an object $\widetilde{\mZ}\circ X\circ\Phi$ of the category $s\F$.
The various adjunctions show that  $\widetilde{\mZ}\circ X\circ\Phi$
and $\Lc(TX)$ represent the same functor, namely the one which sends
an object $G\in s\F$ to the set of $\Gamma$--space from 
$X$ to the underlying $\Gamma$--space of $G^!$.
Hence  $\Lc(TX)$ is isomorphic to $\widetilde{\mZ}\circ X\circ\Phi$ 
in the category $s\F$. Since the free functor $\widetilde{\mZ}$
takes suspension of simplicial sets to bar construction of simplicial
abelian groups,  $B\Lc(TX)$ is then isomorphic to 
$\widetilde{\mZ}\circ \Sigma X\circ\Phi$ in the category $s\F$.

Under the isomorphism 
$B\Lc(TX)\iso\widetilde{\mZ}\circ \Sigma X\circ\Phi$ the map 
$\eta_{TX}\circ i$ corresponds to the composite
\begin{align*}  H\mZ\rscript{c}  \sm  \Sigma X \sm 
H\mZ\rscript{c} \ \to  \
H\mZ\rscript{c} \circ  \Sigma X \circ 
&H\mZ\rscript{c}
\ \to\  
H\mZ \circ \Sigma X \circ  H\mZ  \\
& = \  (\widetilde{\mZ}\circ \Sigma X \circ\Phi)^ ! \ \iso \ (B\Lc(TX))^! \ . 
\end{align*} 
The left map is the assembly map (\ref{assembly map}), 
which is a stable equivalence by \cite[5.23]{Lydakis:Smash} 
since $H\mZ\rscript{c}$ and $\Sigma X$ are cofibrant as $\Gamma$--spaces.
The second map is a weak equivalence since the composition product
of $\Gamma$--spaces preserves stable equivalences 
(Theorem \ref{assembly properties} (a)).

The general case is proved by resolving an arbitrary cofibrant Gamma-ring 
by free Gamma-rings as follows. 
If $R$ is a simplicial object in the category of Gamma-rings, we
denote by $|R|_{\GR}$ its geometric realization  
\cite[VII 3.1]{Goerss-Jardine:book}
in the category of Gamma-rings, ie, the coend \cite[IX.6]{MacL:categories}
\[ |R|_{\GR} \ = \ \int_{n\in\Delta} \ R_n \times_{\GR} \Delta^n \ ; \] 
here $(-)\times_{\GR}\Delta^n $ refers to the enrichment 
of the category of Gamma-rings over simplicial sets, 
which has to be distinguished from the objectwise smash product 
of the underlying $\Gamma$--space with $\Delta^n_+$.  

{\bf Claim}\qua Let $R$ be a simplicial object in the category of
Gamma-rings over $H\mZ\rscript{c}$ such that for all $n\geq 0$ the map
$\eta_{R_n}\co \widetilde\B(R_n) \to (B\Lc(R_n))^!$ is a stable equivalence.
Then the map $\eta_{|R|_{\GR}}\co \widetilde\B(|R|_{\GR}) \to (B\Lc(|R|_{\GR}))^!$
is also a stable equivalence.

To prove the claim it suffices to show 
that the map is a stable equivalence of underlying $\Gamma$--spaces. 
We consider the commutative square
\[\xymatrix@C=15mm{ |\widetilde\B(R)|  \ar_{|\eta_{R}|}[d] \ar^-{\sim}[r] &
\widetilde\B(|R|_{\GR}) \ar^{\eta_{|R|_{\GR}}}[d]  \\
|(B\Lc R)^!| \ar_-{\iso}[r] & (B\Lc|R|_{\GR})^! }\]
On the left the functors $\widetilde\B$ and $(B\Lc -)^!$ are applied
dimensionwise to the simplicial Gamma-ring $R$, and then we form 
the realization of the underlying simplicial $\Gamma$--space.
The two horizontal maps are isomorphisms,
so we may show that the left vertical map is a stable equivalence.
A map of simplicial $\Gamma$--spaces
which is dimensionwise a stable equivalence becomes a stable equivalence
after realization. So the left vertical map in the above square is
a stable equivalence of underlying $\Gamma$--spaces, which proves the claim.

We apply the claim to the cotriple resolution \cite[VII.6]{MacL:categories}
of a given cofibrant Gamma-ring $Q$ over $H\mZ\rscript{c}$. 
The tensor algebra functor $T$ from the category of $\Gamma$--spaces 
to the category of Gamma-rings is left adjoint to the forgetful functor. 
The adjunction gives rise to a cotriple, hence to a simplicial Gamma-ring 
$R$ which in simplicial dimension $n$ consists of the Gamma-ring 
$R_n=T^{n+1}Q$. 
We claim that this simplicial Gamma-ring $R$ is cofibrant in the
Reedy model structure 
(\cite[Thm.\ A]{Reedy:simplicial}, \cite[5.2.5]{Hovey:book}, 
\cite[VII 2.1]{Goerss-Jardine:book}). Indeed, 
the maps from the latching objects to the levels of the simplicial Gamma-ring 
$R$ are freely generated by a wedge summand inclusion of $\Gamma$--spaces 
whose cokernel is a wedge of smash powers of $Q$. 
Since $Q$ is cofibrant as a Gamma-ring, it is cofibrant as a $\Gamma$--space
\cite[4,1 (3)]{SS:monoidal}, hence so are its smash powers.
So the maps from the latching objects to the levels of the $R$ are cofibrations
of Gamma-rings, ie, $R$ is Reedy cofibrant.
In particular, for all $n\geq 0$ the underlying $\Gamma$--space of $T^nQ$ 
is cofibrant, and hence for 
$R_n=T(T^nQ)$ the map $\eta_{R_n}\co \widetilde\B(R_n) \to (B\Lc(R_n))^!$ 
is a stable equivalence by the first part of this proof.
By the claim, the map 
$\eta_{|R|_{\GR}}\co \widetilde\B(|R|_{\GR}) \to (B\Lc(|R|_{\GR}))^!$
is also a stable equivalence.

The cotriple resolution comes with an augmentation $|R|_{\GR}\to Q$. 
After forgetting the multiplication, the augmented simplicial $\Gamma$--space
$|R|_{\GR}\to Q$ has an extra degeneracy, so the  augmentation map 
$|R|_{\GR}\to Q$ is an objectwise equivalence of Gamma-rings 
\cite[III 5.1]{Goerss-Jardine:book}.
Since $R$ is Reedy cofibrant, the realization $|R|_{\GR}$ 
is cofibrant \cite[VII 3.6]{Goerss-Jardine:book}. 
Since the functors $\widetilde\B$ and $(B\Lc -)^!$ 
both preserve stable equivalences between cofibrant Gamma-rings, 
$\eta_Q$ is a stable equivalence as claimed.
\end{proof}

\begin{numbered paragraph} 
{\bf Another bar construction}\qua
The lower left hand corner of the square (\ref{bar square})
arises from a simplicial functor $\widetilde\B(P)\in s\F$ 
which is another reduced bar construction.
We note that the category $\F$ of reduced functors from finitely 
generated free abelian groups to abelian abelian groups has a monoidal
composition product $\circ$ with unit the inclusion functor $I$; 
before composing two functors, one
of them has to be left Kan extended \cite[X.3]{MacL:categories} from 
finitely generated free to all abelian groups.

The functor $P\in\F$ is the composite of the forgetful functor from abelian
groups to pointed sets with its adjoint free functor. Hence $P$ has
the structure of a cotriple on the category of abelian groups.
This cotriple gives rise to a simplicial object $\B(P)$,
augmented over the functor $I$, 
which in simplicial dimension $n$ is given by
\[ \B(P)_n \ = \ \underbrace{P\circ \cdots \circ P}_{n+1} \ , \] 
compare \cite[VII.6]{MacL:categories}.
The augmentation  $\B(P)_0 = P\to I$ is given by the evaluation map
$\epsilon$. 

The $H\mZ$--bimodule $\B(P)_n^!=(P^{\circ (n+1)})^!$ 
is equal to the $(n+2)$--fold
composition product of the Eilenberg--MacLane $\Gamma$--space $H\mZ$.
So the assembly map (\ref{assembly map})
induces a map of simplicial $H\mZ\rscript{c}$--bimodules
\[ \B(H\mZ\rscript{c})_n \ = \ (H\mZ\rscript{c})^{\sm (n+2)} \ 
\varr{1.5cm}{\text{assembly}} \  H\mZ^{\circ (n+2)} 
\ = \  \B(P)_n^! \ . \]
We denote by  $\widetilde\B(P)$ the simplicial
functor obtained from $\B(P)$ by collapsing the simplicial 0-skeleton.
The assembly map passes to a map  
$\widetilde\B(H\mZ\rscript{c}) \to \widetilde\B(P)^!$
on quotients by the respective simplicial 0-skeleta.
\end{numbered paragraph}

\begin{theorem} \label{Thm B} 
The assembly map
\[ \widetilde\B(H\mZ\rscript{c}) \ \varr{1.5cm}{} \ \widetilde\B(P)^! \]
is a stable equivalence of $H\mZ\rscript{c}$--bimodules.
\end{theorem}
\begin{proof}
Since $H\mZ\rscript{c}$ is cofibrant as a $\Gamma$--space, the assembly
map from a smash power of a certain number of copies of 
$H\mZ\rscript{c}$ to the composition power of the same  number of copies of 
$H\mZ$ is a stable equivalence by \cite[5.23]{Lydakis:Smash} 
and Theorem \ref{assembly properties}~(a).
A map of simplicial $\Gamma$--spaces which is dimensionwise a stable
equivalence induces a stable equivalence on realizations.
So the assembly maps 
\[ \B(H\mZ\rscript{c})\ \to \ \B(P)^! \mbox{\qquad and \qquad}
\widetilde\B(H\mZ\rscript{c})\ \to \ \widetilde\B(P)^!  \] 
on realizations are stable equivalences of $H\mZ\rscript{c}$--bimodules.
\end{proof}

The lower horizontal map in the square \eqref{bar square}
arises from an objectwise equivalence
\[ u \,\co\, \B(P) \ \to \ P\oplus_J EJ \] 
of simplicial functors by passage to quotient and application of the
$(-)^!$--con\-struction.
The target $ P\oplus_J EJ$ is the simplicial functor defined as the pushout
of the diagram
\[ P \ \varl{1cm}{\iota} \ J \ \to 
\ \widetilde{\mZ}[\Delta^1]\tensor J \ .  \] 
The map $\iota\co J\to P$ is the inclusion and 
$J\to \widetilde{\mZ}[\Delta^1]\tensor J$ is induced by the inclusion of the
non-basepoint vertex of $\Delta^1$.

We describe the map 
\[ u_n \,\co\,\B(P)_n = P^{\circ (n+1)}
\ \to \  (P\oplus_J EJ)_n \ = \ P \times J^n \]
in simplicial dimension $n$ by giving the components of the various factors
of the target. The projection of $u_n$ to the first factor is the map
\[ \epsilon^{\circ n}\circ 1 \,\co\, P^{\circ (n+1)} \ \to \ P \ .  \] 
The projection of $u_n$ to the $i$-th factor of $J$, for $1\leq i\leq n$
is the map
\[   \epsilon^{\circ (n-i)}\circ 1 \circ \epsilon^{\circ i}  
\ - \ \epsilon^{\circ (n-i+1)}\circ 1 \circ \epsilon^{\circ (i-1)} \,\co\, 
P^{\circ (n+1)} \ \to \ J \ ; \] 
the target of each of the two summands is really the functor $P$, 
but the difference is annihilated by $\epsilon\co P\to I$, so it lands in
$J=\mbox{kernel}(\epsilon)$. 
Note that $P$ is the functor of 0-simplices in  $P\oplus_J EJ$,
and the quotient of  $P\oplus_J EJ$ by $P$ is isomorphic 
to the simplicial functor $BJ$. So the map $u$ passes to quotients
and yields a map of simplicial functors 
$\widetilde u\co  \widetilde\B(P) \to BJ$.

\begin{theorem} \label{Thm C} 
The map 
\[ \widetilde u \,\co\,  \widetilde\B(P) \ \to \ BJ \] 
is an objectwise weak equivalence of simplicial functors.
\end{theorem}
\begin{proof}
The simplicial subfunctor $\widetilde{\mZ}[\Delta^1]\tensor J$
of $P\oplus_J EJ$ is objectwise weakly contractible. So the quotient map 
\[ q \,\co\,P\oplus_J EJ \ \to \  
(P\oplus_J EJ) \, / \, (\widetilde{\mZ}[\Delta^1]\tensor J) \ \iso  \ I \]
is an objectwise weak equivalence of simplicial functors.

The composite map $q\circ u\co \B(P)\to I$ is the augmentation of the cotriple
resolution. Whether or not it is an objectwise weak equivalence
can be checked by looking at the augmented simplicial $\Gamma$--space
$\B(P)^!\to I^!$. However this augmented simplicial $\Gamma$--space
has an extra degeneracy, so the augmentation $\B(P)^!\to I^!$,
and hence the map $u$, is an objectwise weak equivalence 
\cite[III 5.1]{Goerss-Jardine:book}.

The simplicial 0-skeleta of $\B(P)$ and $P\oplus_J EJ$
are both equal to the functor $P$. So if we collapse the 0-skeleta,
then the induced map of quotients $\widetilde u\co  \widetilde\B(P) \to BJ$
is also an  objectwise weak equivalence.
\end{proof}

\section{A homological criterion} \label{der and hyper}

In this section we give a homological condition, 
Theorem \ref{homological criterion} below, for when the map
\[  \lambda_G \,\co\,\Z(G) \  \to \
\map_{\mathcal{GR}}(H\mZ\rscript{c},H\mZ\times G_{st}^!)_{hG_{st}(\mZ)} \] 
defined in (\ref{lambda}) is a weak equivalence. 
Here $\Z(G)$ denotes the groupoid of symmetric 2--cocycles of the functor
$G$ (\ref{def-cocycles}). $G_{st}$ is the Dold--Puppe stabilization 
of $G$, a simplicial functor which corresponds to the cubical
construction $QG$ under the Dold--Kan equivalence between simplicial objects and
non-negative complexes in the category $\F$, compare (\ref{DP and Q}).
In the next section we verify the criterion of 
Theorem \ref{homological criterion} in the case of the symmetric power
functors $B\tensor {\mathcal S}^k$.

As before, $\F$ denotes the abelian category 
of reduced functors from finitely generated free abelian groups 
to all abelian groups and $I\in\F$ is the inclusion functor.
The functor category $\F$ is abelian and  exactness 
can be checked objectwise; $\F$ has enough projectives and injectives. 
The map $G\to QG$ is the inclusion as the object in dimension zero.
Also, $\hExt^*_{\F}(I,-)$ denotes hyper-Ext groups of the inclusion functor $I$
with coefficients in a chain complex of functors, ie, the graded abelian
group of maps out of $I$ in the derived category $\D^+(\F)$ of bounded below
complexes of functors. A priori, these hyper-Ext groups can be non-trivial
in negative dimensions.

\begin{theorem} \label{homological criterion}
Let $G\in \F$ be a functor such that for 
all integers $m\leq 2$ the map
\[ \Ext^m_{\F}(I,G) \ \to \  \hExt^m_{\F}(I,QG) \]
is an isomorphism. Then the map 
\[  \lambda_G \,\co\,\Z(G) \  \to \
\map_{\mathcal{GR}}(H\mZ\rscript{c},H\mZ\times G_{st}^!)_{hG_{st}(\mZ)} \] 
is a weak equivalence of simplicial sets.
\end{theorem}

\begin{remark} \label{MacLane vs TH}
By \cite[6.1]{Sch:stable}, the hyper-cohomology
groups $ \hExt^m_{\F}(I,QG)$ are isomorphic
to the topological Hochschild cohomology groups of $H\mZ$ with coefficients
in the bimodule $G^!$. 

On the other hand, if $A$ is an abelian group, then a theorem of Jibladze and
Pirashvili \cite[Thm.\ A]{JP} identifies the cohomology groups
$\Ext^*_{\F}(I,A\tensor -)$ with the MacLane cohomology groups of $A$
\cite{MacL:ML-homology}. 
Because of this, for an arbitrary functor $G\in\F$ the groups 
$\Ext^*_{\F}(I,G)$ are sometimes referred to as the MacLane cohomology groups
of $\mZ$ with coefficients in the functor $G$. So the criterion 
of Theorem \ref{homological criterion} ask whether the natural map
from the MacLane cohomology to the topological Hochschild cohomology
of the functor $G$ is an isomorphism.

A theorem of Pirashvili and Waldhausen \cite[3.2]{Pira-Wald} 
says that if cohomology is replaced by {\em homology}, then the MacLane theory
coincides with the topological Hochschild theory for arbitrary 
coefficient functors. By \cite[6.7]{Sch:stable}, the cohomological theories
also agree if the coefficient functor is {\em additive}.

However, as the following example shows, the hypothesis of 
Theorem \ref{homological criterion} is not satisfied 
for an arbitrary functor $G\in\F$.
So for general coefficients, MacLane cohomology 
and topological Hochschild cohomology do {\em not} coincide.
\end{remark}

\begin{example} \label{counterexample} 
We give an example of a functor for which the hypothesis of Theorem
\ref{homological criterion} fails.
For a fixed prime $p$ the Frobenius maps 
\[  {\mZ/p}\tensor {\mathcal S}^{p^{h-1}} \ \to \  
{\mZ/p}\tensor {\mathcal S}^{p^h} \ ,
\quad  a\tensor x \ \longmapsto \ a\tensor x^p \]
define a morphism in the category $\F$. 
We consider the functor 
\[ G \ = \ \colim_h \ {\mZ/p}\tensor {\mathcal S}^{p^h} \]
defined as the colimit of the sequence of Frobenius maps.
Since the Frobenius transformations are injective, the natural map 
\[ I\ \varr{2cm}{\text{projection}} \  \mZ/p\tensor I = 
{\mZ/p}\tensor {\mathcal S}^{p^0} \ \varr{2cm}{\text{inclusion}} \ G \] 
is a non-trivial element of $\Hom_{\F}(I,G)$. 
On the other hand, the stable derived functors 
of the symmetric power functor $A\tensor {\mathcal S}^k$
are trivial up to dimension $2k-3$ \cite[12.3]{Dold-Puppe},
so $Q({\mZ/p}\tensor {\mathcal S}^{p^h})$ has trivial homology 
up to dimension \mbox{$(2p^h-3)$.}
Since the $Q$--construction and homology commute with filtered colimits, 
the complex $QG$ is acyclic, so the hyper-cohomology groups 
$\hExt^*_{\F}(I,QG)$ are trivial.
In particular the map
\[ \mbox{Ext}^0_{\F}(I,G) \ \to \ \hExt^0_{\F}(I,QG) \]
is not injective.
\end{example}

\begin{proof}[Proof of Theorem \ref{homological criterion}] 
The map $\lambda_G$ (\ref{lambda}) is obtained from the commutative square
of simplicial abelian groups
\[\xymatrix{ G(\mZ) \ar_{\Theta}[dd] \ar[r] & G_{st}(\mZ) \ar^{\Theta}[d] \\
& \mbox{Z}^2_s(G_{st}) \ar[d] \\
\mbox{Z}^2_s(G) \ar[r] & \map_{\GR}(H\mZ\rscript{c},H\mZ\times G_{st}^!) } \]
by passage to vertical homotopy cofibres in the category 
of simplicial abelian groups.
So it suffices to show that the square is homotopy cocartesian 
(in the category of simplicial abelian groups).
Evaluation at $\mZ$ is represented by the projective functor $P$,
symmetric 2--cocycles are represented by the functor $J$ 
(Lemma \ref{lem-universal cocycle}), and the split extension construction
(\ref{split extension}) has a left adjoint $\Lc$ (Lemma \ref{quillen pair}).
So the square is isomorphic to the square
\[\xymatrix{ \Hom_{\F}(P,G) \ar_{\Hom(\iota,G)}[d] \ar[r] & 
\map_{s\F}(P,G_{st}) \ar[d] \\
\Hom_{\F}(J,G) \ar[r] & \map_{s\F}(\Lc(H\mZ\rscript{c}),G_{st}) }\]
where $\iota\co J\to P$ is the inclusion.

The map $\delta\co \Lc(H\mZ\rscript{c})\to J$ (\ref{delta}) which is adjoint
to the universal derivation is a stable equivalence
by Theorem \ref{thm-adjoint of universal}. 
We let $\alpha\co J\rscript{c}\to J$ be a cofibrant approximation 
of the functor $J$ in the strict model structure of simplicial
functors where the weak equivalences are defined objectwise;
equivalently, the normalized chain complex of $J\rscript{c}$ is a
projective resolution of $J$.
Then $\delta\co \Lc(H\mZ\rscript{c})\to J$ lifts to a map of simplicial
functors $\bar\delta\co \Lc(H\mZ\rscript{c})\to J\rscript{c}$. 
The lift $\bar\delta$ is
then a stable equivalence between cofibrant simplicial functors. 
Since $G_{st}$ is homotopy additive, alias stably fibrant,
the map $\bar\delta$ induces a weak equivalence of simplicial abelian groups 
upon application of $\map_{s\F}(-,G_{st})$. 
Since $G$ is a constant simplicial functor, the map 
\[ \map_{s\F}(\alpha,G) \,\co\, \Hom_{\F}(J,G) = \map_{s\F}(J,G)
\ \to \ \map_{s\F}(J\rscript{c},G) \]
is an isomorphism.
In other words, it suffices to show that
the square in the category of simplicial abelian groups
\[\xymatrix{ \map_{s\F}(P,G) \ar_{\map(\iota\circ\alpha,G)}[d] \ar[r] & 
\map_{s\F}(P,G_{st}) \ar^{\map(\iota\circ\alpha,G_{st})}[d] \\
\map_{s\F}(J\rscript{c},G) \ar[r] & \map_{s\F}(J\rscript{c},G_{st}) }\]
is homotopy cocartesian. For this in turn it is enough to show that the
map on horizontal homotopy cofibres 
\[ \map_{s\F}(\iota\circ\alpha,G_{st}/G) \,\co\,
\map_{s\F}(P,G_{st}/G) \ \to \ \map_{s\F}(J\rscript{c},G_{st}/G) \]
is a weak equivalence where $G_{st}/G$ denotes to cofibre of the stabilization
map.

If $K$ and $K'$ are two simplicial functors such that $K$ is cofibrant,
then the Dold--Kan theorem provides a natural isomorphisms of groups
\[  \pi_n \, \map_{s\F}(K,K') \ \iso \  [NK[n],NK'] \] 
for $n\geq 0$, where $N$ is the normalized chain complex 
and $[-,-]$ denotes maps in the derived category $\D^+(\F)$ 
of bounded below chain complexes of functors.
The normalized chain complex of $J\rscript{c}$ is
quasi-isomorphic to $J$, and the normalized chain complex of $G_{st}$ is the
cubical construction $QG$. So we need to show that the map
\[ [P[n],QG/G] \ \to \ [J[n],QG/G] \]
is an isomorphism for $n\geq 0$.
The short exact sequence of functors $J\to P\to I$ yields a long exact 
sequence after applying $[-,QG/G]$, so it is enough to
show that the groups 
\[ [I[n],QG/G] \ = \ \hExt_{\F}^{-n}(I,QG/G) \]
vanish for $n\geq -1$. This in turn follows from the assumption that
the map
\[ \Ext^m_{\F}(I,G) \ \to \  \hExt^m_{\F}(I,QG) \]
is an isomorphism for all integers $m\leq 2$.
\end{proof}

\section{Cohomology of symmetric power functors} \label{sec-Ext}

The purpose of this section is to prove that the symmetric power functors
satisfy the homological criterion of Theorem \ref{homological criterion};
this completes the proof of the main theorem.
\begin{theorem} \label{finite functors}
For all $m\in\mZ$, all $k\geq 1$ and all abelian groups $A$ the map 
\[ \Ext^m_{\F}(I,A\tensor {\mathcal S}^{k}) \ \to \  
\hExt^m_{\F}(I,Q(A\tensor {\mathcal S}^{k})) \]
is an isomorphism.
\end{theorem}

\begin{remark} 
Theorem \ref{finite functors} can be interpreted as saying 
that MacLane cohomology coincides with topological Hochschild cohomology 
for the symmetric power functors, compare Remark \ref{MacLane vs TH}.
The groups $\Ext^m_{\F}(I,A\tensor {\mathcal S}^k)$ have been calculated for 
$A=\mZ/p$ and for $A=\mZ$ by Franjou-Lannes-Schwartz and Franjou-Pirashvili,
see \cite[Thm.\ 6.6 and Prop.\ 9.1]{Franjou-Lannes-Schwartz:MacLane(F_p)}
and \cite[2.1]{Franjou-Pira:MacLane(Z)}.
So Theorem \ref{finite functors} is a calculation of the
topological Hochschild cohomology groups of $H\mZ$ with coefficients
in the bimodule $(A\tensor {\mathcal S}^k)^!$.
In particular, the topological Hochschild cohomology groups  of $H\mZ$ 
with coefficients in $(A\tensor {\mathcal S}^k)^!$ are trivial in negative dimensions.

Our proof of Theorem \ref{finite functors} is not completely satisfactory
because it uses the explicit calculations of the groups 
$\Ext^m_{\F}(I,{\mathcal S}^{k})$; these enter as the sparseness hypothesis (c)
of Theorem \ref{deduction theorem} below. It would be desirable to
have a direct proof of Theorem \ref{finite functors} which 
would hopefully shed more light on the question for which coefficient functors
MacLane cohomology coincides with topological Hochschild cohomology.
Example \ref{counterexample} shows that some restriction on the functor 
has to be imposed.
\end{remark}

Theorem \ref{finite functors} is a special case of the following
Theorem \ref{deduction theorem}. 
To apply it, we choose a projective resolution $P_*\to I$ of the functor $I$ 
in the abelian category $\F$. 
Then we let $T\co \Ch^+(\F)\to \coCh$ 
be the homomorphism complex out of this resolution,
\[ T(X) \ = \ \Hom^{\bullet}_{\F}(P_*,X)  \ . \]
So $T(X)$ is a (usually unbounded) cochain complex of abelian groups 
and as a functor of $X$ it is additive, exact, 
and preserves inverse limits and quasi-isomorphisms. 
The cohomology groups of $\Hom^{\bullet}_{\F}(P_*,X)$ are
the hyper-co\-homology groups $\hExt^*_{\F}(I,X)$.
By a theorem of Pirashvili (\cite[2.15]{Pira:additivization}, 
see also \cite[0.4]{Franjou-Lannes-Schwartz:MacLane(F_p)} 
or the appendix of \cite{Betley-Pira:stableK}), 
the extension groups $\Ext^*_{\F}(I,-)$ vanish for every  
diagonalizable functor (\ref{DP and Q} (e)), 
so $\Hom^{\bullet}_{\F}(P_*,-)$ takes
diagonalizable functors to acyclic complexes.
The sparseness condition (c) is proved in 
\cite[Prop.\ 2.1]{Franjou-Pira:MacLane(Z)}.

The homomorphism complex  $\Hom^{\bullet}_{\F}(P_*,-)$ is the only functor
to which we apply Theorem \ref{deduction theorem}; nevertheless we state
and prove it in the general form because we think it makes the proof
more understandable.

\begin{theorem} \label{deduction theorem}
Let $T\co \Ch^+(\F)\to \coCh$
be a functor from the category of bounded below chain complexes in the
abelian category $\F$ to the category of (not necessarily bounded)
cochain complexes of abelian groups. Suppose furthermore that
\renewcommand{\labelenumi}{\rm(\alph{enumi})}
\begin{enumerate}
\item $T$ is additive, exact and preserves inverse limits 
and quasi-isomorphisms,
\item the complex $T(D)$ is acyclic 
for every diagonalizable functor $D\in\F$ (\ref{DP and Q} (e)), 
considered as a complex concentrated in dimension 0, and
\item $T$ is {\em sparse on symmetric powers} in the sense that for all 
$k\geq 1$
the cohomology of the complex $T({\mathcal S}^k)$ is concentrated 
in dimensions congruent to 1 modulo $2k$.
\end{enumerate}
Then for all abelian groups $A$, and all $k\geq 1$
the natural map
\[ T(A\tensor {\mathcal S}^k) \ \to  \ T(Q(A\tensor {\mathcal S}^k)) \]
is a quasi-isomorphism.
\end{theorem}

\begin{remark}
The heart of Theorem \ref{deduction theorem} is a convergence issue, 
or a question to what extent the functor $T$ commutes with infinite sums
(up to quasi-isomorphism). 
Indeed, for any functor $G\in\F$, the cokernel of the stabilization 
map $G\to QG$ is a bounded below complex of diagonalizable functors
\cite[7.4]{Johnson-McCarthy},
usually with non-trivial homology in arbitrarily high dimensions. 
So the cokernel $QG/G$ can be written as the colimit of a
sequence of bounded complex of diagonalizable functors. 
So if a functor $T$ as in Theorem \ref{deduction theorem} 
commutes with filtered colimits or infinite sums, then properties (a) and (b) 
already imply that $T(QG/G)$ is acyclic,
and so $T(G)\to T(QG)$ is a quasi-isomorphism, for all $G\in\F$. 
In the case of interest for us, namely 
$T=\Hom^{\bullet}_{\F}(P_*,-),$ the functor $T$ fails to commute with
infinite sums, essentially because the inclusion functor $I$ is not a small
(or compact) object in the derived category of $\F$.
And indeed, $\Hom^{\bullet}_{\F}(P_*,QG/G)$ fails to be
acyclic in general as Example \ref{counterexample} shows.
However the sparseness condition (c) makes it possible to obtain
the desired conclusion for the symmetric power functors. 
\end{remark}

The following observation goes back, at least, 
to Dold and Puppe \cite{Dold-Puppe}.
For every $k\geq 1$ let $d_k$ denote the greatest common divisor 
of the binomial coefficients $k \choose i$ for $1\leq i\leq k-1$. 
Then 
\[ d_k \ = \ \left\lbrace \begin{array}{ll} 
p & \mbox{ if $k=p^h$ for a prime $p$ and $h>0$,} \\
1 & \mbox{ else.} 
\end{array} \right. \]

\begin{lemma}{\em \cite[10.9]{Dold-Puppe}} \label{Dold-Puppe lemma} 
For every $k\geq 1$ and every abelian group $A$, multiplication by the number
$d_k$ on the functor $A\tensor {\mathcal S}^k$ 
factors over a diagonalizable functor (\ref{DP and Q} (e)). 
In particular, if $k$ is not a prime power, then  $A\tensor {\mathcal S}^k$
is a retract of a diagonalizable functor.
\end{lemma}
\begin{proof}
For every $1\leq i\leq k-1$ the comultiplication of the symmetric algebra
gives a map of functors
\[ \Delta_{i,k-i} \,\co\,A\tensor {\mathcal S}^k \ \to \ 
A\tensor {\mathcal S}^i\, \tensor {\mathcal S}^{k-i} \ ; \]
the explicit formula for this map is given by 
\[ a\tensor x_1\cdot \cdots \cdot x_k \ \longmapsto \ 
a\tensor \sum_{T\subset \{1,2,\dots,k\}, \ |T|= i} \quad
(\prod_{j\in T} x_j) \tensor (\prod_{j\not\in T} x_j) \ . \]
The composite of $\Delta_{i,k-i}$  
with the natural projection $A\tensor {\mathcal S}^i \tensor {\mathcal S}^{k-i} \to A\tensor {\mathcal S}^k$
is multiplication by the binomial coefficient $k\choose i$. 
So if we choose a presentation
\[ d_k \ = \ \sum_{i=1}^{k-1} \ \lambda_i \cdot {k\choose i} \]
for suitable integers $\lambda_i$ then the composition 
\[ A\tensor {\mathcal S}^k \ \varr{2.3cm}{\sum \, \lambda_i\cdot\Delta_{i,k-i}} 
\ \bigoplus_{i=1}^{k-1} \ A\tensor {\mathcal S}^i \, \tensor {\mathcal S}^{k-i} \ \to \ 
A\tensor {\mathcal S}^k \]
is multiplication by the number $d_k$.
\end{proof}

Finally, we give the {\bf proof of Theorem \ref{deduction theorem}}.
For a functor $G\in\F$ we use the notation $\bar QG$ for the quotient
complex $QG/G$. By the exactness of $T$ we then have to show that 
for all abelian groups $A$, and all $k\geq 1$
the complex $T(\bar Q(A\tensor {\mathcal S}^k))$ is acyclic.

\medskip

\noindent
{\bf Step 1\qua Diagonalizable functors} 

Suppose $D\in \F$ is diagonalizable (\ref{DP and Q} (e)). 
Then $QD$ is an acyclic complex by 
\ref{DP and Q} (e), and so $T(QD)$ is acyclic. 
By property (b) of the functor $T$ the complex $T(D)$, and hence, 
by exactness, the quotient complex $T(\bar QD)$ is also acyclic.

If $k$ is not a prime power then $A\tensor {\mathcal S}^k$ is a retract 
of a diagonalizable functor by Lemma \ref{Dold-Puppe lemma}, 
so  Theorem \ref{deduction theorem} holds for such exponents.
From now on we assume that the exponent is of the form
$k=p^h$ for a prime number $p$ and some $h\geq 0$.

\medskip 

\noindent
{\bf Step 2\qua Reduction to the case $A=\mZ/p$}

For the course of this proof we call an abelian group $A$
{\em good} if the complex 
\[ T(\bar Q(A\tensor {\mathcal S}^{p^h})) \]
is acyclic. We show that if the group $\mZ/p$ is good,
then every abelian group is good.

Multiplication by the number $p$ is an epimorphism on the functor
$\mQ/\mZ\tensor {\mathcal S}^{p^h}$ with kernel isomorphic to $\mZ/p\tensor {\mathcal S}^{p^h}$.
Since the cubical construction and the functor $T$ are exact, 
multiplication by $p$ on the complex
$T(\bar Q(\mQ/\mZ\tensor {\mathcal S}^{p^h}))$ is surjective and has kernel
isomorphic to $T(\bar Q(\mZ/p\tensor {\mathcal S}^{p^h}))$, which is acyclic
since $\mZ/p$ was assumed to be good.
So multiplication by $p$ is a quasi-isomorphism on the complex
$T(\bar Q(\mQ/\mZ\tensor {\mathcal S}^{p^h}))$.
On the other hand, multiplication by $p$ on $\mQ/\mZ\tensor {\mathcal S}^{p^h}$
factors over a diagonalizable functor $D$, say, 
by Lemma \ref{Dold-Puppe lemma}.
So multiplication by $p$ on $T(\bar Q(\mQ/\mZ\tensor {\mathcal S}^{p^h}))$
factors through the complex $T(\bar QD)$, which is acyclic by Step 1.
Since multiplication by $p$ on $T(\bar Q(\mQ/\mZ\tensor {\mathcal S}^{p^h}))$
is both a quasi-isomorphism and factors through an acyclic complex,
$T(\bar Q(\mQ/\mZ\tensor {\mathcal S}^{p^h}))$ must itself be acyclic,
and so the group $\mQ/\mZ$ is good.

The $\bar Q$--construction and the functor $T$ commute with products. 
So the product of a family of good abelian groups is again good.
In particular a product of any number of copies of the group 
$\mQ/\mZ$ is good. Every injective abelian group is a summand of a product
of copies of $\mQ/\mZ$, hence injective abelian groups are good.

If $A$ is a subgroup of an abelian group $B$, then the sequence
of functors
\[ 0 \ \to \ A\tensor {\mathcal S}^{p^h} \ \to \ B \tensor {\mathcal S}^{p^h}\ 
\to \ (B/A)\tensor {\mathcal S}^{p^h} \ \to \ 0  \]
is exact.
Since the cubical construction and the functor $T$ are also exact,
$A$ is good as soon as $B$ and $B/A$ are.
Since an arbitrary abelian group can be embedded into an 
injective abelian group with injective cokernel, every abelian group is good. 

\medskip

\noindent
{\bf Step 3\qua Reformulation in terms of exterior power functors}
\nopagebreak

Let $\Lambda^{p^h}\in\F$ denote the exterior power functor of degree $p^h$.
The Koszul complex (see \cite[4.3.1.7]{Illusie:cotangentI} or 
\cite[3.2]{Franjou-Lannes-Schwartz:MacLane(F_p)})
is an extension of length $p^h-1$ of the functor 
$\mZ/p\tensor \Lambda^{p^h}$ by the functor $\mZ/p\tensor {\mathcal S}^{p^h}$ with the
special property that all functors occurring in the extension 
are diagonalizable. By Step 1, all these intermediate terms are sent
to acyclic complexes by $T(\bar Q-)$.
Since both $T$ and $\bar Q$ are exact,
the complex $T(\bar Q(\mZ/p\tensor {\mathcal S}^{p^h}))$ is acyclic if and only if
the complex $T(\bar Q(\mZ/p\tensor \Lambda^{p^h}))$ is acyclic.

The spareness assumption \ref{deduction theorem} (c) 
on the cohomology of $T({\mathcal S}^{p^h})$ and the short exact sequence of functors
\[ 0 \ \to \ {\mathcal S}^{p^h} \ \varr{1cm}{\times p} \ {\mathcal S}^{p^h}\
\to \ \mZ/p\tensor {\mathcal S}^{p^h} \ \to \ 0 \]
imply that the cohomology of the complex $T(\mZ/p\tensor {\mathcal S}^{p^h})$ 
is concentrated in dimensions congruent to 0 or 1 mod $2p^h$. 
The existence of the Koszul complex then shows that the cohomology 
of the complex 
$T(\mZ/p\tensor \Lambda^{p^h})$ is concentrated 
in dimensions congruent to $p^h-1$ or $p^h$ mod $2p^h$.

Now we set up an induction on the exponent $h$.
For the inductive step we use a certain complex of functors 
which relates the exterior power functor of degree $p^h$ to that
of degree $p^{h-1}$.
Here and in what follows we extend the cubical construction to 
bounded below chain complexes of functors 
by applying the functor $Q$ dimensionwise and taking the total complex 
of the resulting bicomplex.
This extended $Q$--construction is still exact and preserves quasi-isomorphisms.

\medskip 

\noindent
{\bf Step 4\qua There exists a complex $C$ of functors from $\F$, 
concentrated in non-negative dimensions, 
with the following properties:}
\renewcommand{\labelenumi}{\rm(\alph{enumi})}
\begin{enumerate}
\item  In dimension zero, $C_0=\mZ/p\tensor \Lambda^{p^h}$ and 
the inclusion ${\mZ/p}\tensor\Lambda^{p^{h}}\to C$
induces quasi-isomorphisms
\[ Q({\mZ/p}\tensor\Lambda^{p^{h}})\ \to \ QC
\mbox{\qquad and \qquad} 
T({\mZ/p}\tensor\Lambda^{p^{h}})\ \to \ T(C)  \ . \] 
\item  All non-trivial homology functors of $C$ are isomorphic to
$\mZ/p\tensor\Lambda^{p^{h-1}}$.
\end{enumerate}

\medskip

We let $\F(\mF_p)$ denote the category of reduced functors 
from finitely generated $\mF_p$--vector spaces to $\mF_p$--vector spaces.
We first construct a complex in the category  $\F(\mF_p)$; 
the desired complex is then obtained by composition with 
\[ -\tensor\mF_p \,\co\,(\text{f.g.\ free abelian groups})
\ \to \  (\text{f.g.\ $\mF_p$--vector spaces}) \]
and the inclusion of $\mF_p$--vector spaces into abelian groups.

For every $\mF_p$--vector space $V$ let $L^*(V)$ denote the
quotient of the symmetric algebra on $V$ by the ideal generated by 
all $p$-th powers of elements. Then $L^*(V)$ inherits the grading from 
the symmetric algebra and we let $L^n(V)$ denote the summand of homogenous
degree $n$. If $p=2$, then the functors $L^n$ coincide with the
exterior power functors $\Lambda^n$.
By \cite[1.3.1]{Franjou:puissance} there exists a complex 
$X$ (a quotient of the deRham complex) with
\[ X_i \ = \ \left\lbrace \begin{array}{ll} L^{i}\tensor\Lambda^{p^h-i}_{\mF_p} &
\mbox{for $0\leq i\leq p^h$} \\
0 & \mbox{else.} \end{array} \right. \]
whose only non-trivial homology is in dimension $p^{h-1}(p-1)$ where we have
\[ H_{p^{h-1}(p-1)}X \ \iso  \ \Lambda^{p^{h-1}}_{\mF_p} \]
(we have reversed the grading of \cite[1.3.1]{Franjou:puissance}
so that the differential decreases the dimension and $X$ is a chain complex
as opposed to a cochain complex).
The complex $X$ is part of the complex we are looking for, and we
obtain the other part by dualization as follows.
 
The dual $DF$ of a functor $F\in\F(\mF_p)$ is defined by 
$DF(V) \ = \ F(V^{\vee})^{\vee}$, where $V^{\vee}$ refers 
to the dual vector space of $V$. $DF$ is again an object of the category 
$\F(\mF_p)$. Dualization is contravariant and exact in $F$ and it 
satisfies $D(F\tensor G)\iso DF\tensor DG$. 
The exterior power functors and the functors $L^n$ are self-dual,
ie, there are isomorphisms $D\Lambda^n_{\mF_p}\iso\Lambda^n_{\mF_p}$ 
and $DL^n\iso L^n$.
So if we dualize the complex $X$ we obtain a complex $DX$ which is
concentrated in dimensions $-p^h$ through 0, which satisfies
$(DX)_0\iso \Lambda^{p^h}_{\mF_p}$, $(DX)_{-p^h}\iso L^{p^h}$ and whose
only non-trivial homology is
\[ H_{-p^{h-1}(p-1)}DX \ \iso \ D\Lambda^{p^{h-1}}_{\mF_p} \ \iso  \
\Lambda^{p^{h-1}}_{\mF_p} \ . \]

The desired complex $C$ is now obtained by splicing infinitely many copies
of the complexes $X$ and $DX$ alternatingly 
at $L^{p^h}$ and $\Lambda^{p^h}_{\mF_p}$, and then
passing from the complex in $\F(\mF_p)$ to a complex in $\F$. 
More precisely,
\begin{eqnarray*} C_0 & = & \Lambda^{p^h}_{\mF_p} \\
C_{n(p^h-1)+i} & = & \left\lbrace \begin{array}{ll}
X_i & \mbox{ for $1\leq i\leq p^h-1$ and $n\geq 0$ even,} \\
DX_{p^h-i} & \mbox{  for $1\leq i\leq p^h-1$ and $n\geq 0$ odd.} \\
\end{array} \right.
\end{eqnarray*}

{\bf The map  $Q({\mZ/p}\tensor\Lambda^{p^{h}})\to QC$\,:}\qua
Let $C/C_0$ denote the quotient complex which is equal to
$C$ except in dimension zero, where it is trivial.
The complex $C/C_0$ is bounded below and consists entirely 
of diagonalizable functors. 
Since the $Q$--construction of a diagonalizable functor
is acyclic (\ref{DP and Q} (e)), the complex $Q(C/C_0)$ is also acyclic.
Since $Q$ is exact,
the map  $Q({\mZ/p}\tensor\Lambda^{p^{h}})\to QC$ is a quasi-isomorphism.

\medskip

{\bf The map  $T({\mZ/p}\tensor\Lambda^{p^{h}})\to T(C)$\,:}\qua
The fact that this map is a quasi-isomor\-phism is a consequence 
of the sparseness assumption  \ref{deduction theorem} (c) on $T$.
In more detail: the complex $C/C_0$ can be written as 
the inverse limit of the tower of truncated complexes $K^n$ defined by
\[ (K^n)_i \ = \ \left\lbrace \begin{array}{ll}
C_i & \mbox{ if $1\leq i\leq 2n(p^h-1)$,} \\
\mZ/p \tensor \Lambda^{p^h} & \mbox{ if $i=2n(p^h-1)+1$,}\\
0 & \mbox{ else.}
\end{array} \right. \]
The boundary maps of $K^n$ are those of $C$ and the inclusion of the cycles
\[ \mZ/p \tensor \Lambda^{p^h} \ \iso \ D(X_0) \ \to 
\ D(X_1) \ = \ C_{2n(p^h-1)} \ . \]
The chain map $K^{n+1}\to K^n$ is the identity up to dimension $2n(p^h-1)$
and the epimorphism
\[  C_{2n(p^h-1)+1}  \ = \ X_1 \ \to 
\ X_0 \ = \ \mZ/p \tensor \Lambda^{p^h} \ . \]
in dimension $2n(p^h-1)+1$.
 
For each $n$, the kernel of the projection 
$K^n\to(\mZ/p\tensor\Lambda^{p^h})[2n(p^h-1)+1]$
onto the top functor is a bounded complex of diagonalizable functors.
Since $T$ is exact and takes diagonalizable functors to acyclic complexes,
the induced map
\[ T(K^n)\ \varr{1cm}{\sim} \ T((\mZ/p\tensor\Lambda^{p^h})[2n(p^h-1)+1]) \]
is a quasi-isomorphism. 

By the sparseness assumptions and Step 3, the cohomology of 
$T(\mZ/p\tensor\Lambda^{p^h})$ is concentrated in dimensions congruent to
$p^h-1$ and $p^h$ mod $2p^h$. Hence the cohomology of 
$T(K^n)$ is concentrated in dimensions 
congruent to $-2n+p^h$ and $-2n+p^h+1$ mod $2p^h$.
Hence the surjective map $T(K^{n+1})\to T(K^n)$
induces trivial maps on cohomology groups for dimensional reasons. 
Since $T$ commutes with inverse limits, $T(C/C_0)$ is the inverse limit 
of the tower of complexes $T(K^n)$ for $n\geq 1$, so $T(C/C_0)$ is acyclic.

\begin{example} It might be instructive to describe the complex $C$ 
just constructed in the smallest non-trivial case, 
namely for $p=2$ and $h=1$.
Then $C$ is the mod 2 reduction of an `integral' complex $\widetilde C$
defined by
\[ \widetilde C_i \ = \ \left\lbrace \begin{array}{ll} 
0 & \mbox{if $i<0$,} \\
\Lambda^2 &  \mbox{if $i=0$, and} \\  
I \tensor I  &  \mbox{if $i>0$.}
\end{array} \right. \]
The differential $d_i\co \widetilde C_i\to  \widetilde C_{i-1}$ is
given by
\[ d_i(x\tensor y) \ = \ \left\lbrace \begin{array}{ll} 
x \sm y &  \mbox{if $i=1$,} \\  
x \tensor y + y \tensor x  &  \mbox{if $i\geq 2$ and $i$ is even,}\\
x \tensor y - y \tensor x  &  \mbox{if $i\geq 2$ and $i$ is odd.}
\end{array} \right. \]
The homology of the complex  $\widetilde C$ is 
given by
\[ H_i\, \widetilde C \ \iso \ \left\lbrace \begin{array}{ll} 
\mZ/2 \tensor I &  \mbox{if $i\geq 1$ and $i$ is odd,} \\  
0 &  \mbox{else.}
\end{array} \right. \]
Since the homology functors of  $\widetilde C$ are additive 
and since  $\widetilde C$ is diagonalizable in positive dimensions,
$\widetilde C$ is quasi-isomorphic to the cubical construction
$Q\Lambda^2$. Hence by \ref{DP and Q} (b) we can read off 
the Dold--Puppe stable derived functors of $\Lambda^2$ as
\[ L\rscript{st}_i\, \Lambda^2  \ \iso \  H_i\, \widetilde C(\mZ) \ \iso \ 
\left\lbrace \begin{array}{ll} 
\mZ/2 &  \mbox{if $i\geq 1$ and $i$ is odd,} \\  
0 &  \mbox{else,}
\end{array} \right. \]
compare e.g.\ \cite{simson}.
The complex $C$ constructed in Step 4 for $p=2$ and $h=1$ is isomorphic
to the reduction $\mZ/2\tensor \widetilde C$; by the universal coefficient
theorem, the homology functors of $C$ are thus isomorphic to the functor
$\mZ/2\tensor I$ in every positive dimension, and trivial otherwise.
\end{example}

\medskip 

\noindent
{\bf Step 5\qua The complex $T(\bar Q(\mZ/p\tensor \Lambda^{p^h}))$
is acyclic for all $h\geq 0$.}

We proceed by induction.
For $h=0$ we have $\mZ/p \tensor \Lambda^1=\mZ/p\tensor I$ 
which is an additive functor. Thus all cross-effect vanish and
by property \ref{DP and Q} (c) of the cubical construction, 
the complex $\bar Q(\mZ/p\tensor\Lambda^1)$ is trivial. 
Hence $T(\bar Q(\mZ/p\tensor\Lambda^1))$ is also trivial.

Now suppose that $h\geq 1$ and assume that 
$T(\bar Q(\mZ/p\tensor \Lambda^{p^{h-1}}))$ is already known 
to be acyclic.
Let $C$ be any complex as in Step 4.
We consider the commutative square of bounded below chain complexes 
of functors
\[\xymatrix{ \mZ/p\tensor\Lambda^{p^h} \ar[r] \ar[d] & C \ar[d] \\
Q(\mZ/p\tensor \Lambda^{p^h}) \ar[r] & QC }\]
where we view the functor $\mZ/p\tensor \Lambda^{p^h}$
as a complex concentrated in dimension zero and
the horizontal maps are induced by the inclusion 
$\mZ/p\tensor \Lambda^{p^h}\to C$.

By property (a) of the complex $C$ and since $T$ preserves
quasi-isomorphisms, both rows of the square induce a quasi-isomorphism 
after applying $T$.
Since $T$ is exact, the map
\[ T(\bar Q(\mZ/p\tensor \Lambda^{p^h})) \ \to \ T(\bar QC) \]
is thus a quasi-isomorphism. So in order to finish the induction step, 
it remains to show that the complex $T(\bar QC)$ is acyclic. 

This last step is where the induction hypothesis is used.
The complex $C$ is the inverse limit 
of its Postnikov tower (homological truncations) 
\[ \cdots \ \to \ P_n \ \to \ \cdots \ 
\to \ P_1 \ \to \ P_0  \ \to \ P_{-1}=0 \ . \]
In the tower each map $P_n\to P_{n-1}$ is a surjection whose kernel 
is quasi-isomorphic to the $n$-th homology 
functor of $C$ concentrated in dimension $n$. 
By property (b) of the complex $C$ all non-trivial homology functors 
are isomorphic to the exterior power functor $\mZ/p\tensor\Lambda^{p^{h-1}}$, 
for which we already know that the map
$T(\bar Q(\mZ/p\tensor\Lambda^{p^{h-1}}))$ is acyclic.
Since the $Q$--construction and $T$ are exact we conclude by induction 
that for all $n\geq 0$ the complex $T(\bar QP_n)$ is acyclic.

The Postnikov tower consists of bounded below chain complexes
and it stabilizes in each dimension. So the complex $\bar QC$ is the
inverse limit of the complexes $\bar QP_n$. 
Since $T$ commutes with the inverse limits, $T(\bar QC)$ 
is the inverse limits of the acyclic complexes $T(\bar QP_n)$.
Since $T$ and $\bar Q$ also preserve epimorphisms, this inverse limit
is acyclic.

\section{Perspectives} \label{loose ends}

We end the paper with an application of Theorem \ref{main}
which concerns an interesting homotopical property of the Gamma-ring $DB$.
Then we discuss some variants of the Construction \ref{Gamma-ring map} 
of Gamma-ring maps from formal group laws, and some possible directions
for further investigation. 

For the application we use the conjugation action to obtain an obstruction 
to the existence of $k$--algebra structures on Gamma-rings. 
With this tool we then show that the Gamma-ring $DB$ is not 
stably equivalent to the Eilenberg--MacLane Gamma-ring of any simplicial ring 
(unless $B$ is a $\mathbb Q$--algebra). 
To motivate the criterion we look at the classical case of discrete rings 
first. If $k$ is a commutative ring and $A$ any associative ring, then the 
$k$--algebra structures on $A$ correspond to the central ring maps $k\to A$. 
In particular, the unit map of every such $k$--algebra structure 
gives an element of the set Ring$(k,A)$ of ring maps which is a fixed point 
of the conjugation action of the units of $A$. 
Something similar happens for Gamma-rings. Suppose $R$ is a Gamma-ring 
which is stably equivalent to an algebra over the commutative Gamma-ring $k$. 
Any chain of equivalences to a $k$--algebra determines a homotopy class 
of Gamma-ring maps $[\eta]\in [k,R]_{\Ho\mathcal{GR}}$ 
underlying the unit map of the algebra structure.

\begin{theorem} \label{algebra structures} 
Suppose $\eta\in\mbox{\em Ring}(k,R)$ is a Gamma-ring map whose homotopy class 
underlies a $k$--algebra structure of $R$. Then the conjugation action map
\[ UR^{\times} \ \to \ \mbox{\em Ring}(k,R) \quad ,
\quad u \ \mapsto \ u\cdot \eta \cdot u^{-1} \]
is null-homotopic. So if the conjugation action map is essential 
for every component of the space $\mbox{\em Ring}(k,R)$ of Gamma-ring maps, 
then $R$ is not stably equivalent to any $k$--algebra.
\end{theorem}
\begin{proof} We can assume that $R$ is itself a stably fibrant $k$--algebra. 
Then Construction \ref{conjugation construction} can be done 
in the category of $k$--algebras, as opposed to Gamma-rings, relative
to the trivial group (Lemma \ref{flat pushout proof} is also valid
in the category of $k$--algebras). We obtain a diagram of $k$--algebras
\[\xymatrix{ k \ar[d] & k \ar[d] \ar@{=}[l] \ar[r] & k[UR^{\times}] \ar[d] \\
R & R_1 \ar_-{\sim}[r] \ar@{->>}^-{\sim}[l] & \ R_3\ . }\]
Ignoring the $k$--algebra structure we can use the objects in this diagram 
to model the space of Gamma-ring maps from $k$ to $R$. More precisely, 
the space hom$_{\mathcal{GR}}(k\rscript{c},R_3)$ admits a conjugation action 
by the simplicial group $UR^{\times}$ and this action is equivalent 
to the one in question (here $k\rscript{c}$ is a cofibrant replacement of $k$ 
as a Gamma-ring). But in this model for the conjugation action, 
the composite of the approximation map $k\rscript{c}\to k$ 
with the unit map $k\to R_3$ of the $k$--algebra structure on $R_3$ 
is a point-set level fixed point of the conjugation action of $UR^{\times}$.
\end{proof}

If we combine the previous result with Theorem \ref{main} 
we can deduce that the Gamma-ring $DB$ is not stably equivalent 
to the Eilenberg--MacLane Gamma-ring of any simplicial ring, 
unless $B$ is an algebra over the rational numbers. 
This should be compared to Theorem \ref{DB summary} (b)
which says that as a $\Gamma$--space, 
$DB$ is stably equivalent to the smash product $H\mZ\sm^L HB$. 
In other words, $DB$ `additively' decomposes into a product of 
Eilenberg--MacLane $\Gamma$--spaces. In contrast the following corollary shows 
that the multiplicative structure of $DB$ is genuinely homotopy-theoretic. 
The idea of the proof is that any stable equivalence between $DB$ 
and a simplicial ring would give $DB$ the structure of an $H\mZ$--algebra. 
Such an algebra structure in turn gives rise to a ``central'' Gamma-ring map 
$H\mZ\to DB$ (in the sense of Theorem \ref{algebra structures}). 
But Theorem \ref{main} identifies all Gamma-ring map from $H\mZ$ to $DB$ 
and shows that none of them is ``central''.

\begin{corollary} \label{DB is not simplicial} 
For a commutative ring $B$ the following conditions are equivalent.
\renewcommand{\labelenumi}{\rm(\arabic{enumi})}
\begin{enumerate} 
\item $B$ is an algebra over the rational numbers.
\item The projection $DB\to D_1B\iso HB$ 
is a stable equivalence of Gamma-rings.
\item $DB$ is stably equivalent (as a Gamma-ring) 
to the Eilenberg--MacLane Gamma-ring of a simplicial ring.
\end{enumerate}
\end{corollary}
\begin{proof} Condition (1) is equivalent to condition (2) 
by Theorem \ref{DB summary} (b), and condition (2) 
implies condition (3). The proof that condition (3) implies condition (2) 
is a combination of Theorem \ref{algebra structures}, in the case $k=H\mZ$, 
with Theorem \ref{main}. Assume that condition (3) holds. 
Since the Eilenberg--MacLane Gamma-ring of a simplicial ring is 
an $H\mZ$--algebra, there exist a component of the space $\Rg(H\mZ,DB)$ 
of Gamma-ring maps for which the conjugation action map 
$DB^{\times} \to \Rg(H\mZ,DB)$ is homotopically trivial, 
by Theorem \ref{algebra structures}. 
But this map is part of a homotopy fiber sequence
\[ (DB)_1^{\times} \ \to  \ \Rg(H\mZ,DB) \ \to \ 
\Rg(H\mZ,DB)/\mbox{conj.} \]
Since the base of this fibration is weakly equivalent to the classifying space 
of a groupoid (Theorem \ref{main}), its homotopy groups are trivial above
dimension 1. So the conjugation action map is injective on 
homotopy groups in positive dimensions. Since the map is also null-homotopic, 
the space $(DB)_1^{\times}$ must be weakly contractible,
which implies condition (2).  
\end{proof}

\begin{numbered paragraph} \label{coordinate free} 
{\bf Coordinate free definition}\qua The definition of the Gamma-ring $DB$ 
and the Gamma-ring map $F_{\ast}$ depended on a formal group {\em law} $F$, 
ie, on a 1--dimensional commutative formal group with a choice of coordinate.
We will now describe coordinate-free versions of these constructions which 
at the same time are defined in a more general context. 

As input we consider a category $\C$ which has a zero object and finite 
coproducts. The natural enrichment of $\C$ over the category 
$\Gamma\rscript{op}$ of finite pointed sets is given by
 \[ X\, \sm \, k^+ \ = \ \underbrace{X\amalg \dots \amalg X}_k \qquad 
\mbox{(coproduct in $\C$).} \]
Every object $X$ of $\C$ has an {\em endomorphism Gamma-ring} 
\cite[4.6]{Sch:stable}, denoted $\End_{\C}(X)$ and defined by
\[ \End_{\C}(X)(k^+) \ = \ \Hom_{\C}(X,X\,\sm\, k^+) \ . \]
The unit map $\mS\to \End_{\C}(X)$ comes from the identity map of $X$,
viewed as a point in $\End_{\C}(X)(1^+)$, and the multiplication 
$\End_{\C}(X)\,\sm\,\End_{\C}(X)\to\End_{\C}(X)$ 
is induced by the composition product
\begin{eqnarray*} \End_{\C}(X)(k^+) \ \sm \ \End_{\C}(X)(l^+) & \to & 
\End_{\C}(X)(k^+\sm\, l^+) \ , \\
  f\, \sm \, g \hspace*{2.05cm} & \longmapsto & 
\quad (f\,\sm \,l^+) \, \circ \, g \ . 
\end{eqnarray*}
As an example we can take $\C$ to be the category of commutative, 
complete augmented $B$--algebras. If we choose the object $X$ to be the 
power series ring on one generator, then the endomorphism Gamma-ring of $X$ 
is precisely $DB$.

Now we suppose that the object $X$ of $\C$ is equipped with the structure 
of abelian cogroup object. So there is a given co-addition map 
$X\to X\amalg X$ and a co-inverse map $X\to X$ which make the set 
$\Hom_{\C}(X,Y)$ into an abelian group, natural for all objects $Y$ of $\C$. 
Every abelian cogroup structure on $X$ gives rise to a homomorphism of 
Gamma-rings $H\mZ\to\mbox{End}_{\C}(X)$ as follows. At a finite pointed set 
$k^+$ the map 
\[ H\mZ(k^+) \ = \ \widetilde\mZ[k^+] \ \to \ 
\Hom_{\C}(X,X\,\sm\, k^+) \ = \ \End_{\C}(X)(k^+) \]
is the additive extension of the map that sends $i\in k^+$ to the 
$i$-th coproduct inclusion $X\to X\,\sm\,k^+$. When $\C$ is the category of 
commutative, complete augmented $B$--algebras and $X$ is the power series 
ring on one generator, then making $X$ into an abelian cogroup object is 
the same thing as giving a (1--dimensional, commutative) formal group law 
$F$ over $B$. Furthermore, in this case the map $H\mZ\to \End_{\C}(X)$ 
arising from the abelian cogroup structure corresponds to the map $F_{\ast}$ 
of Construction \ref{Gamma-ring map} under the identification 
$DB\iso \End_{\C}(X)$. So from this point of view construction 
\ref{Gamma-ring map} is just a special case of the fact that every  
abelian cogroup structures gives rise to a homomorphism from the 
Eilenberg--MacLane Gamma-ring $H\mZ$ to an endomorphism Gamma-ring.
\end{numbered paragraph}

\begin{numbered paragraph} \label{non-commutative} 
{\bf Non-commutative formal group laws}\qua 
Construction \ref{Gamma-ring map} can be modified to work for not necessarily 
commutative formal group laws, but this variant does not lead 
to any interesting phenomena. There is a Gamma-ring, denoted by $Gp$, 
which is constructed the same way $H\mathbb Z$ is, but with free groups 
instead of free abelian groups. So as a $\Gamma$--space, 
$Gp$ takes a pointed set to the reduced free group it generates. 
The multiplication again comes from substitution, this times words in 
the generators of the free groups are substituted into each other. 
Abelianization gives a Gamma-ring map $Gp\to H\mZ$. If $F$ is a 1--dimensional 
but not necessarily commutative formal group law over the commutative ring 
$B$, then it gives rise to a Gamma-ring map
\[ F_{\ast}\,\co\, Gp \ \to \ DB \]
in much the same way as in Construction \ref{Gamma-ring map}. 
The Gamma-ring map $F_{\ast}$ factors over $H\mZ$ if and only if the 
formal group law $F$ is commutative. 

While the construction makes sense, $Gp$ is uninteresting as a source of 
Gamma-ring homomorphisms: we claim that the unit map $\mS \to Gp$ from 
the sphere Gamma-ring is a stable equivalence. This claim follows 
from the fact that the map from a high dimensional sphere into the free group 
it generates is an equivalence in the stable range. Since the Gamma-ring $Gp$ 
is stably equivalent to the initial Gamma-ring, the derived
space of homomorphisms into any other Gamma-ring is contractible.
\end{numbered paragraph}

\begin{numbered paragraph} \label{higher dimensional}
{\bf Higher dimensional formal group laws}\qua 
Another variant of Construction \ref{Gamma-ring map} proceeds from 
an $n$--dimensional commutative formal group law $F$. 
This time the construction gives a weak Gamma-ring map 
\[ F_{\ast} \,\co\,H\mZ \ \to \ M_n(DB) \]
into the Gamma-ring of $n\times n$--matrices over $DB$. 
For an arbitrary Gamma-ring $R$ the Gamma-ring $M_n(R)$ of 
$n\times n$--matrices over $R$ is defined as the endomorphism Gamma-ring 
of the free $R$--module on $n$ generators, ie,
\[ M_n(R) \ = \ \Hom_{R\text{-mod}}(R\,\sm\, n^+, R\,\sm\, n^+) \]
(here $\Hom_R$ refers to the internal homomorphism $\Gamma$--space in the 
category of $R$--modules). We define another Gamma-ring 
$\bar{M}_n(DB)$ as the endomorphism Gamma-ring, 
in the sense of \ref{coordinate free},
of the power series ring in $n$ generators 
in the category of augmented, complete $B$--algebras. 
Then there is a stable equivalence $M_n(DB)\to\bar{M}_n(DB)$. 
Since an $n$--dimensional commutative formal group law $F$ is the same thing 
as an abelian cogroup structure on the power series ring in $n$ variables, 
it leads to a map of Gamma-rings $F_{\ast}\co H\mZ\to \bar{M}_n(DB)$.
\end{numbered paragraph}

\begin{numbered paragraph} \label{formal modules} 
{\bf Formal module structures}\qua 
Yet another variation of our main theme consists in considering formal module 
structures over an associative ring $R$. A {\em (1--dimensional) formal 
$R$--module (law)} over a commutative ring $B$ consists of a 
(1--dimensional and commutative) formal group law $F$ and a ring homomorphism 
from $R$ into the endomorphism ring of the formal group law $F$. 
In the spirit of Construction \ref{Gamma-ring map}, any formal $R$--module 
structure $F$ over $B$ gives rise to a Gamma-ring map $F_{\ast}\co HR\to DB$ 
with source the Eilenberg--MacLane Gamma-ring of $R$. 

By the same method as in Section \ref{comparison map} we obtain a map
\[ \kappa_R \,\co\,\widetilde{\mathcal{F}R\text{-mod}}\rscript{str}(B)
\ \to\ \Rg(HR,DB)/\mbox{conj.} \]
from the groupoid of formal $R$--module structures over $B$ and strict
isomorphisms to the homotopy orbits of the derived space of Gamma-ring maps
by the connected component of the homotopy units $DB^{\times}$.
It seems reasonable to expect that the map $\kappa_R$ is
again a weak equivalence; whether this is the case depends on whether
the appropriate analog of Theorem \ref{finite functors} holds
over the ring $R$.

We denote by $R\rscript{com}= R/(rs-sr)$ the quotient
of $R$ by the commutator ideal. Then the arguments of Sections
\ref{filtration} through \ref{der and hyper} can be adapted 
to show:

\begin{theorem} \label{generalization}
Suppose that the ring $R$ has the following property:
for all $m\in\mZ$, all $k\geq 0$ and $R\rscript{\em com}$--modules $A$
the map
\[ \Ext^m_{\F(R\rscript{\em com})}(I,A\tensor_{R^{\text{\em com}}} {\mathcal S}^k) \ \to \  
\hExt^m_{\F(R\rscript{\em com})}(I,Q(A\tensor_{R^{\text{\em com}}} {\mathcal S}^k)) \]
is an isomorphism. Then the  map 
\[  \kappa_R \,\co\,\widetilde{\mathcal{F}R\text{\em -mod}}\rscript{\em str}(B)
\ \to\ \Rg(HR,DB)/\mbox{\em conj.} \] 
is a weak equivalence of simplicial sets.
\end{theorem}

The hypothesis is true for $R=\mZ$ --- this is the content of
Theorem \ref{finite functors}. The proof of Theorem \ref{finite functors}
can be adapted to establish the hypothesis for $R=\mZ/n$.
We conjecture that indeed the hypothesis holds  in general;
we expect that a `good' proof of Theorem \ref{finite functors}, 
ie, a proof that does not use the calculations of the
MacLane cohomology groups $\Ext^*_{\F}(I,{\mathcal S}^k)$ as input,
would also work in the more general context.
If this is the case, the map $\kappa_R$ is a weak equivalence 
for any ring $R$. 
\end{numbered paragraph}

\begin{example} Suppose $B$ is an $\mF_p$--algebra and $F$ a formal
group law over $F$. In this case we can reinterpret the height
of $F$ in terms of the homotopy class of the Gamma-ring map $F_*$.

The $p$--series of $F$ is either trivial or of the form
\[ [p]_F(x) \ = \ u\cdot x^{p^h} \ + \ \mbox{terms of higher degree} \]
for some $h\geq 1$ and some non-zero $u\in B$. 
The number $h$ is called the {\em height} of $F$. If $[p]_F=0$, then
$F$ is isomorphic to the additive formal group law 
\cite[III.1 Cor.\ 2]{Froehlich}, and the height of $F$ is infinite.

\medskip

{\bf Claim}\qua The height of $F$ is equal to the largest number $h$ such that
\[ H\mZ \ \varr{1cm}{F_*} \ DB \ \to \ D_{p^h-1}B  \] 
can be factored, in the homotopy category of Gamma-rings, 
over the Eilenberg--MacLane Gamma-ring for $\mF_p$.

\medskip

Indeed, if $F$ has height $h$, then for every pointed set $K$ the map
\[ F_{\ast}(K) \,\co\, \widetilde\mZ[K] = H\mZ(K) \ \to \ 
 DB(K) \ \subseteq \ \widetilde{B}[\![K]\!] \]
satisfies
\[ \qquad F_*(p\cdot x) \ = \ [p]_F(F_*(x)) \ \equiv \ 0 
\mbox{\qquad modulo degree $p^h$.}\]
So the composite map
\[ H\mZ \ \varr{1cm}{F_*} \ DB \ \to \ D_{p^h-1}B  \] 
factors uniquely over $H\mF_p$ on the point-set level.
Conversely, suppose that there exists a commutative square 
in the homotopy category of Gamma-rings
\[\xymatrix@C=15mm{ H\mZ \ar^-{F_*}[r] \ar[d] & DB \ar[d] \\
H\mF_p \ar@{-->}[r] & D_{p^h-1}B }\]
By the analog of Theorem \ref{k-bud comparison} 
for buds of formal $\mF_p$--modules (which holds since the hypothesis
of Theorem \ref{generalization} are satisfied for $R=\mF_p$), the maps
from $H\mF_p \to D_{p^h-1}B$ in the homotopy category of Gamma-rings
are in bijective correspondence with strict isomorphism classes of
$p^h$--buds of formal $\mF_p$--module structures on $B$.
But over an $\mF_p$--algebra every formal $\mF_p$--module bud
is strictly isomorphic to the additive formal
group law \cite[III.1 Cor.\ 2]{Froehlich}. 
Hence the $(p^h-1)$--buds of $F$ and of the additive formal group law over $B$ 
induce the same maps $H\mZ\to D_{p^h-1}B$ 
in  the homotopy category of Gamma-rings.
By Theorem \ref{k-bud comparison}, the $(p^h-1)$--bud of $F$ is thus strictly
isomorphic to the $(p^h-1)$--bud of the additive formal group law, and so
$[p]_F\equiv 0$ modulo degree $p^h-1$. So the height of $F$ is at least $h$.
\end{example}


\begin{thebibliography}

\bibitem{Basterra:AQHomology}
\textbf{M Basterra}, \emph{Andr\'e-{Q}uillen cohomology of commutative
  ${S}$--algebras}, J. Pure Appl. Algebra 144 (1999) 111--143

\bibitem{Basterra-McCarthy}
\textbf{M Basterra}, \textbf{R McCarthy}, \emph{{$\Gamma$}-homology,
  topological {A}ndr\'e-{Q}uillen homology and stabilization}, Topology Appl.
  121 (2002) 551--566

\bibitem{Basterra-Richter}
\textbf{M Basterra}, \textbf{B Richter}, \emph{({C}o-)homology theories for
  commutative \mbox{($S$-)}alg\-eb\-ras}, from: ``Structured Ring Spectra'', (A
  Baker, B Richter, editors), to appear in London Mathematical Society Lecture
  Notes, Cambridge University Press

\bibitem{Betley-Pira:stableK}
\textbf{S Betley}, \textbf{T Pirashvili}, \emph{Stable ${K}$--theory as a
  derived functor}, J. Pure Appl. Algebra 96 (1994) 245--258

\bibitem{Boek:THH}
\textbf{M B{\"o}kstedt}, \emph{Topological {Hochschild} homology} (1985),
  preprint, Bielefeld

\bibitem{Borceux:2}
\textbf{F Borceux}, \emph{Handbook of Categorical Algebra 2: Categories and
  Structures}, volume~51 of \emph{Encyclopedia of Mathematics and its
  Applications}, Cambridge University Press (1994)

\bibitem{Bous:operations}
\textbf{A\,K Bousfield}, \emph{Operations on Derived Functors of Non-additive
  Functors} (1967), mimeographed Notes, Brandeis University

\bibitem{BF}
\textbf{A\,K Bousfield}, \textbf{E\,M Friedlander}, \emph{Homotopy theory of
  ${\Gamma}$--spaces, spectra, and bisimplicial sets}, from: ``Geometric
  Applications of Homotopy Theory II'', {L}ecture Notes in Mathematics 658,
  Springer (1978)  80--130

\bibitem{Bous-Kan}
\textbf{A\,K Bousfield}, \textbf{D\,M Kan}, \emph{Homotopy limits, completions
  and localizations}, volume 304 of \emph{{L}ecture Notes in Mathematics},
  Springer (1972)

\bibitem{Cartan}
\textbf{H Cartan}, \emph{Alg{\'e}bres d'{E}ilenberg--{M}ac{L}ane et homotopie}
  (1954-55), seminaire Henri Cartan

\bibitem{Dold-Puppe}
\textbf{A Dold}, \textbf{D Puppe}, \emph{Homologie nicht-additiver {F}unktoren.
  {A}nwendungen}, Ann. Inst. Fourier Grenoble 11 (1961) 201--312

\bibitem{Dwyer:oper}
\textbf{W\,G Dwyer}, \emph{Homotopy operations for simplicial commutative
  algebras}, Trans. Amer. Math. Soc. 260 (1980) 421--435

\bibitem{Dwyer-Spalinski}
\textbf{W\,G Dwyer}, \textbf{J Spalinski}, \emph{Homotopy theories and model
  categories}, from: ``Handbook of algebraic topology'', (I\,M James, editor),
  North-Holland, Amsterdam (1995)  73--126

\bibitem{EM:cubical}
\textbf{S Eilenberg}, \textbf{S MacLane}, \emph{Homology theories for
  multiplicative systems}, Trans. Amer. Math. Soc. 71 (1951) 294--330

\bibitem{EM:K(pi)-II}
\textbf{S Eilenberg}, \textbf{S MacLane}, \emph{On the groups {$H(\pi,n)$},
  {II}}, Ann. Math. 60 (1954) 49--139

\bibitem{Franjou:puissance}
\textbf{V Franjou}, \emph{Extensions entre puissances ext\'erieures et entre
  puissances sym\'et\-riques}, J. Algebra 179 (1996) 501--522

\bibitem{Franjou-Lannes-Schwartz:MacLane(F_p)}
\textbf{V Franjou}, \textbf{J Lannes}, \textbf{L Schwartz}, \emph{Autour de la
  cohomologie de {M}ac {L}ane des corps finis}, Invent. Math. 115 (1994)
  513--538

\bibitem{Franjou-Pira:MacLane(Z)}
\textbf{V Franjou}, \textbf{T Pirashvili}, \emph{On the {M}ac{L}ane cohomology
  for the ring of integers}, Topology 37 (1998) 109--114

\bibitem{Froehlich}
\textbf{A Fr{\"o}hlich}, \emph{Formal groups}, volume~74 of \emph{{L}ecture
  Notes in Mathematics}, Springer-Verlag, Berlin (1968)

\bibitem{Goerss:Hilton-Milnor}
\textbf{P\,G Goerss}, \emph{A {Hilton-Milnor} theorem for categories of
  simplicial algebras}, Amer. J. Math. 111 (1989) 927--971

\bibitem{Goerss:AQ-homology}
\textbf{P\,G Goerss}, \emph{On the {Andr\'e--Quillen} cohomology of commutative
  {${\mathbb F}_2$}-algebras}, Asterisque 186 (1990)

\bibitem{Goerss-Jardine:book}
\textbf{P\,G Goerss}, \textbf{J\,F Jardine}, \emph{Simplicial homotopy theory},
  Birkh\"auser Verlag, Basel (1999)

\bibitem{Hovey:book}
\textbf{M Hovey}, \emph{Model categories}, American Mathematical Society,
  Providence, RI (1999)

\bibitem{HSS}
\textbf{M Hovey}, \textbf{B\,E Shipley}, \textbf{J\,H Smith}, \emph{Symmetric
  spectra}, J. Amer. Math. Soc. 13 (2000) 149--208

\bibitem{Illusie:cotangentI}
\textbf{L Illusie}, \emph{Complexe cotangent et d\'eformations. {I}}, volume
  239 of \emph{{L}ecture Notes in Mathematics}, Springer-Verlag, Berlin (1971)

\bibitem{JP}
\textbf{M Jibladze}, \textbf{T Pirashvili}, \emph{Cohomology of algebraic
  theories}, J. Algebra 137 (1991) 253--296

\bibitem{Johnson-McCarthy}
\textbf{B Johnson}, \textbf{R McCarthy}, \emph{Linearization, {D}old-{P}uppe
  stabilization, and {M}ac {L}ane's ${Q}$--construction}, Trans. Amer. Math.
  Soc. 350 (1998) 1555--1593

\bibitem{Lazard:formels}
\textbf{M Lazard}, \emph{Sur les groupes de {L}ie formels \`a un param\`etre},
  Bull. Soc. Math. France 83 (1955) 251--274

\bibitem{Lydakis:Smash}
\textbf{M Lydakis}, \emph{Smash products and {$\Gamma$}-spaces}, Math. Proc.
  Cambridge Philos. Soc. 126 (1999) 311--328

\bibitem{MacL:ML-homology}
\textbf{S MacLane}, \emph{Homologie des anneaux et des modules}, from:
  ``Colloque de topologie alg\'ebrique, Louvain, 1956'', Georges Thone, Li\`ege
  (1957)  55--80

\bibitem{MacL:categories}
\textbf{S MacLane}, \emph{Categories for the working mathematician}, Springer
  (1971)

\bibitem{MMSS}
\textbf{M\,A Mandell}, \textbf{J\,P May}, \textbf{S Schwede}, \textbf{B
  Shipley}, \emph{Model categories of diagram spectra}, Proc. London Math. Soc.
  (3) 82 (2001) 441--512

\bibitem{McDuff:monoids}
\textbf{D McDuff}, \emph{On the classifying spaces of discrete monoids},
  Topology 18 (1979) 313--320

\bibitem{Pira:additivization}
\textbf{T Pirashvili}, \emph{Higher additivizations}, Trudy Tbiliss. Mat. Inst.
  Razmadze Akad. Nauk Gruzin. SSR 91 (1988) 44--54

\bibitem{Pira:approximation}
\textbf{T Pirashvili}, \emph{Polynomial approximation of {E}xt and {T}or groups
  in functor categories}, Comm. Algebra 21 (1993) 1705--1719

\bibitem{Pira-Richter}
\textbf{T Pirashvili}, \textbf{B Richter}, \emph{Robinson-{W}hitehouse complex
  and stable homotopy}, Topology 39 (2000) 525--530

\bibitem{Pira-Wald}
\textbf{T Pirashvili}, \textbf{F Waldhausen}, \emph{Mac{L}ane homology and
  topological {H}ochschild homology}, J. Pure Appl. Algebra 82 (1992) 81--98

\bibitem{Quillen:comrings}
\textbf{D Quillen}, \emph{On the (co-) homology of commutative rings}, from:
  ``Applications of Categorical Algebra (Proc. Sympos. Pure Math., Vol. XVII,
  New York, 1968)'', Amer. Math. Soc., Providence, R.I. (1970)  65--87

\bibitem{Quillen:HA}
\textbf{D\,G Quillen}, \emph{Homotopical Algebra}, volume~43 of \emph{{L}ecture
  Notes in Mathematics}, Springer (1967)

\bibitem{Reedy:simplicial}
\textbf{C\,L Reedy}, \emph{Homotopy theory of model categories}, preprint,
  1973.

\bibitem{Robinson:Gamma-Einfty}
\textbf{A Robinson}, \emph{Gamma homology, {L}ie representations and {$E\sb
  \infty$} multiplications}, Invent. Math. 152 (2003) 331--348

\bibitem{Robinson-Whitehouse:GammaHomology}
\textbf{A Robinson}, \textbf{S Whitehouse}, \emph{Operads and
  {$\Gamma$}-homology of commutative rings}, Math. Proc. Cambridge Philos. Soc.
  132 (2002) 197--234

\bibitem{Sch:modelspec}
\textbf{S Schwede}, \emph{Spectra in model categories and applications to the
  algebraic cotangent complex}, J. Pure Appl. Algebra 120 (1997) 77--104

\bibitem{Sch:Gamma}
\textbf{S Schwede}, \emph{Stable homotopical algebra and {$\Gamma$}-spaces},
  Math. Proc. Cambridge Philos. Soc. 126 (1999) 329--356

\bibitem{Sch:comparison}
\textbf{S Schwede}, \emph{${S}$--modules and symmetric spectra}, Math. Ann. 319
  (2001) 517--532

\bibitem{Sch:stable}
\textbf{S Schwede}, \emph{Stable homotopy of algebraic theories}, Topology 40
  (2001) 1--41

\bibitem{SS:monoidal}
\textbf{S Schwede}, \textbf{B\,E Shipley}, \emph{Algebras and modules in
  monoidal model categories}, Proc. London Math. Soc. (3) 80 (2000) 491--511

\bibitem{Segal:categories}
\textbf{G Segal}, \emph{Categories and cohomology theories}, Topology 13 (1974)
  293--312

\bibitem{Shipley:THH}
\textbf{B Shipley}, \emph{Symmetric spectra and topological {H}ochschild
  homology}, $K$--Theory 19 (2000) 155--183

\bibitem{simson}
\textbf{D Simson}, \emph{Stable derived functors of the second symmetric power
  functor, second exterior power functor and {W}hitehead gamma functor},
  Colloq. Math. 32 (1974) 49--55

\bibitem{Turner:looking_glass}
\textbf{J\,M Turner}, \emph{Simplicial commutative $\mathbf{F}_p$--algebras
  through the looking-glass of $\mathbf{F}_p$--local spaces}, from: ``Homotopy
  invariant algebraic structures (Baltimore, MD, 1998)'', Amer. Math. Soc.,
  Providence, RI (1999)  353--363

\end{thebibliography}
\end{document}